\documentstyle[12pt,fleqn]{article}

\begin{document}

\title{Commutative Complex Numbers in Four Dimensions}

\author{Silviu Olariu
\thanks{e-mail: olariu@ifin.nipne.ro}\\
Institute of Physics and Nuclear Engineering, Tandem Laboratory\\
76900 Magurele, P.O. Box MG-6, Bucharest, Romania}

\date{4 August 2000}

\maketitle

\abstract

Commutative complex numbers of the form $u=x+\alpha y+\beta z+\gamma t$ in 4
dimensions are studied, the variables $x, y, z$ and $t$ being real numbers.
Four distinct types of multiplication rules for the complex bases $\alpha,
\beta$ and $\gamma$ are investigated, which correspond to hypercomplex entities
called in this paper circular fourcomplex numbers, hyperbolic fourcomplex
numbers, planar fourcomplex numbers, and polar fourcompex numbers.  Exponential
and trigonometric forms for the fourcomplex numbers are given in all these
cases.  Expressions are given for the elementary functions of the fourcomplex
variables mentioned above.  Relations of equality exist between the partial
derivatives of the real components of the functions of fourcomplex variables.
The integral of a fourcomplex function between two points is independent of the
path connecting the points.  The concepts of poles and residues can be
introduced for the circular, planar, and polar fourcomplex numbers, for which
the exponential forms depend on cyclic variables.  A hypercomplex polynomial
can be written as a product of linear factors for circular and planar
fourcomplex numbers, and as a product of linear or quadratic factors for the
hyperbolic and polar fourcomplex numbers.

\endabstract

\section{Introduction}

A regular, two-dimensional complex number $x+iy$ 
can be represented geometrically by the modulus $\rho=(x^2+y^2)^{1/2}$ and 
by the polar angle $\theta=\arctan(y/x)$. The modulus $\rho$ is multiplicative
and the polar angle $\theta$ is additive upon the multiplication of ordinary 
complex numbers.

The quaternions of Hamilton are a system of hypercomplex numbers
defined in four dimensions, the
multiplication being a noncommutative operation, \cite{1} 
and many other hypercomplex systems are
possible, \cite{2a}-\cite{2b} but these hypercomplex systems 
do not have all the required properties of regular, 
two-dimensional complex numbers which rendered possible the development of the 
theory of functions of a complex variable.

This paper belongs to a series of studies on commutative complex numbers in $n$
dimensions. \cite{2c}
Systems of hypercomplex numbers in 4 dimensions of the form 
$u=x+\alpha y+\beta z+\gamma t$ are described in this work,
where the variables $x, y,z$ and $t$ are real numbers,
for which the multiplication is both associative and commutative.
The product of two
fourcomplex numbers is equal to zero if both numbers are 
equal to zero, or if the numbers belong to certain four-dimensional hyperplanes
as discussed further in this work. 
The fourcomplex numbers are
rich enough in properties such that an exponential and trigonometric forms
exist 
and the concepts of analytic fourcomplex 
function,  contour integration and residue can be defined. 
Expressions are given for the elementary functions of fourcomplex variable.
The functions $f(u)$ of fourcomplex variable defined by power series, 
have derivatives $\lim_{u\rightarrow u_0} [f(u)-f(u_0)]/(u-u_0)$ independent of
the direction of approach of $u$ to $u_0$. If the fourcomplex function $f(u)$
of the fourcomplex variable $u$ is written in terms of 
the real functions $P(x,y,z,t),Q(x,y,z,t),R(x,y,z,t),S(x,y,z,t)$, then
relations of equality  
exist between partial derivatives of the functions $P,Q,R,S$. 
The integral $\int_A^B f(u) du$ of a fourcomplex
function between two points $A,B$ is independent of the path connecting $A,B$.

Four distinct types
of hypercomplex numbers are studied, as discussed further.
In Sec. II, the multiplication rules for the complex units $\alpha, \beta$ and
$\gamma$ are 
$\alpha^2=-1, \:\beta^2=-1, \:\gamma^2=1, \alpha\beta=\beta\alpha=-\gamma,\:
\alpha\gamma=\gamma\alpha=\beta, \:\beta\gamma=\gamma\beta=\alpha$.
The exponential form of a fourcomplex number is
$x+\alpha y+\beta z+\gamma t
=\rho\exp\left[\gamma \ln\tan\psi+
(1/2)\alpha(\phi+\chi)+(1/2)\beta(\phi-\chi) 
\right]$ , where the amplitude is
$\rho^4=\left[(x+t)^2+(y+z)^2\right]\left[(x-t)^2+(y-z)^2\right]$ ,
$\phi, \chi$ are azimuthal angles, $0\leq\phi<2\pi, 0\leq\chi<2\pi$, 
and $\psi$ is a planar angle, $0<\psi\leq\pi/2$.
The trigonometric form of a fourcomplex number is
$x+\alpha y+\beta z+\gamma t=
d[\cos(\psi-\pi/4)+\gamma\sin(\psi-\pi/4)]$
$\exp\left[(1/2)\alpha(\phi+\chi)+(1/2)\beta(\phi-\chi)\right]$, 
where $d^2=x^2+y^2+z^2+t^2$.
The amplitude $\rho$ and $\tan\psi$ are multiplicative 
and the angles $\phi, \chi$ 
are additive upon the multiplication of fourcomplex numbers.
Since there are two cyclic variables, $\phi$ and $\chi$, these fourcomplex
numbers are called circular fourcomplex numbers.
If $f(u)$ is an analytic fourcomplex function, then 
$\oint_\Gamma f(u)du/(u-u_0)=
\pi[(\alpha+\beta) \:\;{\rm int}(u_{0\xi\upsilon},\Gamma_{\xi\upsilon})
+ (\alpha-\beta)\:\;{\rm int}(u_{0\tau\zeta},\Gamma_{\tau\zeta})]\;f(u_0)$,
where the functional ${\rm int}$ takes the values 0 or 1 depending on the
relation between the projections of the point $u_0$ and of the curve $\Gamma$
on certain planes.
A polynomial  can be written as a 
product of linear or quadratic factors, although the factorization may not be
unique.

In Sec. III, the multiplication rules for the complex units $\alpha, \beta$ and
$\gamma$ are
$\alpha^2=1, \:\beta^2=1, \:\gamma^2=1, \alpha\beta=\beta\alpha=\gamma,\:
\alpha\gamma=\gamma\alpha=\beta, \:\beta\gamma=\gamma\beta=\alpha$. 
The exponential form of a fourcomplex number, which can be defined for
$s=x+y+z+t>0, \: s^\prime= x-y+z-t>0 , \: s^{\prime\prime}=x+y-z-t>0,  \:
s^{\prime\prime\prime}=x-y-z+t>0$, is
$x+\alpha y+\beta z+\gamma t=
\mu\exp(\alpha y_1+\beta z_1+\gamma t_1)$, where the amplitude is
$\mu=(s s^\prime s^{\prime\prime}s^{\prime\prime\prime})^{1/4}$ and the
arguments are
$y_1=(1/4)\ln(ss^{\prime\prime}/s^\prime s^{\prime\prime\prime}),\:
z_1=(1/4)\ln(ss^{\prime\prime\prime}/s^\prime s^{\prime\prime}),\:
t_1=(1/4)\ln(ss^\prime/s^{\prime\prime}s^{\prime\prime\prime})$.
Since there is no cyclic variable, these fourcomplex numbers are called
hyperbolic fourcomplex numbers. 
The amplitude $\mu$ is multiplicative and the arguments $y_1, z_1, t_1$
are additive upon the multiplication of fourcomplex numbers.
A polynomial can be
written as a 
product of linear or quadratic factors, although the factorization may not be
unique. 

In Sec. IV, the multiplication rules for the complex units $\alpha, \beta$ and
$\gamma$ are
$\alpha^2=\beta, \:\beta^2=-1, \:\gamma^2=-\beta,
\alpha\beta=\beta\alpha=\gamma,\: 
\alpha\gamma=\gamma\alpha=-1, \:\beta\gamma=\gamma\beta=-\alpha$ .
The exponential function of a fourcomplex number can be expanded in terms of
the four-dimensional cosexponential functions
$f_{40}(x)=1-x^4/4!+x^8/8!-\cdots ,\:\:
f_{41}(x)=x-x^5/5!+x^9/9!-\cdots,\:\:
f_{42}(x)=x^2/2!-x^6/6!+x^{10}/10!-\cdots ,\:\:
f_{43}(x)=x^3/3!-x^7/7!+x^{11}/11!-\cdots$. Expressions are obtained for the
four-dimensional cosexponential functions in terms of elementary functions.
The exponential form of a fourcomplex number is
$x+\alpha y+\beta z+\gamma t
=\rho\exp\left\{(1/2)(\alpha-\gamma)\ln\tan\psi
+(1/2)[\beta+(\alpha+\gamma)/\sqrt{2}]\phi
-(1/2)[\beta-(\alpha+\gamma)/\sqrt{2}]\chi
\right\}$,
where the amplitude is $\rho^4=\left\{\left[x+(y-t)/\sqrt{2}\right]^2
+\left[z+(y+t)/\sqrt{2}\right]^2\right\}
\left\{
\left[x-(y-t)/\sqrt{2}\right]^2+\left[z-(y+t)/\sqrt{2}\right]^2\right\}$,  
$\phi, \chi$ are azimuthal angles, $0\leq\phi<2\pi, 0\leq\chi<2\pi$, 
and $\psi$ is a planar angle, $0\leq\psi\leq\pi/2$.
The trigonometric form of a fourcomplex number is
$x+\alpha y+\beta z+\gamma t
=d\left[\cos\left(\psi-\pi/4\right)\right.$
$\left.+(1/\sqrt{2})(\alpha-\gamma)\sin\left(\psi-\pi/4\right)\right]$
$\exp\left\{
(1/2)[\beta+(\alpha+\gamma)/\sqrt{2}]\phi
-(1/2)[\beta-(\alpha+\gamma)/\sqrt{2}]\chi\right\}$,
where $d^2=x^2+y^2+z^2+t^2$.
The amplitude $\rho$ and $\tan\psi$ are multiplicative 
and the angles $\phi, \chi$ 
are additive upon the multiplication of fourcomplex numbers.
There are two cyclic variables, $\phi$ and $\chi$, so that these fourcomplex
numbers are also of a circular type. In order to distinguish
them from the circular hypercomplex numbers, these are called
planar fourcomplex numbers.
If $f(u)$ is an analytic fourcomplex function, then 
$\oint_\Gamma f(u)du/(u-u_0)=
\pi\left[\left(\beta+(\alpha+\gamma)/\sqrt{2}\right) \;{\rm
int}(u_{0\xi\upsilon},\Gamma_{\xi\upsilon}) \right.$
$\left.+ \left(\beta-(\alpha+\gamma)/\sqrt{2}\right)\;{\rm
int}(u_{0\tau\zeta},\Gamma_{\tau\zeta})\right]\;f(u_0)$ , 
where the functional ${\rm int}$ takes the values 0 or 1 depending on the
relation between the projections of the point $u_0$ and of the curve $\Gamma$
on certain planes.
A polynomial can be
written as a 
product of linear or quadratic factors, although the factorization may not be
unique. 
The fourcomplex numbers described in this section are a particular case for 
$n=4$ of the planar hypercomplex numbers in $n$ dimensions.\cite{2c},\cite{2d}

In Sec. V, the multiplication rules for the complex units $\alpha, \beta$ and
$\gamma$ are
$\alpha^2=\beta, \:\beta^2=1, \:\gamma^2=\beta,
\alpha\beta=\beta\alpha=\gamma,\: 
\alpha\gamma=\gamma\alpha=1, \:\beta\gamma=\gamma\beta=\alpha$.
The product of two fourcomplex numbers is equal to zero if both numbers are
equal to zero, or if the numbers belong to certain four-dimensional hyperplanes
described further in this work. 
The exponential function of a fourcomplex number can be expanded in terms of
the four-dimensional cosexponential functions
$g_{40}(x)=1+x^4/4!+x^8/8!+\cdots ,\:\:
g_{41}(x)=x+x^5/5!+x^9/9!+\cdots,\:\:
g_{42}(x)=x^2/2!+x^6/6!+x^{10}/10!+\cdots ,\:\:
g_{43}(x)=x^3/3!+x^7/7!+x^{11}/11!+\cdots$. Addition theorems and other
relations are obtained for these four-dimensional cosexponential functions.
The exponential form of a fourcomplex number, which can be defined for
$x+y+z+t>0, x-y+z-t>0$, is
$u=\rho\exp\left[(1/4)(\alpha+\beta+\gamma) \ln(\sqrt{2}/\tan\theta_+)
-(1/4)(\alpha-\beta+\gamma) \ln(\sqrt{2}/\tan\theta_-)
+(\alpha-\gamma)\phi/2\right]$,
where $\rho=(\mu_+\mu_-)^{1/2},\:\: \mu_+^2=(x-z)^2+(y-t)^2,\:\:
\mu_-^2=(x+z)^2-(y+t)^2$, $e_+=(1+\alpha+\beta+\gamma)/4,
e_-=(1-\alpha+\beta-\gamma)/4$,
$e_1=(1-\beta)/2, \tilde e_1=(\alpha-\gamma)/2$,
the polar angles are 
$\tan\theta_+=\sqrt{2}\mu_+/v_+,  
\tan\theta_-=\sqrt{2}\mu_+/v_-$, $0\leq\theta_+\leq\pi, 
0\leq\theta_-\leq\pi$,
and the azimuthal angle $\phi$ is defined by the relations 
$x-y=\mu_+\cos\phi,\:\:z-t=\mu_+\sin\phi$, $0\leq\phi<2\pi$. 
The trigonometric form of the fourcomplex number $u$ is
$u=d\sqrt{2}$
$\left(1+1/\tan^2\theta_++1/\tan^2\theta_-\right)^{-1/2}
\left\{e_1+e_+\sqrt{2}/\tan\theta_+
+e_-\sqrt{2}/\tan\theta_-\right\}$
$\exp\left[\tilde e_1\phi\right]$.
The amplitude $\rho$ and $\tan\theta_+/\sqrt{2}, \tan\theta_-/\sqrt{2}$,
are multiplicative, and the azimuthal angle $\phi$ is additive
upon the multiplication of fourcomplex numbers.
There is only one cyclic variable, $\phi$, 
and there are two axes $v_+, v_-$ which play an important role in the
description of these numbers, so that
these hypercomplex numbers are called polar fourcomplex numbers.
If $f(u)$ is an analytic fourcomplex function, then 
$\oint_\Gamma f(u)du/(u-u_0)=
\pi(\beta-\gamma) \;{\rm int}(u_{0\xi\upsilon},\Gamma_{\xi\upsilon})f(u_0)$ ,
where the functional ${\rm int}$ takes the values 0 or 1 depending on the
relation between the projections of the point $u_0$ and of the curve $\Gamma$
on certain planes.
A polynomial can be
written as a 
product of linear or quadratic factors, although the factorization may not be
unique.
The fourcomplex numbers described in this section are a particular case for 
$n=4$ of the polar hypercomplex numbers in $n$ dimensions.\cite{2c},\cite{2d}

\section{Circular Complex Numbers in Four Dimensions}

\subsection{Operations with circular fourcomplex numbers}

A circular fourcomplex number is determined by its four components $(x,y,z,t)$.
The sum of the circular fourcomplex numbers $(x,y,z,t)$ and
$(x^\prime,y^\prime,z^\prime,t^\prime)$ is the circular fourcomplex
number $(x+x^\prime,y+y^\prime,z+z^\prime,t+t^\prime)$. 
The product of the circular fourcomplex numbers
$(x,y,z,t)$ and $(x^\prime,y^\prime,z^\prime,t^\prime)$ 
is defined in this work to be the circular fourcomplex
number
$(xx^\prime-yy^\prime-zz^\prime+tt^\prime,
xy^\prime+yx^\prime+zt^\prime+tz^\prime,
xz^\prime+zx^\prime+yt^\prime+ty^\prime,
xt^\prime+tx^\prime-yz^\prime-zy^\prime)$.

Circular fourcomplex numbers and their operations can be represented by
writing the 
circular fourcomplex number $(x,y,z,t)$ as  
$u=x+\alpha y+\beta z+\gamma t$, where $\alpha, \beta$ and $\gamma$ 
are bases for which the multiplication rules are 
\begin{equation}
\alpha^2=-1, \:\beta^2=-1, \:\gamma^2=1, \alpha\beta=\beta\alpha=-\gamma,\:
\alpha\gamma=\gamma\alpha=\beta, \:\beta\gamma=\gamma\beta=\alpha .
\label{1}
\end{equation}
Two circular fourcomplex numbers $u=x+\alpha y+\beta z+\gamma t, 
u^\prime=x^\prime+\alpha y^\prime+\beta z^\prime+\gamma t^\prime$ are equal, 
$u=u^\prime$, if and only if $x=x^\prime, y=y^\prime,
z=z^\prime, t=t^\prime$. 
If 
$u=x+\alpha y+\beta z+\gamma t, 
u^\prime=x^\prime+\alpha y^\prime+\beta z^\prime+\gamma t^\prime$
are circular fourcomplex numbers, 
the sum $u+u^\prime$ and the 
product $uu^\prime$ defined above can be obtained by applying the usual
algebraic rules to the sum 
$(x+\alpha y+\beta z+\gamma t)+ 
(x^\prime+\alpha y^\prime+\beta z^\prime+\gamma t^\prime)$
and to the product 
$(x+\alpha y+\beta z+\gamma t)
(x^\prime+\alpha y^\prime+\beta z^\prime+\gamma t^\prime)$,
and grouping of the resulting terms,
\begin{equation}
u+u^\prime=x+x^\prime+\alpha(y+y^\prime)+\beta(z+z^\prime)+\gamma(t+t^\prime),
\label{1a}
\end{equation}
\begin{eqnarray}
\lefteqn{uu^\prime=
xx^\prime-yy^\prime-zz^\prime+tt^\prime+
\alpha(xy^\prime+yx^\prime+zt^\prime+tz^\prime)+
\beta(xz^\prime+zx^\prime+yt^\prime+ty^\prime)}\nonumber\\
&&+\gamma(xt^\prime+tx^\prime-yz^\prime-zy^\prime)
\label{1b}
\end{eqnarray}

If $u,u^\prime,u^{\prime\prime}$ are circular fourcomplex numbers, the
multiplication  
is associative
\begin{equation}
(uu^\prime)u^{\prime\prime}=u(u^\prime u^{\prime\prime})
\label{2}
\end{equation}
and commutative
\begin{equation}
u u^\prime=u^\prime u ,
\label{3}
\end{equation}
as can be checked through direct calculation.
The circular fourcomplex zero is $0+\alpha\cdot 0+\beta\cdot 0+\gamma\cdot 0,$ 
denoted simply 0, 
and the circular fourcomplex unity is $1+\alpha\cdot 0+\beta\cdot 0+\gamma\cdot
0,$  
denoted simply 1.

The inverse of the circular fourcomplex number 
$u=x+\alpha y+\beta z+\gamma t$ is a circular fourcomplex number
$u^\prime=x^\prime+\alpha y^\prime+\beta z^\prime+\gamma t^\prime$
having the property that
\begin{equation}
uu^\prime=1 .
\label{4}
\end{equation}
Written on components, the condition, Eq. (\ref{4}), is
\begin{equation}
\begin{array}{c}
xx^\prime-yy^\prime-zz^\prime+tt^\prime=1,\\
yx^\prime+xy^\prime+tz^\prime+zt^\prime=0,\\
zx^\prime+ty^\prime+xz^\prime+yt^\prime=0,\\
tx^\prime-zy^\prime-yz^\prime+xt^\prime=0 .
\end{array}
\label{5}
\end{equation}
The system (\ref{5}) has the solution
\begin{equation}
x^\prime=\frac{x(x^2+y^2+z^2-t^2)-2yzt}
{\rho^4} ,
\label{6a}
\end{equation}
\begin{equation}
y^\prime=
\frac{y(-x^2-y^2+z^2-t^2)+2xzt}
{\rho^4} ,
\label{6b}
\end{equation}
\begin{equation}
z^\prime=\frac{z(-x^2+y^2-z^2-t^2)+2xyt}
{\rho^4} ,
\label{6c}
\end{equation}
\begin{equation}
t^\prime=\frac{t(-x^2+y^2+z^2+t^2)-2xyz}
{\rho^4} ,
\label{6d}
\end{equation}
provided that $\rho\not=0, $ where
\begin{equation}
\rho^4=x^4+y^4+z^4+t^4+2(x^2y^2+x^2z^2-x^2t^2-y^2z^2+y^2t^2+z^2t^2)-8xyzt .
\label{6e}
\end{equation}
The quantity $\rho$ will be called amplitude of the circular fourcomplex number
$x+\alpha y+\beta z +\gamma t$.
Since
\begin{equation}
\rho^4=\rho_+^2\rho_-^2 ,
\label{7a}
\end{equation}
where
\begin{equation}
\rho_+^2=(x+t)^2+(y+z)^2, \: \rho_-^2=(x-t)^2+(y-z)^2 ,
\label{7b}
\end{equation}
a circular fourcomplex number $u=x+\alpha y+\beta z+\gamma t$ has an inverse,
unless 
\begin{equation}
x+t=0,\: y+z=0 ,  
\label{8}
\end{equation}
or
\begin{equation}
x-t=0,\: y-z=0 .  
\label{9}
\end{equation}

Because of conditions (\ref{8})-(\ref{9}) these 2-dimensional surfaces
will be called nodal hyperplanes. 
It can be shown that 
if $uu^\prime=0$ then either $u=0$, or $u^\prime=0$, 
or one of the circular fourcomplex numbers is of the form $x+\alpha y+\beta
y+\gamma x$ 
and the other of the form $x^\prime+\alpha y^\prime-\beta y^\prime -\gamma
x^\prime$.

\subsection{Geometric representation of circular fourcomplex numbers}

The circular fourcomplex number $x+\alpha y+\beta z+\gamma t$ can be
represented by  
the point $A$ of coordinates $(x,y,z,t)$. 
If $O$ is the origin of the four-dimensional space $x,y,z,t,$ the distance 
from $A$ to the origin $O$ can be taken as
\begin{equation}
d^2=x^2+y^2+z^2+t^2 .
\label{10}
\end{equation}
The distance $d$ will be called modulus of the circular fourcomplex number
$x+\alpha 
y+\beta z +\gamma t$, $d=|u|$. 
The orientation in the four-dimensional space of the line $OA$ can be specified
with the aid of three angles $\phi, \chi, \psi$
defined with respect to the rotated system of axes
\begin{equation}
\xi=\frac{x+t}{\sqrt{2}}, \: \tau=\frac{x-t}{\sqrt{2}}, \:
\upsilon=\frac{y+z}{\sqrt{2}}, \: \zeta=\frac{y-z}{\sqrt{2}} .
\label{11}
\end{equation}
The variables $\xi, \upsilon, \tau, \zeta$ will be called canonical 
circular fourcomplex variables.
The use of the rotated axes $\xi, \upsilon, \tau, \zeta$ 
for the definition of the angles $\phi, \chi, \psi$ 
is convenient for the expression of the circular fourcomplex numbers
in exponential and trigonometric forms, as it will be discussed further.
The angle $\phi$ is the angle between the projection of $A$ in the plane
$\xi,\upsilon$ and the $O\xi$ axis, $0\leq\phi<2\pi$,  
$\chi$ is the angle between the projection of $A$ in the plane $\tau,\zeta$ and
the  $O\tau$ axis, $0\leq\chi<2\pi$,
and $\psi$ is the angle between the line $OA$ and the plane $\tau O \zeta
$, $0\leq \psi\leq\pi/2$, 
as shown in Fig. 1. The angles $\phi$ and $\chi$ will be called azimuthal
angles, the angle $\psi$ will be called planar angle .
The fact that $0\leq \psi\leq\pi/2$ means that $\psi$ has
the same sign on both faces of the two-dimensional hyperplane $\upsilon O
\zeta$. The components of the point $A$
in terms of the distance $d$ and the angles $\phi, \chi, \psi$ are thus
\begin{equation}
\frac{x+t}{\sqrt{2}}=d\cos\phi \sin\psi , 
\label{12a}
\end{equation}
\begin{equation}
\frac{x-t}{\sqrt{2}}=d\cos\chi \cos\psi , 
\label{12b}
\end{equation}
\begin{equation}
\frac{y+z}{\sqrt{2}}=d\sin\phi \sin\psi , 
\label{12c}
\end{equation}
\begin{equation}
\frac{y-z}{\sqrt{2}}=d\sin\chi \cos\psi .
\label{12d}
\end{equation}
It can be checked that $\rho_+=\sqrt{2}d\sin\psi, \rho_-=\sqrt{2}d\cos\psi$.
The coordinates $x,y,z,t$ in terms of the variables $d, \phi, \chi,
\psi$ are
\begin{equation}
x=\frac{d}{\sqrt{2}}(\cos\psi\cos\chi+\sin\psi\cos\phi),
\label{12e}
\end{equation}
\begin{equation}
y=\frac{d}{\sqrt{2}}(\cos\psi\sin\chi+\sin\psi\sin\phi),
\label{12g}
\end{equation}
\begin{equation}
z=\frac{d}{\sqrt{2}}(-\cos\psi\sin\chi+\sin\psi\sin\phi),
\label{12h}
\end{equation}
\begin{equation}
t=\frac{d}{\sqrt{2}}(-\cos\psi\cos\chi+\sin\psi\cos\phi).
\label{12f}
\end{equation}
The angles $\phi, \chi, \psi$ can be expressed in terms of the coordinates
$x,y,z,t$ as 
\begin{equation}
\sin\phi = (y+z)/\rho_+ ,\: \cos\phi = (x+t)/\rho_+ ,
\label{13a}
\end{equation}
\begin{equation}
\sin\chi = (y-z)/\rho_- ,\: \cos\chi = (x-t)/\rho_- ,
\label{13b}
\end{equation}
\begin{equation}
\tan\psi=\rho_+/\rho_- .
\label{13c}
\end{equation}
The nodal hyperplanes are $\xi O\upsilon$, for which $\tau=0, \zeta=0$, and
$\tau O\zeta$, for which $\xi=0, \upsilon=0$.
For points in 
the nodal hyperplane $\xi O\upsilon$ the planar angle is $\psi=\pi/2$, 
for points in the nodal hyperplane $\tau O\zeta$ the planar angle is $\psi=0$.

It can be shown that if $u_1=x_1+\alpha y_1+\beta z_1+\gamma t_1, 
u_2=x_2+\alpha y_2+\beta z_2+\gamma t_2$ are circular fourcomplex
numbers of amplitudes and angles $\rho_1, \phi_1, \chi_1, \psi_1$ and
respectively $\rho_2, \phi_2, \chi_2, \psi_2$, then the amplitude $\rho$ and
the angles $\phi, \chi, \psi$ of the product circular fourcomplex number
$u_1u_2$ 
are 
\begin{equation}
\rho=\rho_1\rho_2, 
\label{14a}
\end{equation}
\begin{equation}
 \phi=\phi_1+\phi_2, \: \chi=\chi_1+\chi_2, \: \tan\psi=\tan\psi_1\tan\psi_2 . 
\label{14b}
\end{equation}
The relations (\ref{14a})-(\ref{14b}) are consequences of the definitions
(\ref{6e})-(\ref{7b}), (\ref{13a})-(\ref{13c}) and of the identities
\begin{eqnarray}
\lefteqn{[(x_1x_2-y_1y_2-z_1z_2+t_1t_2)+(x_1t_2+t_1x_2-y_1z_2-z_1y_2)]^2
\nonumber}\\
&& +[(x_1y_2+y_1x_2+z_1t_2+t_1z_2)+(x_1z_2+z_1x_2+y_1t_2+t_1y_2)]^2\nonumber\\
&&=[(x_1+t_1)^2+(y_1+z_1)^2][(x_2+t_2)^2+(y_2+z_2)^2] ,
\label{15}
\end{eqnarray}
\begin{eqnarray}
\lefteqn{[(x_1x_2-y_1y_2-z_1z_2+t_1t_2)-(x_1t_2+t_1x_2-y_1z_2-z_1y_2)]^2
\nonumber}\\
&& +[(x_1y_2+y_1x_2+z_1t_2+t_1z_2)-(x_1z_2+z_1x_2+y_1t_2+t_1y_2)]^2\nonumber\\
&&=[(x_1-t_1)^2+(y_1-z_1)^2][(x_2-t_2)^2+(y_2-z_2)^2] ,
\label{16}
\end{eqnarray}
\begin{eqnarray}
\lefteqn{(x_1x_2-y_1y_2-z_1z_2+t_1t_2)+(x_1t_2+t_1x_2-y_1z_2-z_1y_2)
\nonumber}\\
&& =(x_1+t_1)(x_2+t_2)-(y_1+z_1)(y_2+z_2) ,
\label{17}
\end{eqnarray}
\begin{eqnarray}
\lefteqn{(x_1x_2-y_1y_2-z_1z_2+t_1t_2)-(x_1t_2+t_1x_2-y_1z_2-z_1y_2)
\nonumber}\\
&& =(x_1-t_1)(x_2-t_2)-(y_1-z_1)(y_2-z_2) ,
\label{18}
\end{eqnarray}
\begin{eqnarray}
\lefteqn{(x_1y_2+y_1x_2+z_1t_2+t_1z_2)+(x_1z_2+z_1x_2+y_1t_2+t_1y_2)
\nonumber}\\
&&=(y_1+z_1)(x_2+t_2)+(y_2+z_2)(x_2+t_2) ,
\label{19}
\end{eqnarray}
\begin{eqnarray}
\lefteqn{(x_1y_2+y_1x_2+z_1t_2+t_1z_2)-(x_1z_2+z_1x_2+y_1t_2+t_1y_2)
\nonumber}\\
&&=(y_1-z_1)(x_2-t_2)+(y_2-z_2)(x_2-t_2) .
\label{20}
\end{eqnarray}
The identities (\ref{15}) and (\ref{16}) can also be written as
\begin{equation}
\rho_+^2=\rho_{1+}\rho_{2+} ,
\label{21a}
\end{equation}
\begin{equation}
\rho_-^2=\rho_{1-}\rho_{2-} ,
\label{21b}
\end{equation}
where
\begin{equation}
\rho_{j+}^2=(x_j+t_j)^2+(y_j+z_j)^2, \: \rho_{j-}^2=(x_j-t_j)^2+(y_j-z_j)^2 ,\:
j=1,2. 
\label{22}
\end{equation}

The fact that the amplitude of the product is equal to the product of the 
amplitudes, as written in Eq. (\ref{14a}), can 
be demonstrated also by using a representation of the multiplication of the 
circular fourcomplex numbers by matrices, in which the circular fourcomplex
number $u=x+\alpha 
y+\beta z+\gamma t$ is represented by the matrix
\begin{equation}
A=\left(\begin{array}{cccc}
x&y&z&t\\
-y&x&t&-z\\
-z&t&x&-y\\
t&z&y&x 
\end{array}\right) .
\label{23}
\end{equation}
The product $u=x+\alpha y+\beta z+\gamma t$ of the circular fourcomplex numbers
$u_1=x_1+\alpha y_1+\beta z_1+\gamma t_1, u_2=x_2+\alpha y_2+\beta z_2+\gamma
t_2$, can be represented by the matrix multiplication 
\begin{equation}
A=A_1A_2.
\label{24}
\end{equation}
It can be checked that the determinant ${\rm det}(A)$ of the matrix $A$ is
\begin{equation}
{\rm det}A = \rho^4 .
\label{25}
\end{equation}
The identity (\ref{14a}) is then a consequence of the fact the determinant 
of the product of matrices is equal to the product of the determinants 
of the factor matrices.

\subsection{The exponential and trigonometric forms of circular 
fourcomplex numbers}

The exponential function of the fourcomplex variable $u$ is defined by the
series
\begin{equation}
\exp u = 1+u+u^2/2!+u^3/3!+\cdots . 
\label{g1}
\end{equation}
It can be checked by direct multiplication of the series that
\begin{equation}
\exp(u+u^\prime)=\exp u \cdot \exp u^\prime . 
\label{g2}
\end{equation}
These relations have the same form for all hypercomplex systems discussed in
this work.

If $u=x+\alpha y+\beta z+\gamma t$, then  $\exp u$ can be calculated as
$\exp u=\exp x \cdot \exp (\alpha y) \cdot \exp (\beta z) \cdot \exp (\gamma
t)$. According to Eq. (\ref{1}), 
\begin{equation}
\alpha^{2m}=(-1)^m, \alpha^{2m+1}=(-1)^m\alpha, 
\beta^{2m}=(-1)^m, \beta^{2m+1}=(-1)^m\beta, 
\gamma^m=1,  
\label{28}
\end{equation}
where $m$ is a natural number,
so that $\exp (\alpha y), \: \exp(\beta z)$ and $\exp(\gamma t)$ can be written
as 
\begin{equation}
\exp (\alpha y) = \cos y +\alpha \sin y , \:
\exp (\beta z) = \cos z +\beta \sin z , 
\label{29}
\end{equation}
and
\begin{equation}
\exp (\gamma t) = \cosh t +\gamma \sinh t . 
\label{30}
\end{equation}
From Eqs. (\ref{29})-(\ref{30}) it can be inferred that
\begin{equation}
\begin{array}{l}
(\cos y +\alpha \sin y)^m=\cos my +\alpha \sin my , \\
(\cos z +\beta \sin z)^m=\cos mz +\beta \sin mz , \\
(\cosh t +\gamma \sinh t)^m=\cosh mt +\gamma \sinh mt . 
\end{array}
\label{30b}
\end{equation}

Any circular fourcomplex number $u=x+\alpha y+\beta z+\gamma t$ can be writen
in the 
form 
\begin{equation}
x+\alpha y+\beta z+\gamma t=e^{x_1+\alpha y_1+\beta z_1+\gamma t_1} .
\label{31}
\end{equation}
The expressions of $x_1, y_1, z_1, t_1$ as functions of 
$x, y, z, t$ can be obtained by
developing $e^{\alpha y_1}, e^{\beta z_1}$ and $e^{\gamma t_1}$ with the aid of
Eqs. (\ref{29}) and (\ref{30}), by multiplying these expressions and separating
the hypercomplex components, 
\begin{equation}
x=e^{x_1}(\cos y_1\cos z_1\cosh t_1-\sin y_1\sin z_1\sinh t_1) ,
\label{32}
\end{equation}
\begin{equation}
y=e^{x_1}(\sin y_1\cos z_1\cosh t_1+\cos y_1\sin z_1\sinh t_1) ,
\label{33}
\end{equation}
\begin{equation}
z=e^{x_1}(\cos y_1\sin z_1\cosh t_1+\sin y_1\cos z_1\sinh t_1) ,
\label{34}
\end{equation}
\begin{equation}
t=e^{x_1}(-\sin y_1\sin z_1\cosh t_1+\cos y_1\cos z_1\sinh t_1) ,
\label{35}
\end{equation}
From Eqs. (\ref{32})-(\ref{35}) it can be shown by direct calculation that
\begin{equation}
x^2+y^2+z^2+t^2=e^{2x_1}\cosh 2t_1 ,
\label{36}
\end{equation}
\begin{equation}
2(xt+yz)=e^{2x_1}\sinh 2t_1 ,
\label{37}
\end{equation}
so that
\begin{equation}
e^{4x_1}=(x^2+y^2+z^2+t^2)^2-4(xt+yz)^2 .
\label{38}
\end{equation}
By comparing the expression in the right-hand side of Eq. (\ref{38}) with the
expression of $\rho$, Eq. ({\ref{7a}), it can be seen that
\begin{equation}
e^{x_1}=\rho . 
\label{39}
\end{equation}
The variable $t_1$ is then given by
\begin{equation}
\cosh 2t_1=\frac{d^2}{\rho^2}, \: \sinh 2t_1=\frac{2(xt+yz)}{\rho^2} .
\label{40}
\end{equation}
From the fact that  $\rho^4=d^4-4(xt+yz)^2$ it follows that $d^2/\rho^2\geq 1$,
so that Eq. (\ref{40}) has always a real solution, and $t_1=0$ for $xt+yz=0$.
It can be shown similarly that
\begin{equation}
\cos 2y_1 = \frac{x^2-y^2+z^2-t^2}{\rho^2} , \: 
\sin 2y_1 = \frac{2(xy-zt)}{\rho^2} ,
\label{41}
\end{equation}
\begin{equation}
\cos 2z_1 = \frac{x^2+y^2-z^2-t^2}{\rho^2} , \: 
\sin 2z_1 = \frac{2(xz-yt)}{\rho^2} .
\label{42}
\end{equation}
It can be shown that $(x^2-y^2+z^2-t^2)^2\leq\rho^4$, the equality taking place
for $xy=zt$, and $(x^2+y^2-z^2-t^2)^2\leq\rho^4$, the equality taking place
for $xz=yt$, so that Eqs. (\ref{41}) and Eqs. (\ref{42}) have always real
solutions. 

The variables
\begin{equation}
y_1 =\frac{1}{2}\arcsin \frac{2(xy-zt)}{\rho^2} ,\;
z_1 =\frac{1}{2}\arcsin \frac{2(xz-yt)}{\rho^2} ,\;
t_1 =\frac{1}{2}{\rm argsinh}\frac{2(xt+yz)}{\rho^2} 
\label{40b}
\end{equation}
are additive upon the multiplication of circular fourcomplex numbers, as can be
seen from the identities
\begin{eqnarray}
\lefteqn{(x x^{\prime}-y y^{\prime}-z z^{\prime}+t t^{\prime})(x y^{\prime}+y
x^{\prime}+z t^{\prime}+t z^{\prime})\nonumber}\\ 
&&-(x z^{\prime}+z x^{\prime}+y t^{\prime}+t y^{\prime})(x t^{\prime}-y
z^{\prime}-z y^{\prime}+t x^{\prime})\nonumber\\ 
&&=(x y -z t )(x^{\prime 2}-y^{\prime 2}+z^{\prime 2}-t^{\prime 2})+(x^{ 2}-y^{
2}+z^{ 2}-t^{ 2})(x^{\prime}y^{\prime}-z^{\prime}t^{\prime}), 
\label{51ca}
\end{eqnarray}
\begin{eqnarray}
\lefteqn{(x x^{\prime}-y y^{\prime}-z z^{\prime}+t t^{\prime})(x z^{\prime}+z
x^{\prime}+y t^{\prime}+t y^{\prime}) 
\nonumber}\\
&&-(x y^{\prime}+y x^{\prime}+z t^{\prime}+t z^{\prime})(x t^{\prime}-y
z^{\prime}-z y^{\prime}+t x^{\prime})\nonumber\\ 
&&=(x z -y t )(x^{\prime 2}+y^{\prime 2}-z^{\prime 2}-t^{\prime 2})+(x^{ 2}+y^{
2}-z^{ 2}-t^{ 2})(x^{\prime}z^{\prime}-y^{\prime}t^{\prime}), 
\label{51cb}
\end{eqnarray}
\begin{eqnarray}
\lefteqn{(x x^{\prime}-y y^{\prime}
-z z^{\prime}+t t^{\prime})(x t^{\prime}
-y z^{\prime}-z y^{\prime}+t x^{\prime})
\nonumber}\\
&&+(x y^{\prime}+y x^{\prime}+z t^{\prime}+t z^{\prime})
(x z^{\prime}+z x^{\prime}+y t^{\prime}+t y^{\prime})\nonumber\\
&&=(x t +y z )(x^{\prime 2}+y^{\prime 2}+z^{\prime 2}+t^{\prime 2})
+(x^{ 2}+y^{ 2}+z^{ 2}+t^{ 2})(x^{\prime}t^{\prime}+y^{\prime}z^{\prime}).
\label{51cc}
\end{eqnarray}

The expressions appearing in Eqs. (\ref{40})-(\ref{42}) can be calculated in
terms of the angles $\phi, \chi, \psi$ with the aid of Eqs.
(\ref{12e})-(\ref{12f}) as
\begin{equation}
\frac{d^2}{\rho^2}=\frac{1}{\sin
2\psi},\;\frac{2(xt+yz)}{\rho^2}=-\frac{1}{\tan 2\psi},
\label{51cd}
\end{equation}
\begin{equation}
\frac{x^2-y^2+z^2-t^2}{\rho^2}=\cos(\phi+\chi),
\;\frac{2(xy-zt)}{\rho^2}=\sin(\phi+\chi),
\label{51ce}
\end{equation}
\begin{equation}
\frac{x^2+y^2-z^2-t^2}{\rho^2}=\cos(\phi-\chi),
\;\frac{2(xz-yt)}{\rho^2}=\sin(\phi-\chi).
\label{51cf}
\end{equation}
Then from Eqs. (\ref{40})-(\ref{42}) and (\ref{51cd})-(\ref{51cf}) it results
that 
\begin{equation}
y_1 =\frac{\phi+\chi}{2},\;
z_1=\frac{\phi-\chi}{2},\;t_1=\frac{1}{2}\ln\tan\psi, 
\label{51cg}
\end{equation}
so that the circular fourcomplex number $u$, Eq. (\ref{31}), can be written as
\begin{equation}
u=\rho\exp\left[\alpha\frac{\phi+\chi}{2}+\beta\frac{\phi-\chi}{2}
+\gamma \frac{1}{2}\ln\tan\psi\right] .
\label{51}
\end{equation}
In Eq. (\ref{51}) the circular fourcomplex number $u=x+\alpha y+\beta z+\gamma
t$ is 
written as the product of the amplitude $\rho$ and of an exponential function,
and therefore this form of $u$ will be called the exponential form of the
circular fourcomplex number.
It can be checked that
\begin{equation}
\exp\left(\frac{\alpha+\beta}{2}\phi\right)=\frac{1-\gamma}{2}
+\frac{1+\gamma}{2}\cos\phi+\frac{\alpha+\beta}{2}\sin\phi ,
\label{51b}
\end{equation}
\begin{equation}
\exp\left(\frac{\alpha-\beta}{2}\chi\right)=\frac{1+\gamma}{2}
+\frac{1-\gamma}{2}\cos\chi+\frac{\alpha-\beta}{2}\sin\phi ,
\label{51c}
\end{equation}
which shows that $e^{(\alpha+\beta)\phi/2}$ and $e^{(\alpha-\beta)\chi/2}$ are
periodic 
functions of $\phi$ and respectively $\chi$, with period $2\pi$.

The relations between the variables $y_1, z_1, t_1$ and the angles $\phi, \chi,
\psi$ can be obtained alternatively by
substituting in Eqs. (\ref{32})-(\ref{35}) the expression $e^{x_1}=d/(\cosh
2t_1)^{1/2} $, Eq. (\ref{36}), and summing and subtracting of the relations,
\begin{equation}
\frac{x+t}{\sqrt{2}}=d\cos(y_1+z_1)\sin(\eta+\pi/4) ,
\label{43}
\end{equation}
\begin{equation}
\frac{x-t}{\sqrt{2}}=d\cos(y_1-z_1)\cos(\eta+\pi/4) ,
\label{44}
\end{equation}
\begin{equation}
\frac{y+z}{\sqrt{2}}=d\sin(y_1+z_1)\sin(\eta+\pi/4) ,
\label{45}
\end{equation}
\begin{equation}
\frac{y-z}{\sqrt{2}}=d\sin(y_1-z_1)\cos(\eta+\pi/4) ,
\label{46}
\end{equation}
where the variable $\eta$ is defined by the relations
\begin{equation}
\frac{\cosh t_1}{(\cosh 2t_1)^{1/2}}=\cos\eta, \: 
\frac{\sinh t_1}{(\cosh 2t_1)^{1/2}}=\sin\eta ,
\label{47}
\end{equation}
and when $-\infty<t_1<\infty$, the range of the variable $\eta$ is
$-\pi/4\leq\eta\leq \pi/4$.
The comparison of Eqs. (\ref{12a})-(\ref{12d}) and (\ref{43})-(\ref{46}) shows
that 
\begin{equation}
\phi=y_1+z_1, \: \chi=y_1-z_1, \: \psi=\eta+\pi/4 .
\label{48}
\end{equation}

It can be shown with the aid of Eq. (\ref{30}) that
\begin{equation}
\exp \left(\frac{1}{2}\gamma\ln\tan\psi \right) 
=\frac{1}{(\sin 2\psi)^{1/2}}\left[
\cos(\psi-\pi/4)+\gamma\sin(\psi-\pi/4)\right] . 
\label{30bb}
\end{equation}
The circular fourcomplex number $u$, Eq. (\ref{51}), can then be written
equivalently as 
\begin{equation}
u=d\{\cos(\psi-\pi/4)+\gamma\sin(\psi-\pi/4)\}
\exp\left[\alpha\frac{\phi+\chi}{2}+\beta\frac{\phi-\chi}{2}\right].
\label{52}
\end{equation}
In Eq. (\ref{52}), the circular fourcomplex number 
$u=x+\alpha y+\beta z+\gamma t$ is
written as the product of the modulus $d$ and of factors depending on the
geometric angles $\phi, \chi$ and $\psi$, and this form will be called the
trigonometric form of the circular fourcomplex number.

If $u_1, u_2$ are circular fourcomplex numbers of moduli and angles $d_1,
\phi_1, 
\chi_1, \psi_1$ and respectively $d_2, \phi_2, \chi_2, \psi_2$, the product of
the planar factors can be calculated to be
\begin{eqnarray}
\lefteqn{[\cos(\psi_1-\pi/4)+\gamma\sin(\psi_1-\pi/4)] 
[\cos(\psi_2-\pi/4)+\gamma\sin(\psi_2-\pi/4)]}\nonumber\\
&&=[\cos(\psi_1-\psi_2)-\gamma\cos(\psi_1+\psi_2)] .
\label{53}
\end{eqnarray}
The right-hand side of Eq. (\ref{53}) can be written as
\begin{eqnarray}
\lefteqn{\cos(\psi_1-\psi_2)-\gamma\cos(\psi_1+\psi_2)}\nonumber\\
&&=[2(\cos^2\psi_1\cos^2\psi_2+\sin^2\psi_1\sin^2\psi_2)]^{1/2}
[\cos(\psi-\pi/4)+\gamma\sin(\psi-\pi/4)] ,
\label{54}
\end{eqnarray}
where the angle $\psi$, determined by the condition that
\begin{equation}
\tan(\psi-\pi/4)=-\cos(\psi_1+\psi_2)/\cos(\psi_1-\psi_2)
\label{55}
\end{equation}
is given by $\tan\psi=\tan\psi_1\tan\psi_2$ ,
which is consistent with Eq. (\ref{14b}).
It can be checked that the modulus $d$ of the product $u_1u_2$ is 
\begin{equation}
d=\sqrt{2}d_1d_2\left(\cos^2\psi_1\cos^2\psi_2
+\sin^2\psi_1\sin^2\psi_2\right)^{1/2}. 
\label{56}
\end{equation}

\subsection{Elementary functions of a circular fourcomplex variable}

The logarithm $u_1$ of the circular fourcomplex number $u$, $u_1=\ln u$, can be
defined 
as the solution of the equation
\begin{equation}
u=e^{u_1} ,
\label{57}
\end{equation}
written explicitly previously in Eq. (\ref{31}), for $u_1$ as a function of
$u$. From Eq. (\ref{51}) it results that 
\begin{equation}
\ln u=\ln \rho+\frac{1}{2}\gamma \ln\tan\psi
+\alpha\frac{\phi+\chi}{2}+\beta\frac{\phi-\chi}{2} .
\label{58}
\end{equation}
It can be inferred from Eqs. (\ref{14a}) and (\ref{14b}) that
\begin{equation}
\ln(uu^\prime)=\ln u+\ln u^\prime ,
\label{59}
\end{equation}
up to multiples of $\pi(\alpha+\beta)$ and $\pi(\alpha-\beta)$.

The power function $u^n$ can be defined for real values of $n$ as
\begin{equation}
u^m=e^{m\ln u} .
\label{60}
\end{equation}
The power function is multivalued unless $n$ is an integer. 
For integer $n$, it can be inferred from Eq. (\ref{59}) that
\begin{equation}
(uu^\prime)^m=u^n\:u^{\prime m} .
\label{61}
\end{equation}
If, for example, $m=2$, it can be checked with the aid of Eq. (\ref{52})
that Eq. (\ref{60}) gives indeed $(x+\alpha y+\beta z+\gamma t)^2=
x^2-y^2-z^2+t^2+2\alpha(xy+zt)+2\beta(xz+yt)+2\gamma(2xt-yz)$.

The trigonometric functions of the fourcomplex variable $u$
are defined by the series
\begin{equation}
\cos u = 1 - u^2/2!+u^4/4!+\cdots, 
\label{g3}
\end{equation}
\begin{equation}
\sin u=u-u^3/3!+u^5/5! +\cdots .
\label{g4}
\end{equation}
It can be checked by series multiplication that the usual addition theorems
hold also for the circular fourcomplex numbers $u, u^\prime$,
\begin{equation}
\cos(u+u^\prime)=\cos u\cos u^\prime - \sin u \sin u^\prime ,
\label{g5}
\end{equation}
\begin{equation}
\sin(u+u^\prime)=\sin u\cos u^\prime + \cos u \sin u^\prime .
\label{g6}
\end{equation}
These relations have the same form for all systems of hypercomplex numbers
discussed in this work.
The cosine and sine functions of the hypercomplex variables $\alpha y, 
\beta z$ and $ \gamma t$ can be expressed as
\begin{equation}
\cos\alpha y=\cosh y, \: \sin\alpha y=\alpha\sinh y, 
\label{66}
\end{equation}
\begin{equation}
\cos\beta y=\cosh y, \: \sin\beta y=\beta\sinh y, 
\label{67}
\end{equation}
\begin{equation}
\cos\gamma y=\cos y, \: \sin\gamma y=\gamma\sin y .
\label{68}
\end{equation}
The cosine and sine functions of a circular fourcomplex number $x+\alpha
y+\beta 
z+\gamma t$ can then be
expressed in terms of elementary functions with the aid of the addition
theorems Eqs. (\ref{g5}), (\ref{g6}) and of the expressions in  Eqs. 
(\ref{66})-(\ref{68}).

The hyperbolic functions of the fourcomplex variable
$u$ are defined by the series
\begin{equation}
\cosh u = 1 + u^2/2!+u^4/4!+\cdots, 
\label{g7}
\end{equation}
\begin{equation}
\sinh u=u+u^3/3!+u^5/5! +\cdots .
\label{g8}
\end{equation}
It can be checked by series multiplication that the usual addition theorems
hold also for the circular fourcomplex numbers $u, u^\prime$,
\begin{equation}
\cosh(u+u^\prime)=\cosh u\cosh u^\prime + \sinh u \sinh u^\prime ,
\label{g9}
\end{equation}
\begin{equation}
\sinh(u+u^\prime)=\sinh u\cosh u^\prime + \cosh u \sinh u^\prime .
\label{g10}
\end{equation}
These relations have the same form for all systems of hypercomplex numbers
discussed in this work.
The hyperbolic cosine and sine functions of the hypercomplex variables 
$\alpha y,  \beta z$ and $ \gamma t$ can be expressed as
\begin{equation}
\cosh\alpha y=\cos y, \: \sinh\alpha y=\alpha\sin y, 
\label{73}
\end{equation}
\begin{equation}
\cosh\beta y=\cos y, \: \sinh\beta y=\beta\sin y, 
\label{74}
\end{equation}
\begin{equation}
\cosh\gamma y=\cosh y, \: \sinh\gamma y=\gamma\sinh y .
\label{75}
\end{equation}
The hyperbolic cosine and sine functions of a circular fourcomplex number 
$x+\alpha y+\beta z+\gamma t$ can then be
expressed in terms of elementary functions with the aid of the addition
theorems Eqs. (\ref{g9}), (\ref{g10}) and of the expressions in  Eqs. 
(\ref{73})-(\ref{75}).

\subsection{Power series of circular fourcomplex variables}

A circular fourcomplex series is an infinite sum of the form
\begin{equation}
a_0+a_1+a_2+\cdots+a_n+\cdots , 
\label{76}
\end{equation}
where the coefficients $a_n$ are circular fourcomplex numbers. The convergence
of  
the series (\ref{76}) can be defined in terms of the convergence of its 4 real
components. The convergence of a circular fourcomplex series can however be
studied 
using circular fourcomplex variables. The main criterion for absolute
convergence  
remains the comparison theorem, but this requires a number of inequalities
which will be discussed further.

The modulus of a circular fourcomplex number $u=x+\alpha y+\beta z+\gamma t$
can be 
defined as 
\begin{equation}
|u|=(x^2+y^2+z^2+t^2)^{1/2} ,
\label{77}
\end{equation}
so that, according to Eq. (\ref{10}), $d=|u|$. Since $|x|\leq |u|, |y|\leq |u|,
|z|\leq |u|, |t|\leq |u|$, a property of 
absolute convergence established via a comparison theorem based on the modulus
of the series (\ref{76}) will ensure the absolute convergence of each real
component of that series.

The modulus of the sum $u_1+u_2$ of the circular fourcomplex numbers $u_1, u_2$
fulfils 
the inequality
\begin{equation}
||u_1|-|u_2||\leq |u_1+u_2|\leq |u_1|+|u_2| .
\label{78}
\end{equation}
For the product the relation is 
\begin{equation}
|u_1u_2|\leq \sqrt{2}|u_1||u_2| ,
\label{79}
\end{equation}
which replaces the relation of equality extant for regular complex numbers.
The equality in Eq. (\ref{79}) takes place for $x_1=t_1, y_1=z_1,
x_2=t_2, y_2=z_2$ or $x_1=-t_1, y_1=-z_1,
x_2=-t_2, y_2=-z_2$. In Eq. (\ref{56}), this corresponds to $\psi_1=0,\psi_2=0$
or $\psi_1=\pi/2, \psi_2=\pi/2$.
The modulus of a product, which has the property that
$0\leq|u_1u_2|$, becomes equal to zero for
$x_1=t_1, y_1=z_1, x_2=-t_2, y_2=-z_2$ or $x_1=-t_1, y_1=-z_1,
x_2=t_2, y_2=z_2$, as discussed after Eq. (\ref{9}).
In Eq. (\ref{56}), the latter situation corresponds to $\psi_1=0,\psi_2=\pi/2$
or $\psi_1=0, \psi_2=\pi/2$.

It can be shown that
\begin{equation}
x^2+y^2+z^2+t^2\leq|u^2|\leq \sqrt{2}(x^2+y^2+z^2+t^2) .
\label{80}
\end{equation}
The left relation in Eq. (\ref{80}) becomes an equality, 
$x^2+y^2+z^2+t^2=|u^2|$, for $xt+yz=0$. This condition corresponds to
$\psi_1=\psi_2=\pi/4$ in Eq. (\ref{56}).
The inequality in Eq. (\ref{79}) implies that
\begin{equation}
|u^l|\leq 2^{(l-1)/2}|u|^l .
\label{81}
\end{equation}
From Eqs. (\ref{79}) and (\ref{81}) it results that
\begin{equation}
|au^l|\leq 2^{l/2} |a| |u|^l .
\label{82}
\end{equation}

A power series of the circular fourcomplex variable $u$ is a series of the form
\begin{equation}
a_0+a_1 u + a_2 u^2+\cdots +a_l u^l+\cdots .
\label{83}
\end{equation}
Since
\begin{equation}
\left|\sum_{l=0}^\infty a_l u^l\right| \leq  \sum_{l=0}^\infty
2^{l/2}|a_l| |u|^l ,
\label{84}
\end{equation}
a sufficient condition for the absolute convergence of this series is that
\begin{equation}
\lim_{l\rightarrow \infty}\frac{\sqrt{2}|a_{l+1}||u|}{|a_l|}<1 .
\label{85}
\end{equation}
Thus the series is absolutely convergent for 
\begin{equation}
|u|<c,
\label{86}
\end{equation}
where 
\begin{equation}
c=\lim_{l\rightarrow\infty} \frac{|a_l|}{\sqrt{2}|a_{l+1}|} .
\label{87}
\end{equation}

The convergence of the series (\ref{83}) can be also studied with the aid of
the transformation 
\begin{equation}
x+\alpha y+\beta z+\gamma t=\sqrt{2}(e_1\xi+\tilde e_1 \upsilon+e_2\tau
+\tilde e_2\zeta) , 
\label{87b}
\end{equation}
where $\xi,\upsilon, \tau, \zeta$ have been defined in Eq. (\ref{11}),
and
\begin{equation}
e_1=\frac{1+\gamma}{2},\:\:\tilde e_1=\frac{\alpha+\beta}{2},\:\:
e_2=\frac{1-\gamma}{2},\:\:\tilde e_2=\frac{\alpha-\beta}{2}.
\label{87c}
\end{equation}
The ensemble $e_1, \tilde e_1, e_2, \tilde e_2$ will be called the canonical
circular fourcomplex base, and Eq. (\ref{87b}) gives the canonical form of the
circular fourcomplex number.
It can be checked that
\begin{eqnarray}
\lefteqn{e_1^2=e_1, \:\:\tilde e_1^2=-e_1,\:\: e_1\tilde e_1=\tilde e_1,\:\:
e_2^2=e_2, \:\:\tilde e_2^2=-e_2,\:\: e_2\tilde e_2=\tilde e_2,\:\:\nonumber}\\
&&e_1e_2=0,\:\: \tilde e_1\tilde e_2=0, \:\:e_1\tilde e_2=0, \:\:
e_2\tilde e_1=0.  
\label{87d}
\end{eqnarray}
The moduli of the bases in Eq. (\ref{87c}) are
\begin{equation}
|e_1|=\frac{1}{\sqrt{2}},\;|\tilde e_1|=\frac{1}{\sqrt{2}},\;
|e_2|=\frac{1}{\sqrt{2}},\;|\tilde e_2|=\frac{1}{\sqrt{2}},
\label{87e}
\end{equation}
and it can be checked that
\begin{equation}
|x+\alpha y+\beta z+\gamma t|^2=\xi^2+\upsilon^2+\tau^2+\zeta^2.
\label{87f}
\end{equation}

If $u=u^\prime u^{\prime\prime}$, the components $\xi,\upsilon, \tau, \zeta$
are related, according to Eqs. (\ref{17})-(\ref{20}) by 
\begin{equation}
\xi=\sqrt{2}(\xi^\prime \xi^{\prime\prime}-\upsilon^\prime
\upsilon^{\prime\prime}), \:\: 
\upsilon=\sqrt{2}(\xi^\prime \upsilon^{\prime\prime}+\upsilon^\prime
\xi^{\prime\prime}), \:\: 
\tau=\sqrt{2}(\tau^\prime \tau^{\prime\prime}-\zeta^\prime
\zeta^{\prime\prime}), \:\: 
\zeta=\sqrt{2}(\tau^\prime \zeta^{\prime\prime}+\zeta^\prime
\tau^{\prime\prime}), \:\: 
\label{87g}
\end{equation}
which show that, upon multiplication, the components $\xi,\upsilon$ and $\tau,
\zeta$ obey, up to a normalization constant, the same
rules as the real and imaginary components of usual, two-dimensional complex
numbers.

If the coefficients in Eq. (\ref{83}) are 
\begin{equation}
a_l= a_{l0}+\alpha a_{l1}+\beta a_{l2}+\gamma a_{l3}, 
\label{n88a}
\end{equation}
and
\begin{equation}
A_{l1}=a_{l0}+a_{l3},\;
\tilde A_{l1}=a_{l1}+a_{l2},\;
A_{l2}=a_{l0}-a_{l3},\;
\tilde A_{l2}=a_{l1}-a_{l2},
\label{n88b}
\end{equation}
the series (\ref{83}) can be written as
\begin{equation}
\sum_{l=0}^\infty 2^{l/2}\left[
(e_1 A_{l1}+\tilde e_1\tilde A_{l1})(e_1 \xi+\tilde e_1 \upsilon)^l 
+(e_2 A_{l2}+\tilde e_2\tilde A_{l2})(e_2 \tau+\tilde e_2 \zeta)^l 
\right].
\label{n89a}
\end{equation}
Thus, the series in Eqs. (\ref{83}) and (\ref{n89a}) are
absolutely convergent for   
\begin{equation}
\rho_+<c_1, \;\rho_-<c_2,
\label{n90}
\end{equation}
where 
\begin{equation}
c_1=\lim_{l\rightarrow\infty} \frac
{\left[A_{l1}^2+\tilde A_{l1}^2\right]^{1/2}}
{\sqrt{2}\left[A_{l+1,1}^2+\tilde A_{l+1,1}^2\right]^{1/2}},\;\;
c_2=\lim_{l\rightarrow\infty} \frac
{\left[A_{l2}^2+\tilde A_{l2}^2\right]^{1/2}}
{\sqrt{2}\left[A_{l+1,2}^2+\tilde A_{l+1,2}^2\right]^{1/2}}.
\label{n91}
\end{equation}

It can be shown that $c=(1/\sqrt{2}){\rm
min}(c_1,c_2)$, where ${\rm min}$ designates the smallest of
the numbers $c_1,c_2$. Using the expression of $|u|$ in
Eq. (\ref{87f}),  it can be seen that the spherical region of
convergence defined in Eqs. (\ref{86}), (\ref{87}) is included in the
cylindrical region of convergence defined in Eqs. (\ref{n90}) and (\ref{n91}).

\subsection{Analytic functions of circular fourcomplex variables}

The derivative  
of a function $f(u)$ of the fourcomplex variables $u$ is
defined as a function $f^\prime (u)$ having the property that
\begin{equation}
|f(u)-f(u_0)-f^\prime (u_0)(u-u_0)|\rightarrow 0 \:\:{\rm as} 
\:\:|u-u_0|\rightarrow 0 . 
\label{gs88}
\end{equation}
If the difference $u-u_0$ is not parallel to one of the nodal hypersurfaces,
the definition in Eq. (\ref{gs88}) can also 
be written as
\begin{equation}
f^\prime (u_0)=\lim_{u\rightarrow u_0}\frac{f(u)-f(u_0)}{u-u_0} .
\label{g12}
\end{equation}
The derivative of the function $f(u)=u^m $, with $m$ an integer, 
is $f^\prime (u)=mu^{m-1}$, as can be seen by developing $u^m=[u_0+(u-u_0)]^m$
as
\begin{equation}
u^m=\sum_{p=0}^{m}\frac{m!}{p!(m-p)!}u_0^{m-p}(u-u_0)^p,
\label{gs90}
\end{equation}
and using the definition (\ref{gs88}).

If the function $f^\prime (u)$ defined in Eq. (\ref{gs88}) is independent of
the 
direction in space along which $u$ is approaching $u_0$, the function $f(u)$ 
is said to be analytic, analogously to the case of functions of regular complex
variables. \cite{3} 
The function $u^m$, with $m$ an integer, 
of the fourcomplex variable $u$ is analytic, because the
difference $u^m-u_0^m$ is always proportional to $u-u_0$, as can be seen from
Eq. (\ref{gs90}). Then series of
integer powers of $u$ will also be analytic functions of the fourcomplex
variable $u$, and this result holds in fact for any commutative algebra. 

If an analytic function is defined by a series around a certain point, for
example $u=0$, as
\begin{equation}
f(u)=\sum_{k=0}^\infty a_k u^k ,
\label{gs91a}
\end{equation}
an expansion of $f(u)$ around a different point $u_0$,
\begin{equation}
f(u)=\sum_{k=0}^\infty c_k (u-u_0)^k ,
\label{g14}
\end{equation}
can be obtained by
substituting in Eq. (\ref{gs91a}) the expression of $u^k$ according to Eq.
(\ref{gs90}). Assuming that the series are absolutely convergent so that the
order of the terms can be modified and ordering the terms in the resulting
expression according to the increasing powers of $u-u_0$ yields
\begin{equation}
f(u)=\sum_{k,l=0}^\infty \frac{(k+l)!}{k!l!}a_{k+l} u_0^l (u-u_0)^k .
\label{gs91b}
\end{equation}
Since the derivative of order $k$ at $u=u_0$ of the function $f(u)$ , Eq.
(\ref{gs91a}), is 
\begin{equation}
f^{(k)}(u_0)=\sum_{l=0}^\infty \frac{(k+l)!}{l!}a_{k+l} u_0^l ,
\label{gs91c}
\end{equation}
the expansion of $f(u)$ around $u=u_0$, Eq. (\ref{gs91b}), becomes
\begin{equation}
f(u)=\sum_{k=0}^\infty \frac{1}{k!} f^{(k)}(u_0)(u-u_0)^k ,
\label{gs91d}
\end{equation}
which has the same form as the series expansion of 2-dimensional complex
functions. 
The relation (\ref{gs91d}) shows that the coefficients in the series expansion,
Eq. (\ref{g14}), are
\begin{equation}
c_k=\frac{1}{k!}f^{(k)}(u_0) .
\label{g15}
\end{equation}

The rules for obtaining the derivatives and the integrals of the basic
functions can 
be obtained from the series of definitions and, as long as these series
expansions have the same form as the corresponding series for the
2-dimensional complex functions, the rules of derivation and integration remain
unchanged.

If the fourcomplex function $f(u)$
of the fourcomplex variable $u$ can be expressed in terms of 
the real functions $P(x,y,z,t),Q(x,y,z,t),R(x,y,z,t), S(x,y,z,t)$ of real
variables $x,y,z,t$ as 
\begin{equation}
f(u)=P(x,y,z,t)+\alpha Q(x,y,z,t)+\beta R(x,y,z,t)+\gamma S(x,y,z,t),
\label{g16}
\end{equation}
then relations of equality 
exist between partial derivatives of the functions $P,Q,R,S$. These relations
can be obtained by writing the derivative of the function $f$ as
\begin{eqnarray}
\lim_{u\rightarrow u_0}
\lefteqn{\frac{1}{\Delta x+\alpha \Delta y +\beta\Delta z+\gamma\Delta t} 
\left[\frac{\partial P}{\partial x}\Delta x+
\frac{\partial P}{\partial y}\Delta y+
\frac{\partial P}{\partial z}\Delta z+
\frac{\partial P}{\partial t}\Delta t\right.\nonumber}\\
&&+\alpha\left(\frac{\partial Q}{\partial x}\Delta x+
\frac{\partial Q}{\partial y}\Delta y+
\frac{\partial Q}{\partial z}\Delta z 
+\frac{\partial Q}{\partial t}\Delta t \right)
+\beta\left(\frac{\partial R}{\partial x}\Delta x+
\frac{\partial R}{\partial y}\Delta y+
\frac{\partial R}{\partial z}\Delta z
+\frac{\partial R}{\partial t}\Delta t
\right)\nonumber\\
&&\left.+\gamma\left(\frac{\partial S}{\partial x}\Delta x+
\frac{\partial S}{\partial y}\Delta y+
\frac{\partial S}{\partial z}\Delta z
+\frac{\partial S}{\partial t}\Delta t\right)\right] ,
\label{g17}
\end{eqnarray}
where the difference appearing in Eq. (\ref{g12}) is
$u-u_0=\Delta x+\alpha\Delta y +\beta\Delta z+\gamma\Delta t$. 
These relations have the same form for all systems of hypercomplex numbers
discussed in this work.

For the present system of hypercomplex numbers, the
relations between the partial derivatives of the functions $P, Q, R, S$ are
obtained by setting succesively in   
Eq. (\ref{g17}) $\Delta x\rightarrow 0, \Delta y=\Delta z=\Delta t=0$;
then $ \Delta y\rightarrow 0, \Delta x=\Delta z=\Delta t=0;$  
then $  \Delta z\rightarrow 0,\Delta x=\Delta y=\Delta t=0$; and finally
$ \Delta t\rightarrow 0,\Delta x=\Delta y=\Delta z=0 $. 
The relations are 
\begin{equation}
\frac{\partial P}{\partial x} = \frac{\partial Q}{\partial y} =
\frac{\partial R}{\partial z} = \frac{\partial S}{\partial t},
\label{95}
\end{equation}
\begin{equation}
\frac{\partial Q}{\partial x} = -\frac{\partial P}{\partial y} =
-\frac{\partial S}{\partial z} = \frac{\partial R}{\partial t},
\label{96}
\end{equation}
\begin{equation}
\frac{\partial R}{\partial x} = -\frac{\partial S}{\partial y} =
-\frac{\partial P}{\partial z} = \frac{\partial Q}{\partial t},
\label{97}
\end{equation}
\begin{equation}
\frac{\partial S}{\partial x} = \frac{\partial R}{\partial y} =
\frac{\partial Q}{\partial z} = \frac{\partial P}{\partial t}.
\label{98}
\end{equation}
The relations (\ref{95})-(\ref{98}) are analogous to the Riemann relations
for the real and imaginary components of a complex function. It can be shown
from Eqs. (\ref{95})-(\ref{98}) that the component $P$ is a solution
of the equations 
\begin{equation}
\frac{\partial^2 P}{\partial x^2}+\frac{\partial^2 P}{\partial y^2}=0,
\:\: 
\frac{\partial^2 P}{\partial x^2}+\frac{\partial^2 P}{\partial z^2}=0,
\:\:
\frac{\partial^2 P}{\partial y^2}+\frac{\partial^2 P}{\partial t^2}=0,
\:\:
\frac{\partial^2 P}{\partial z^2}+\frac{\partial^2 P}{\partial t^2}=0,
\:\:
\label{99}
\end{equation}
\begin{equation}
\frac{\partial^2 P}{\partial x^2}-\frac{\partial^2 P}{\partial t^2}=0,
\:\:
\frac{\partial^2 P}{\partial y^2}-\frac{\partial^2 P}{\partial z^2}=0,
\label{100}
\end{equation}
and the components $Q, R, S$ are solutions of similar equations.

As can be seen from Eqs. (\ref{99})-(\ref{100}), the components $P, Q, R, S$ of
an analytic function of circular fourcomplex variable are harmonic 
with respect to the pairs of variables $x,y; x,z; y,t$ and $ z,t$, and are
solutions of the wave 
equation with respect to the pairs of variables $x,t$ and $y,z$.
The components $P, Q, R, S$ are also solutions of the mixed-derivative
equations 
\begin{equation}
\frac{\partial^2 P}{\partial x\partial y}=\frac{\partial^2 P}{\partial
z\partial t} ,
\:\: 
\frac{\partial^2 P}{\partial x\partial z}=\frac{\partial^2 P}{\partial
y\partial t} ,
\:\: 
\frac{\partial^2 P}{\partial x\partial t}=-\frac{\partial^2 P}{\partial
y\partial z} ,
\label{107}
\end{equation}
and the components $Q, R, S$ are solutions of similar equations.

\subsection{Integrals of functions of circular fourcomplex variables}

The singularities of circular fourcomplex functions arise from terms of the
form 
$1/(u-u_0)^m$, with $m>0$. Functions containing such terms are singular not
only at $u=u_0$, but also at all points of the two-dimensional hyperplanes
passing through $u_0$ and which are parallel to the nodal hyperplanes. 

The integral of a circular fourcomplex function between two points $A, B$ along
a path 
situated in a region free of singularities is independent of path, which means
that the integral of an analytic function along a loop situated in a region
free from singularities is zero,
\begin{equation}
\oint_\Gamma f(u) du = 0,
\label{111}
\end{equation}
where it is supposed that a surface $\Sigma$ spanning 
the closed loop $\Gamma$ is not intersected by any of
the two-dimensional hyperplanes associated with the
singularities of the function $f(u)$. Using the expression, Eq. (\ref{g16}),
for $f(u)$ and the fact that $du=dx+\alpha  dy+\beta dz+\gamma dt$, the
explicit form of the integral in Eq. (\ref{111}) is
\begin{eqnarray}
\lefteqn{\oint _\Gamma f(u) du = \oint_\Gamma
[(Pdx-Qdy-Rdz+Sdt)+\alpha(Qdx+Pdy+Sdz+Rdt)\nonumber}\\
&&+\beta(Rdx+Sdy+Pdz+Qdt)+\gamma(Sdx-Rdy-Qdz+Pdt)] .
\label{112}
\end{eqnarray}
If the functions $P, Q, R, S$ are regular on a surface $\Sigma$
spanning the loop $\Gamma$,
the integral along the loop $\Gamma$ can be transformed with the aid of the
theorem of Stokes in an integral over the surface $\Sigma$ of terms of the form
$\partial P/\partial y +  \partial Q/\partial x, \:\:
\partial P/\partial z + \partial R/\partial x,\:\:
\partial P/\partial t - \partial S/\partial x, \:\:
\partial Q/\partial z -  \partial R/\partial y, \:\:
\partial Q/\partial t + \partial S/\partial y,\:\:
\partial R/\partial t + \partial S/\partial z$ and of similar terms arising
from the $\alpha, \beta$ and $\gamma$ components, 
which are equal to zero by Eqs. (\ref{95})-(\ref{98}), and this proves Eq.
(\ref{111}). 

The integral of the function $(u-u_0)^m$ on a closed loop $\Gamma$ is equal to
zero for $m$ a positive or negative integer not equal to -1,
\begin{equation}
\oint_\Gamma (u-u_0)^m du = 0, \:\: m \:\:{\rm integer},\: m\not=-1 .
\label{112b}
\end{equation}
This is due to the fact that $\int (u-u_0)^m du=(u-u_0)^{m+1}/(m+1), $ and to
the 
fact that the function $(u-u_0)^{m+1}$ is singlevalued for $m$ an integer.

The integral $\oint du/(u-u_0)$ can be calculated using the exponential form 
(\ref{51}),
\begin{equation}
u-u_0=\rho\exp\left(\alpha\frac{\phi+\chi}{2}+\beta\frac{\phi-\chi}{2}
+\gamma \ln\tan\psi\right) ,
\label{113}
\end{equation}
so that 
\begin{equation}
\frac{du}{u-u_0}=\frac{d\rho}{\rho}+\frac{\alpha+\beta}{2}d\phi
+\frac{\alpha-\beta}{2}d\chi+\gamma d\ln\tan\psi .
\label{114}
\end{equation}
Since $\rho$ and $\ln\tan\psi$ are singlevalued variables, it follows that
$\oint_\Gamma d\rho/\rho =0, \oint_\Gamma d\ln\tan\psi=0$. On the other hand,
$\phi$ 
and $\chi$ are cyclic variables, so that they may give a contribution to the
integral around the closed loop $\Gamma$.
Thus, if $C_+$ is a circle of radius $r$
parallel to the $\xi O\upsilon$ plane, and the
projection of the center of this circle on the $\xi O\upsilon$ plane
coincides with the projection of the point $u_0$ on this plane, the points
of the circle $C_+$ are described according to Eqs.
(\ref{11})-(\ref{12d}) by the equations
\begin{eqnarray}
\lefteqn{\xi=\xi_0+r \sin\psi\cos\phi , \:
\upsilon=\upsilon_0+r \sin\psi\sin\phi , \:
\tau=\tau_0+r\cos\psi \cos\chi , \nonumber}\\
&&\zeta=\zeta_0+r \cos\psi\sin\chi , 
\label{115}
\end{eqnarray}
for constant values of $\chi$ and $\psi, \:\psi\not=0, \pi/2$, where
$u_0=x_0+\alpha y_0+\beta 
z_0+\gamma t_0$,  and $\xi_0, \upsilon_0, \tau_0, \zeta_0$ are calculated from
$x_0, y_0, z_0, t_0$ according to Eqs. (\ref{11}).
Then
\begin{equation}
\oint_{C_+}\frac{du}{u-u_0}=\pi(\alpha+\beta) .
\label{116}
\end{equation}
If $C_-$ is a circle of radius $r$
parallel to the $\tau O\zeta$ plane,
and the projection of the center of this circle on the $\tau O\zeta$ plane
coincides with the projection of the point $u_0$ on this plane, the points
of the circle $C_-$ are described by the same Eqs. (\ref{115}) 
but for constant values of $\phi$ and $\psi, \:\psi\not=0, \pi/2$
Then
\begin{equation}
\oint_{C_-}\frac{du}{u-u_0}=\pi(\alpha-\beta) .
\label{117}
\end{equation}
The expression of $\oint_\Gamma du/(u-u_0)$ can be written as a single
equation with the aid of a functional int($M,C$) defined for a point $M$ and a
closed curve $C$ in a two-dimensional plane, such that 
\begin{equation}
{\rm int}(M,C)=\left\{
\begin{array}{l}
1 \;\:{\rm if} \;\:M \;\:{\rm is \;\:an \;\:interior \;\:point \;\:of} \;\:C
,\\  
0 \;\:{\rm if} \;\:M \;\:{\rm is \;\:exterior \;\:to}\:\; C .\\
\end{array}\right.
\label{118}
\end{equation}
With this notation the result of the integration along a closed path
$\Gamma$ can be written as 
\begin{equation}
\oint_\Gamma\frac{du}{u-u_0}=
\pi(\alpha+\beta) \:\;{\rm int}(u_{0\xi\upsilon},\Gamma_{\xi\upsilon})
+\pi (\alpha-\beta)\:\;{\rm int}(u_{0\tau\zeta},\Gamma_{\tau\zeta}),
\label{119}
\end{equation}
where $u_{0\xi\upsilon}, u_{0\tau\zeta}$ and $\Gamma_{\xi\upsilon},
\Gamma_{\tau\zeta}$ are respectively the projections of the point $u_0$ and of
the loop $\Gamma$ on the planes $\xi \upsilon$ and $\tau \zeta$.

If $f(u)$ is an analytic circular fourcomplex function which can be expanded in
a 
series as written in Eq. (\ref{g14}), and the expansion holds on the curve
$\Gamma$ and on a surface spanning $\Gamma$, then from Eqs. (\ref{112b}) and
(\ref{119}) it follows that
\begin{equation}
\oint_\Gamma \frac{f(u)du}{u-u_0}=
\pi[(\alpha+\beta) \:\;{\rm int}(u_{0\xi\upsilon},\Gamma_{\xi\upsilon})
+ (\alpha-\beta)\:\;{\rm int}(u_{0\tau\zeta},\Gamma_{\tau\zeta})]\;f(u_0) ,
\label{120}
\end{equation}
where $\Gamma_{\xi\upsilon}, \Gamma_{\tau\zeta}$ are the projections of 
the curve $\Gamma$ on the planes $\xi \upsilon$ and respectively $\tau \zeta$,
as shown in Fig. 2.
Substituting in the right-hand side of 
Eq. (\ref{120}) the expression of $f(u)$ in terms of the real 
components $P, Q, R, S$, Eq. (\ref{g16}), yields
\begin{eqnarray}
\lefteqn{\oint_\Gamma \frac{f(u)du}{u-u_0}=
\pi [-(1+\gamma)(Q+R)+(\alpha+\beta)(P+S)] 
\:\;{\rm int}(u_{0\xi\upsilon},\Gamma_{\xi\upsilon})\nonumber}\\
&&+\pi [-(1-\gamma)(Q-R)+(\alpha-\beta)(P-S)] 
\:\;{\rm int}(u_{0\tau\zeta},\Gamma_{\tau\zeta}) ,
\label{121}
\end{eqnarray}
where $P, Q, R, S$ are the values of the components of $f$ at $u=u_0$.

If $f(u)$ can be expanded as written in Eq. (\ref{g14}) on 
$\Gamma$ and on a surface spanning $\Gamma$, then from Eqs. (\ref{112b}) and
(\ref{119}) it also results that
\begin{equation}
\oint_\Gamma \frac{f(u)du}{(u-u_0)^{m+1}}=
\frac{\pi}{m!}[(\alpha+\beta) \:\;{\rm
int}(u_{0\xi\upsilon},\Gamma_{\xi\upsilon}) 
+ (\alpha-\beta)\:\;{\rm int}(u_{0\tau\zeta},\Gamma_{\tau\zeta})]\;
f^{(m)}(u_0) ,
\label{122}
\end{equation}
where it has been used the fact that the derivative $f^{(m)}(u_0)$ of order $m$
of $f(u)$ at $u=u_0$ is related to the expansion coefficient in Eq. (\ref{g14})
according to Eq. (\ref{g15}).

If a function $f(u)$ is expanded in positive and negative powers of $u-u_j$,
where $u_j$ are circular fourcomplex constants, $j$ being an index, the
integral of $f$ 
on a closed loop $\Gamma$ is determined by the terms in the expansion of $f$
which are of the form $a_j/(u-u_j)$,
\begin{equation}
f(u)=\cdots+\sum_j\frac{a_j}{u-u_j}+\cdots
\label{123}
\end{equation}
Then the integral of $f$ on a closed loop $\Gamma$ is
\begin{equation}
\oint_\Gamma f(u) du = 
\pi(\alpha+\beta) \sum_j{\rm int}(u_{j\xi\upsilon},\Gamma_{\xi\upsilon})a_j
+ \pi(\alpha-\beta)\sum_j{\rm int}(u_{j\tau\zeta},\Gamma_{\tau\zeta})a_j.
\label{124}
\end{equation}

\subsection{Factorization of circular fourcomplex polynomials}

A polynomial of degree $m$ of the circular fourcomplex variable 
$u=x+\alpha y+\beta z+\gamma t$ has the form
\begin{equation}
P_m(u)=u^m+a_1 u^{m-1}+\cdots+a_{m-1} u +a_m ,
\label{125}
\end{equation}
where the constants are in general circular fourcomplex numbers.

It can be shown that any circular fourcomplex polynomial has a circular
fourcomplex root, whence 
it follows that a polynomial of degree $m$ can be written as a product of
$m$ linear factors of the form $u-u_j$, where the circular fourcomplex numbers
$u_j$ are 
the roots of the polynomials, although the factorization may not be unique, 
\begin{equation}
P_m(u)=\prod_{j=1}^m (u-u_j) .
\label{126}
\end{equation}

The fact that any circular fourcomplex polynomial has a root can be shown by
considering 
the transformation of a fourdimensional sphere with the center at the origin by
the function $u^m$. The points of the hypersphere of radius $d$ are of the form
written in Eq. (\ref{52}), with $d$ constant and $0\leq\phi<2\pi,
0\leq\chi<2\pi, 0\leq\psi\leq \pi/2$. The point $u^m$ is
\begin{equation}
u^m=d^m\exp\left(\alpha m\frac{\phi+\chi}{2}+\beta m\frac{\phi-\chi}{2}\right)
[\cos(\psi-\pi/4)+\gamma\sin(\psi-\pi/4)]^m .
\label{127}
\end{equation}
It can be shown with the aid of Eq. (\ref{56}) that
\begin{equation}
\left|u\exp
\left(\alpha \frac{\phi+\chi}{2}+\beta \frac{\phi-\chi}{2}\right)\right|
=|u|,
\label{127bb}
\end{equation}
so that
\begin{eqnarray}
\lefteqn{\left|[\cos(\psi-\pi/4)+\gamma\sin(\psi-\pi/4)]^m
\exp\left(\alpha m\frac{\phi+\chi}{2}+\beta m\frac{\phi-\chi}{2}\right)\right|
\nonumber}\\
&&=\left|[\cos(\psi-\pi/4)+\gamma\sin(\psi-\pi/4)]^m\right| .
\label{127b}
\end{eqnarray}
The right-hand side of Eq. (\ref{127b}) is
\begin{equation}
|(\cos\epsilon+\gamma\sin\epsilon)^m|^2
=\sum_{k=0}^m C_{2m}^{2k}\cos^{2m-2k}\epsilon\sin^{2k}\epsilon ,
\label{128}
\end{equation}
where $\epsilon=\psi-\pi/4$, 
and since $C_{2m}^{2k}\geq C_m^k$, it can be concluded that
\begin{equation}
|(\cos\epsilon+\gamma\sin\epsilon)^m|^2\geq 1 .
\label{129}
\end{equation}
Then
\begin{equation}
d^m\leq |u^m|\leq 2^{(m-1)/2} d^m ,
\label{129b}
\end{equation}
which shows that the image of a four-dimensional sphere via the transformation
operated by the function $u^m$ is a finite hypersurface.

If $u^\prime=u^m$, and
\begin{equation}
u^\prime=d^\prime
[\cos(\psi^\prime-\pi/4)+\gamma\sin(\psi^\prime-\pi/4)]
\exp\left(\alpha \frac{\phi^\prime+\chi^\prime}{2}+\beta
\frac{\phi^\prime-\chi^\prime}{2}\right),
\label{130}
\end{equation}
then 
\begin{equation}
\phi^\prime=m\phi, \: \chi^\prime=m\chi, \: \tan\psi^\prime=\tan^m\psi .
\label{131}
\end{equation}
Since for any values of the angles $\phi^\prime, \chi^\prime, \psi^\prime$
there is a set of solutions $\phi, \chi, \psi$ of Eqs. (\ref{131}), and since
the image of the hypersphere is a finite hypersurface, it follows that the
image of the four-dimensional sphere via the function $u^m$ is also a closed
hypersurface. A continuous hypersurface is called closed when any ray issued
from the 
origin intersects that surface at least once in the finite part of the space.

A transformation of the four-dimensional space by the polynomial $P_m(u)$
will be considered further. By this transformation, a hypersphere of radius $d$
having the center at the origin is changed into a certain finite closed
surface, as discussed previously. 
The transformation of the four-dimensional space by the polynomial $P_m(u)$
associates to the point $u=0$ the point $f(0)=a_m$, and the image of a
hypersphere of very large radius $d$ can be represented with good approximation
by the image of that hypersphere by the function $u^m$. 
The origin of the axes is an inner
point of the latter image. If the radius of the hypersphere is now reduced
continuously from the initial very large values to zero, the image hypersphere
encloses initially the origin, but the image shrinks to $a_m$ when the radius
approaches the value zero.  Thus, the
origin is initially inside the image hypersurface, and it lies outside the
image hypersurface when the radius of the hypersphere tends to zero. Then since
the image hypersurface is closed, the image surface must intersect at some
stage the origin of the axes, which means that there is a point $u_1$ such that
$f(u_1)=0$. The factorization in Eq. (\ref{126}) can then be obtained by
iterations.

The roots of the polynomial $P_m$ can be obtained by the following method.
If the constants in Eq. (\ref{125}) are $a_l=a_{l0}+\alpha a_{l1}
+\beta a_{l2}+\gamma a_{l3}$, and
with the 
notations of Eq. (\ref{n88b}), the polynomial $P_m(u)$ can be written as
\begin{eqnarray}
\lefteqn{P_m=\sum_{l=0}^{m} 2^{(m-l)/2}
(e_1 A_{l1}+\tilde e_1\tilde A_{l1})(e_1 \xi+\tilde e_1
\upsilon)^{m-l}\nonumber}\\ 
&&+\sum_{l=0}^{m} 2^{(m-l)/2}
(e_2 A_{l2}+\tilde e_2\tilde A_{l2})(e_2 \tau+\tilde e_2 \zeta)^{m-l} ,
\label{126a}
\end{eqnarray}
where the constants $A_{lk}, \tilde A_{lk}, k=1,2,$ are real numbers.
Each of the polynomials of degree $m$ in $e_1 \xi+\tilde e_1\upsilon, 
e_2 \tau+\tilde e_2\zeta$
in Eq. (\ref{126a}) 
can always be written as a product of linear factors of the form
$e_1 (\xi-\xi_p)+\tilde e_1(\upsilon- \upsilon_p)$ and respectively
$e_2 (\tau-\tau_p)+\tilde e_2(\zeta- \zeta_p)$, where the
constants $\xi_p, \upsilon_p, \tau_p, \zeta_p$ are real,
\begin{eqnarray}
\sum_{l=0}^{m} 2^{(m-l)/2}
(e_1 A_{l1}+\tilde e_1\tilde A_{l1})(e_1 \xi+\tilde e_1 \upsilon)^{m-l}
=\prod_{p=1}^{m}2^{m/2}\left\{e_1 (\xi-\xi_p)+\tilde e_1(\upsilon- \upsilon_p)
\right\},
\label{126bb}
\end{eqnarray}
\begin{eqnarray}
\sum_{l=0}^{m} 2^{(m-l)/2}
(e_2 A_{l2}+\tilde e_2\tilde A_{l2})(e_2 \tau+\tilde e_2 \zeta)^{m-l}
=\prod_{p=1}^{m}2^{m/2}\left\{e_2 (\tau-\tau_p)+\tilde e_2(\zeta- \zeta_p)
\right\}.
\label{126bc}
\end{eqnarray}

Due to the relations  (\ref{87d}),
the polynomial $P_m(u)$ can be written as a product of factors of
the form 
\begin{eqnarray}
P_m(u)=\prod_{p=1}^m 2^{m/2}\left\{e_1 (\xi-\xi_p)+\tilde e_1(\upsilon-
\upsilon_p) 
+e_2 (\tau-\tau_p)+\tilde e_2(\zeta- \zeta_p)\right\}.
\label{128b}
\end{eqnarray}
This relation can be written with the aid of Eq. (\ref{87b}) in the form 
(\ref{126}), where
\begin{eqnarray}
u_p=\sqrt{2}(e_1 \xi_p+\tilde e_1 \upsilon_p
+e_2 \tau_p+\tilde e_2 \zeta_p) .
\label{129bx}
\end{eqnarray}
The roots  
$e_1 \xi_p+\tilde e_1 \upsilon_p$ and $e_2 \tau_p+\tilde e_2 \zeta_p$
defined in Eqs. (\ref{126bb}) and respectively (\ref{126bc}) may be ordered
arbitrarily. This means that Eq. 
(\ref{129bx}) gives sets of $m$ roots 
$u_1,...,u_m$ of the polynomial $P_m(u)$, 
corresponding to the various ways in which the roots 
$e_1 \xi_p+\tilde e_1 \upsilon_p$ and $e_2 \tau_p+\tilde e_2 \zeta_p$
are ordered according to $p$ for each polynomial. 
Thus, while the hypercomplex components in Eqs. (\ref{126bb}), Eqs.
(\ref{126bc}) taken separately have unique factorizations, the polynomial
$P_m(u)$ can be written in many different ways as a product of linear factors. 
The result of the circular fourcomplex integration, Eq. (\ref{124}), is however
unique.  

If, for example, $P(u)=u^2+1$, the possible factorizations are
$P=(u-\tilde e_1-\tilde e_2)(u+\tilde e_1+\tilde e_2)$ and
$P=(u-\tilde e_1+\tilde e_2)(u+\tilde e_1-\tilde e_2)$, which can also be
written as $u^2+1=(u-\alpha)(u+\alpha)$ or as
$u^2+1=(u-\beta)(u+\beta)$. The result of the circular fourcomplex integration,
Eq. (\ref{124}), is however unique. 
It can be checked
that $(\pm \tilde e_1\pm\tilde e_2)^2=
-e_1-e_2=-1$.

\subsection{Representation of circular 
fourcomplex numbers by irreducible matrices}

If $T$ is the unitary matrix,
\begin{equation}
T =\left(
\begin{array}{cccc}
\frac{1}{\sqrt{2}}&0                 &0                 &\frac{1}{\sqrt{2}}\\
0                 &\frac{1}{\sqrt{2}}&\frac{1}{\sqrt{2}}&   0              \\
\frac{1}{\sqrt{2}}& 0                & 0                &-\frac{1}{\sqrt{2}} \\
0                 &\frac{1}{\sqrt{2}}&-\frac{1}{\sqrt{2}}&   0               \\
\end{array}
\right),
\label{129x}
\end{equation}
it can be shown 
that the matrix $T U T^{-1}$ has the form 
\begin{equation}
T U T^{-1}=\left(
\begin{array}{cc}
V_1      &     0    \\
0        &     V_2  \\
\end{array}
\right),
\label{129y}
\end{equation}
where $U$ is the matrix in Eq. (\ref{23}) used to represent the circular
fourcomplex 
number $u$. In Eq. (\ref{129y}), $V_1, V_2$ are
the matrices
\begin{equation}
V_1=\left(
\begin{array}{cc}
x+t    &   y+z   \\
-y-z   &   x+t   \\
\end{array}\right),\;\;
V_2=\left(
\begin{array}{cc}
x-t    &   y-z   \\
-y+z   &   x-t   \\
\end{array}\right).
\label{130x}
\end{equation}
In Eq. (\ref{129y}), the symbols 0 denote the matrix
\begin{equation}
\left(
\begin{array}{cc}
0   &  0   \\
0   &  0   \\
\end{array}\right).
\label{131x}
\end{equation}
The relations between the variables $x+t,y+z$,$x-t,y-z$ for the multiplication
of circular fourcomplex numbers have been written in Eqs.
(\ref{17})-(\ref{20}). The 
matrix 
$T U T^{-1}$ provides an irreducible representation
\cite{4} of the circular fourcomplex number $u$ in terms of matrices with real
coefficients.

\section{Hyperbolic Complex Numbers in Four Dimensions}

\subsection{Operations with hyperbolic fourcomplex numbers}

A hyperbolic fourcomplex number is determined by its four components
$(x,y,z,t)$. The sum 
of the hyperbolic fourcomplex numbers $(x,y,z,t)$ and
$(x^\prime,y^\prime,z^\prime,t^\prime)$ is the hyperbolic fourcomplex
number $(x+x^\prime,y+y^\prime,z+z^\prime,t+t^\prime)$. 
The product of the hyperbolic fourcomplex numbers
$(x,y,z,t)$ and $(x^\prime,y^\prime,z^\prime,t^\prime)$ 
is defined in this work to be the hyperbolic fourcomplex
number
$(xx^\prime+yy^\prime+zz^\prime+tt^\prime,
xy^\prime+yx^\prime+zt^\prime+tz^\prime,
xz^\prime+zx^\prime+yt^\prime+ty^\prime,
xt^\prime+tx^\prime+yz^\prime+zy^\prime)$.

Hyperbolic fourcomplex numbers and their operations can be represented by
writing the 
hyperbolic fourcomplex number $(x,y,z,t)$ as  
$u=x+\alpha y+\beta z+\gamma t$, where $\alpha, \beta$ and $\gamma$ 
are bases for which the multiplication rules are 
\begin{equation}
\alpha^2=1, \:\beta^2=1, \:\gamma^2=1, \alpha\beta=\beta\alpha=\gamma,\:
\alpha\gamma=\gamma\alpha=\beta, \:\beta\gamma=\gamma\beta=\alpha .
\label{h1}
\end{equation}
Two hyperbolic fourcomplex numbers $u=x+\alpha y+\beta z+\gamma t, 
u^\prime=x^\prime+\alpha y^\prime+\beta z^\prime+\gamma t^\prime$ are equal, 
$u=u^\prime$, if and only if $x=x^\prime, y=y^\prime,
z=z^\prime, t=t^\prime$. 
If 
$u=x+\alpha y+\beta z+\gamma t, 
u^\prime=x^\prime+\alpha y^\prime+\beta z^\prime+\gamma t^\prime$
are hyperbolic fourcomplex numbers, 
the sum $u+u^\prime$ and the 
product $uu^\prime$ defined above can be obtained by applying the usual
algebraic rules to the sum 
$(x+\alpha y+\beta z+\gamma t)+ 
(x^\prime+\alpha y^\prime+\beta z^\prime+\gamma t^\prime)$
and to the product 
$(x+\alpha y+\beta z+\gamma t)
(x^\prime+\alpha y^\prime+\beta z^\prime+\gamma t^\prime)$,
and grouping of the resulting terms,
\begin{equation}
u+u^\prime=x+x^\prime+\alpha(y+y^\prime)+\beta(z+z^\prime)+\gamma(t+t^\prime),
\label{h1a}
\end{equation}
\begin{equation}
uu^\prime=
xx^\prime+yy^\prime+zz^\prime+tt^\prime+
\alpha(xy^\prime+yx^\prime+zt^\prime+tz^\prime)+
\beta(xz^\prime+zx^\prime+yt^\prime+ty^\prime)+
\gamma(xt^\prime+tx^\prime+yz^\prime+zy^\prime)
\label{h1b}
\end{equation}

If $u,u^\prime,u^{\prime\prime}$ are hyperbolic fourcomplex numbers, the
multiplication is associative 
\begin{equation}
(uu^\prime)u^{\prime\prime}=u(u^\prime u^{\prime\prime})
\label{h2}
\end{equation}
and commutative
\begin{equation}
u u^\prime=u^\prime u ,
\label{h3}
\end{equation}
as can be checked through direct calculation.
The hyperbolic fourcomplex zero is $0+\alpha\cdot 0+\beta\cdot 0+\gamma\cdot
0,$  
denoted simply 0, 
and the hyperbolic fourcomplex unity is $1+\alpha\cdot 0+\beta\cdot
0+\gamma\cdot 0,$  
denoted simply 1.

The inverse of the hyperbolic fourcomplex number 
$u=x+\alpha y+\beta z+\gamma t$ is a hyperbolic fourcomplex number
$u^\prime=x^\prime+\alpha y^\prime+\beta z^\prime+\gamma t^\prime$
having the property that
\begin{equation}
uu^\prime=1 .
\label{h4}
\end{equation}
Written on components, the condition, Eq. (\ref{h4}), is
\begin{equation}
\begin{array}{c}
xx^\prime+yy^\prime+zz^\prime+tt^\prime=1,\\
yx^\prime+xy^\prime+tz^\prime+zt^\prime=0,\\
zx^\prime+ty^\prime+xz^\prime+yt^\prime=0,\\
tx^\prime+zy^\prime+yz^\prime+xt^\prime=0 .
\end{array}
\label{h5}
\end{equation}
The system (\ref{h5}) has the solution
\begin{equation}
x^\prime=\frac{x(x^2-y^2-z^2-t^2)+2yzt}{\nu} ,
\label{h6a}
\end{equation}
\begin{equation}
y^\prime=
\frac{y(-x^2+y^2-z^2-t^2)+2xzt}{\nu} ,
\label{h6b}
\end{equation}
\begin{equation}
z^\prime=\frac{z(-x^2-y^2+z^2-t^2)+2xyt}{\nu} ,
\label{h6c}
\end{equation}
\begin{equation}
t^\prime=\frac{t(-x^2-y^2-z^2+t^2)+2xyz}{\nu} ,
\label{h6d}
\end{equation}
provided that $\nu\not= 0$, where
\begin{equation}
\nu=x^4+y^4+z^4+t^4-2(x^2y^2+x^2z^2+x^2t^2+y^2z^2+y^2t^2+z^2t^2)+8xyzt .
\label{h6e}
\end{equation}

The quantity $\nu$ can be written as
\begin{equation}
\nu=ss^\prime s^{\prime\prime}s^{\prime\prime\prime} ,
\label{h7}
\end{equation}
where
\begin{equation}
s=x+y+z+t, \: s^\prime= x-y+z-t , \: s^{\prime\prime}=x+y-z-t,  \:
s^{\prime\prime\prime}=x-y-z+t .
\label{h8}
\end{equation}
The variables $s, s^\prime, s^{\prime\prime}, s^{\prime\prime\prime}$ will be
called canonical  
hyperbolic fourcomplex variables.
Then a hyperbolic fourcomplex number $u=x+\alpha y+\beta z+\gamma t$ has an
inverse, 
unless 
\begin{equation}
s=0 ,\:\:{\rm or}\:\: s^\prime=0, \:\:{\rm or}\:\: s^{\prime\prime}=0, \:\:
{\rm or}\:\: s^{\prime\prime\prime}=0 . 
\label{h9}
\end{equation}

For arbitrary values of the variables $x,y,z,t$, the quantity $\nu$ can be
positive or negative. If $\nu\geq 0$, the quantity $\mu=\nu^{1/4}$
will be called amplitude of the hyperbolic fourcomplex number
$x+\alpha y+\beta z +\gamma t$.
The normals of the hyperplanes in Eq. (\ref{h9}) are orthogonal to
each other. Because of conditions (\ref{h9}) these hyperplanes will
be also called the nodal hyperplanes. 
It can be shown that 
if $uu^\prime=0$ then either $u=0$, or $u^\prime=0$, or $q, q^\prime$ belong to
pairs of orthogonal hypersurfaces as described further.
Thus, divisors of zero exist if
one of the hyperbolic fourcomplex
numbers $u, u^\prime$ belongs to one of the nodal hyperplanes 
and the other hyperbolic fourcomplex 
number belongs to the straight line through the origin
which is normal to that hyperplane,
\begin{equation}
x+y+z+t=0,\:\:{\rm and} \:\: x^\prime=y^\prime=z^\prime=t^\prime,  
\label{h10a}
\end{equation}
or
\begin{equation}
x-y+z-t=0, \:\:{\rm and}\:\: x^\prime=-y^\prime=z^\prime=-t^\prime, 
\label{h10b}
\end{equation}
or
\begin{equation}
x+y-z-t=0, \:\:{\rm and}\:\: x^\prime=y^\prime=-z^\prime=-t^\prime, 
\label{h10c}
\end{equation}
or
\begin{equation}
x-y-z+t=0, \:\:{\rm and}\:\: x^\prime=-y^\prime=-z^\prime=t^\prime.
\label{h10d}
\end{equation}
Divisors of zero also exist if the hyperbolic fourcomplex numbers $u,u^\prime$
belong to 
different members of the pairs 
of two-dimensional hypersurfaces listed further,
\begin{equation}
x+y=0,\:z+t=0\:\: {\rm and}\:\:  x^\prime-y^\prime=0 , \: z^\prime-t^\prime=0,
\label{h11a}
\end{equation}
or
\begin{equation}
x+z=0,\:y+t=0\:\: {\rm and}\:\:  x^\prime-z^\prime=0 , \: y^\prime-t^\prime=0,
\label{h11b}
\end{equation}
or
\begin{equation}
y+z=0,\:x+t=0\:\: {\rm and}\:\:  y^\prime-z^\prime=0 , \: x^\prime-t^\prime=0.
\label{h11c}
\end{equation}

\subsection{Geometric representation of hyperbolic fourcomplex numbers}

The hyperbolic fourcomplex number $x+\alpha y+\beta z+\gamma t$ can be
represented by the point $A$ of coordinates $(x,y,z,t)$. 
If $O$ is the origin of the four-dimensional space $x,y,z,t,$ the distance 
from $A$ to the origin $O$ can be taken as
\begin{equation}
d^2=x^2+y^2+z^2+t^2 .
\label{h12}
\end{equation}
The distance $d$ will be called modulus of the hyperbolic fourcomplex number
$x+\alpha 
y+\beta z +\gamma t$, $d=|u|$. 

If $u=x+\alpha y+\beta z +\gamma t, u_1=x_1+\alpha y_1+\beta z_1 +\gamma t_1,
u_2=x_2+\alpha y_2+\beta z_2 +\gamma t_2$, and $u=u_1u_2$, and if
\begin{equation}
s_j=x_j+y_j+z_j+t_j, \: s_j^\prime= x_j-y_j+z_j-t_j , \:
s_j^{\prime\prime}=x_j+y_j-z_j-t_j,  \: 
s_j^{\prime\prime\prime}=x_j-y_j-z_j+t_j , j=1,2,
\label{h13}
\end{equation}
it can be shown that
\begin{equation}
s=s_1s_2 ,\:\:
s^\prime=s_1^\prime s_2^\prime, \:\:
s^{\prime\prime}=s_1^{\prime\prime}s_2^{\prime\prime}, \:\:
s^{\prime\prime\prime}=s_1^{\prime\prime\prime}s_2^{\prime\prime\prime} . 
\label{h14}
\end{equation}
The relations (\ref{h14}) are a consequence of the identities
\begin{eqnarray}
\lefteqn{(x_1x_2+y_1y_2+z_1z_2+t_1t_2)+(x_1y_2+y_1x_2+z_1t_2+t_1z_2)
\nonumber}\\
&&+(x_1z_2+z_1x_2+y_1t_2+t_1y_2)+(x_1t_2+t_1x_2+y_1z_2+z_1y_2)\nonumber\\
&&=(x_1+y_1+z_1+t_1)(x_2+y_2+z_2+t_2)
\label{h15}
\end{eqnarray}
\begin{eqnarray}
\lefteqn{(x_1x_2+y_1y_2+z_1z_2+t_1t_2)-(x_1y_2+y_1x_2+z_1t_2+t_1z_2)
\nonumber}\\
&&+(x_1z_2+z_1x_2+y_1t_2+t_1y_2)-(x_1t_2+t_1x_2+y_1z_2+z_1y_2)\nonumber\\
&&=(x_1-y_1+z_1-t_1)(x_2-y_2+z_2-t_2)
\label{h16}
\end{eqnarray}
\begin{eqnarray}
\lefteqn{(x_1x_2+y_1y_2+z_1z_2+t_1t_2)+(x_1y_2+y_1x_2+z_1t_2+t_1z_2)
\nonumber}\\
&&-(x_1z_2+z_1x_2+y_1t_2+t_1y_2)-(x_1t_2+t_1x_2+y_1z_2+z_1y_2)\nonumber\\
&&=(x_1+y_1-z_1-t_1)(x_2+y_2-z_2-t_2)
\label{h17}
\end{eqnarray}
\begin{eqnarray}
\lefteqn{(x_1x_2+y_1y_2+z_1z_2+t_1t_2)-(x_1y_2+y_1x_2+z_1t_2+t_1z_2)
\nonumber}\\
&&-(x_1z_2+z_1x_2+y_1t_2+t_1y_2)+(x_1t_2+t_1x_2+y_1z_2+z_1y_2)\nonumber\\
&&=(x_1-y_1-z_1+t_1)(x_2-y_2-z_2+t_2)
\label{h18}
\end{eqnarray}
A consequence of the relations (\ref{h14}) is that if $u=u_1u_2$, then
\begin{equation}
\nu=\nu_1\nu_2 ,
\label{h19}
\end{equation}
where
\begin{equation}
\nu_j=s_j s_j^\prime s_j^{\prime\prime} s_j^{\prime\prime\prime} , \:\:j=1,2.
\label{h20}
\end{equation}

The hyperbolic fourcomplex numbers
\begin{equation}
e=\frac{1+\alpha+\beta+\gamma}{4},\:
e^\prime=\frac{1-\alpha+\beta-\gamma}{4},\:
e^{\prime\prime}=\frac{1+\alpha-\beta-\gamma}{4},\:
e^{\prime\prime\prime}=\frac{1-\alpha-\beta+\gamma}{4}
\label{h21}
\end{equation}
are orthogonal,
\begin{equation}
ee^\prime=0,\:ee^{\prime\prime}=0,\:ee^{\prime\prime\prime}=0,\:e^\prime
e^{\prime\prime}=0,\:e^\prime
e^{\prime\prime\prime}=0,\:e^{\prime\prime}e^{\prime\prime\prime}=0,  
\label{h22}
\end{equation}
and have also the property that
\begin{equation}
e^2=e, \: e^{\prime 2}=e^{\prime}, \:
e^{\prime\prime 2}=e^{\prime\prime}, \:e^{\prime\prime\prime
2}=e^{\prime\prime\prime} . 
\label{h23}
\end{equation}
The hyperbolic fourcomplex number $u=x+\alpha y+\beta z+\gamma t$ can be
written as 
\begin{equation}
x+\alpha y+\beta z+\gamma t
=(x+y+z+t)e+(x-y+z-t)e^\prime+(x+y-z-t)e^{\prime\prime}
+(x-y-z+t)e^{\prime\prime\prime},
\label{h24}
\end{equation}
or, by using Eq. (\ref{h8}),
\begin{equation}
u=se+s^\prime
e^\prime+s^{\prime\prime}e^{\prime\prime}
+s^{\prime\prime\prime}e^{\prime\prime\prime}. 
\label{h25}
\end{equation}
The ensemble $e, e^\prime, e^{\prime\prime}, e^{\prime\prime\prime}$ will be
called the canonical 
hyperbolic fourcomplex base, and Eq. (\ref{h25}) gives the canonical form of
the 
hyperbolic fourcomplex number.
Thus, if $u_j=s_je+s_j^\prime e^\prime+s_j^{\prime\prime}e^{\prime\prime}
+s_j^{\prime\prime\prime}e^{\prime\prime\prime}, \:j=1,2$, and $u=u_1u_2$, then
the multiplication of the hyperbolic fourcomplex numbers is expressed by the
relations (\ref{h14}).
The moduli of the bases $e, e^\prime, e^{\prime\prime}, e^{\prime\prime\prime}$
are
\begin{equation}
|e|=\frac{1}{2},\; |e^\prime|=\frac{1}{2},\; 
|e^{\prime\prime}|=\frac{1}{2},\; |e^{\prime\prime\prime}|=\frac{1}{2}.
\label{h25b}
\end{equation}
The distance $d$, Eq. (\ref{h12}), is given by
\begin{equation}
d^2=\frac{1}{4}\left(s^2+s^{\prime 2}
+s^{\prime\prime 2}+s^{\prime\prime\prime 2}\right). 
\label{h25c}
\end{equation}

The relation (\ref{h25c}) shows that the variables $s, s^{\prime},
s^{\prime\prime}, s^{\prime\prime\prime}$ can be written as 
\begin{equation}
s=2d\cos\psi\cos\phi,\; s^\prime=2d\cos\psi\sin\phi,\; 
s^{\prime\prime}=2d\sin\psi\cos\chi,\; s=2d\sin\psi\sin\chi,
\label{h25d}
\end{equation}
where $\phi$ is the azimuthal angle in the $s,s^\prime$ plane,
$0\leq\phi<2\pi$, 
$\chi$ is the azimuthal angle in the $s^{\prime\prime},
s^{\prime\prime\prime}$ plane, $0\leq\chi<2\pi$,
and $\psi$ is the angle between the line $OA$ and the plane $ss^\prime$,
$0\leq\psi\leq\pi/2$.
The variables $x,y,z,t$ can be expressed in terms of the distance $d$ and the
angles $\phi, \chi, \psi$ as
\begin{equation}
\begin{array}{c}
x=\frac{d}{2}(\cos\psi\cos\phi+\cos\psi\sin\phi 
+\sin\psi\cos\chi+\sin\psi\sin\chi),\\
y=\frac{d}{2}(\cos\psi\cos\phi-\cos\psi\sin\phi 
+\sin\psi\cos\chi-\sin\psi\sin\chi),\\
z=\frac{d}{2}(\cos\psi\cos\phi+\cos\psi\sin\phi 
-\sin\psi\cos\chi-\sin\psi\sin\chi),\\
t=\frac{d}{2}(\cos\psi\cos\phi-\cos\psi\sin\phi 
-\sin\psi\cos\chi+\sin\psi\sin\chi).
\end{array}
\label{h25e}
\end{equation}

If $u=u_1u_2$, and the hypercomplex numbers $u_1, u_2$ are described by the
variables $d_1, \phi_1, \chi_1, \psi_1$ and respectively $d_2, \phi_2, \chi_2,
\psi_2$, then from Eq. (\ref{h25d}) it results that
\begin{eqnarray}
\lefteqn{\tan\phi=\tan\phi_1\tan\phi_2,\;
\tan\chi=\tan\chi_1\tan\chi_2,\nonumber}\\
&&\frac{\tan^2\psi\sin 2\chi}{\sin 2\phi}=
\frac{\tan^2\psi_1\sin 2\chi_1}{\sin 2\phi_1}
\frac{\tan^2\psi_2\sin 2\chi_2}{\sin 2\phi_2}.
\label{h25f}
\end{eqnarray}

The relation (\ref{h19}) for the product of hyperbolic fourcomplex numbers can 
be demonstrated also by using a representation of the multiplication of the 
hyperbolic fourcomplex numbers by matrices, in which the hyperbolic fourcomplex
number $u=x+\alpha 
y+\beta z+\gamma t$ is represented by the matrix
\begin{equation}
A=\left(\begin{array}{cccc}
x&y&z&t\\
y&x&t&z\\
z&t&x&y\\
t&z&y&x 
\end{array}\right) .
\label{h26}
\end{equation}
The product $u=x+\alpha y+\beta z+\gamma t$ of the hyperbolic fourcomplex
numbers 
$u_1=x_1+\alpha y_1+\beta z_1+\gamma t_1, u_2=x_2+\alpha y_2+\beta z_2+\gamma
t_2$, can be represented by the matrix multiplication 
\begin{equation}
A=A_1A_2.
\label{h27}
\end{equation}
It can be checked that the determinant ${\rm det}(A)$ of the matrix $A$ is
\begin{equation}
{\rm det}A = \nu .
\label{h28}
\end{equation}
The identity (\ref{h19}) is then a consequence of the fact the determinant 
of the product of matrices is equal to the product of the determinants 
of the factor matrices.

\subsection{Exponential form of a hyperbolic fourcomplex number}

The exponential function of a fourcomplex variable $u$ and the addition
theorem for the exponential function have been written in Eqs. (\ref{g1}) and
(\ref{g2}).
If $u=x+\alpha y+\beta z+\gamma t$, then  $\exp u$ can be calculated as
$\exp u=\exp x \cdot \exp (\alpha y) \cdot \exp (\beta z) \cdot \exp (\gamma
t)$. According to Eqs. (\ref{h1}), 
\begin{equation}
\alpha^{2m}=1, \alpha^{2m+1}=\alpha, 
\beta^{2m}=1, \beta^{2m+1}=\beta, 
\gamma^{2m}=1, \gamma^{2m+1}=\gamma, 
\label{h31}
\end{equation}
where $m$ is a natural number,
so that $\exp (\alpha y), \: \exp(\beta z)$ and $\exp(\gamma t)$ can be
written as 
\begin{equation}
\exp (\alpha y) = \cosh y +\alpha \sinh y , \:
\exp (\beta z) = \cosh y +\beta \sinh z , \:
\exp (\gamma t) = \cosh t +\gamma \sinh t . 
\label{h32}
\end{equation}
From Eqs. (\ref{h32}) it can be inferred that
\begin{eqnarray}
\lefteqn{(\cosh t +\alpha \sinh t)^m=\cosh mt +\alpha \sinh mt ,\:
(\cosh t +\beta \sinh t)^m=\cosh mt +\beta \sinh mt ,\nonumber}\\
&&(\cosh t +\gamma \sinh t)^m=\cosh mt +\gamma \sinh mt . 
\label{h33}
\end{eqnarray}

The hyperbolic fourcomplex numbers $u=x+\alpha y+\beta z+\gamma t$ for which
$s=x+y+z+t>0, \: s^\prime= x-y+z-t>0 , \: s^{\prime\prime}=x+y-z-t>0,  \:
s^{\prime\prime\prime}=x-y-z+t>0$ can be written in the form 
\begin{equation}
x+\alpha y+\beta z+\gamma t=e^{x_1+\alpha y_1+\beta z_1+\gamma t_1} .
\label{h34}
\end{equation}
The conditions $s=x+y+z+t>0, \: s^\prime= x-y+z-t>0 , \:
s^{\prime\prime}=x+y-z-t>0,  \: s^{\prime\prime\prime}=x-y-z+t>0$ 
correspond in Eq. (\ref{h25d}) to a range of angles $0<\phi<\pi/2,
0<\chi<\pi/2, 0<\psi\leq\pi/2$. 
The expressions of $x_1, y_1, z_1, t_1$ as functions of 
$x, y, z, t$ can be obtained by
developing $e^{\alpha y_1}, e^{\beta z_1}$and $e^{\gamma t_1}$ with the aid of
Eqs. (\ref{h32}), by multiplying these expressions and separating
the hypercomplex components, 
\begin{equation}
x=e^{x_1}(\cosh y_1\cosh z_1\cosh t_1+\sinh y_1\sinh z_1\sinh t_1) ,
\label{h35}
\end{equation}
\begin{equation}
y=e^{x_1}(\sinh y_1\cosh z_1\cosh t_1+\cosh y_1\sinh z_1\sinh t_1) ,
\label{h36}
\end{equation}
\begin{equation}
z=e^{x_1}(\cosh y_1\sinh z_1\cosh t_1+\sinh y_1\cosh z_1\sinh t_1) ,
\label{h37}
\end{equation}
\begin{equation}
t=e^{x_1}(\sinh y_1\sinh z_1\cosh t_1+\cosh y_1\cosh z_1\sinh t_1) ,
\label{h38}
\end{equation}
It can be shown from Eqs. (\ref{h35})-(\ref{h38}) that
\begin{equation}
x_1=\frac{1}{4} \ln(s s^\prime s^{\prime\prime}s^{\prime\prime\prime}) , \:
y_1=\frac{1}{4}\ln\frac{ss^{\prime\prime}}{s^\prime s^{\prime\prime\prime}},\:
z_1=\frac{1}{4}\ln\frac{ss^\prime}{s^{\prime\prime}s^{\prime\prime\prime}},\:
t_1=\frac{1}{4}\ln\frac{ss^{\prime\prime\prime}}{s^\prime s^{\prime\prime}}.
\label{h39}
\end{equation}
The exponential form of the hyperbolic fourcomplex number $u$ can be written
as
\begin{equation}
u=\mu\exp\left(
\frac{1}{4}\alpha \ln\frac{ss^{\prime\prime}}{s^\prime s^{\prime\prime\prime}}
+\frac{1}{4}\beta \ln\frac{ss^\prime}{s^{\prime\prime}s^{\prime\prime\prime}}
+\frac{1}{4}\gamma \ln\frac{ss^{\prime\prime\prime}}{s^\prime s^{\prime\prime}}
\right),
\label{h40}
\end{equation}
where
\begin{equation}
\mu=(s s^\prime s^{\prime\prime}s^{\prime\prime\prime})^{1/4}.
\label{h40b}
\end{equation}
The exponential form of the hyperbolic fourcomplex number $u$ can be written
with the aid of the relations (\ref{h25d}) as
\begin{equation}
u=\mu\exp\left(\frac{1}{4}\alpha \ln\frac{1}{\tan\phi\tan\chi}
+\frac{1}{4}\beta\ln\frac{\sin 2\phi}{\tan^2\psi \sin 2\chi}
+\frac{1}{4}\gamma \ln\frac{\tan \chi}{\tan\phi}\right).
\label{h40x}
\end{equation}
The amplitude $\mu$ can be expressed in terms of the distance $d$ with the aid
of Eqs. (\ref{h25d}) as
\begin{equation}
\mu=d\sin^{1/2}2\psi\sin^{1/4}2\phi\sin^{1/4}2\chi.
\label{h40c}
\end{equation}
The hypercomplex number can be written as
\begin{equation}
u=d\sin^{1/2}2\psi\sin^{1/4}2\phi\sin^{1/4}2\chi
\exp\left(\frac{1}{4}\alpha \ln\frac{1}{\tan\phi\tan\chi}
+\frac{1}{4}\beta\ln\frac{\sin 2\phi}{\tan^2\psi \sin 2\chi}
+\frac{1}{4}\gamma \ln\frac{\tan \chi}{\tan\phi}\right),
\label{h40d}
\end{equation}
which is the trigonometric form of the hypercomplex number $u$.

\subsection{Elementary functions of a hyperbolic fourcomplex variable}

The logarithm $u_1$ of the hyperbolic fourcomplex number $u$, $u_1=\ln u$, can
be defined 
for $s>0, s^\prime>0, s^{\prime\prime}>0, s^{\prime\prime\prime}>0$ 
as the solution of the equation
\begin{equation}
u=e^{u_1} ,
\label{h41}
\end{equation}
for $u_1$ as a function of
$u$. From Eq. (\ref{h40}) it results that 
\begin{equation}
\ln u=\frac{1}{4}\ln\mu+
+\frac{1}{4}\alpha \ln\frac{ss^{\prime\prime}}{s^\prime s^{\prime\prime\prime}}
+\frac{1}{4}\beta \ln\frac{ss^\prime}{s^{\prime\prime}s^{\prime\prime\prime}}
+\frac{1}{4}\gamma \ln\frac{ss^{\prime\prime\prime}}{s^\prime
s^{\prime\prime}}. 
\label{h42}
\end{equation}
Using the expression in Eq. (\ref{h40x}), the logarithm can be written as
\begin{equation}
\ln u=\frac{1}{4}\ln\mu
+\frac{1}{4}\alpha \ln\frac{1}{\tan\phi\tan\chi}
+\frac{1}{4}\beta\ln\frac{\sin 2\phi}{\tan^2\psi \sin 2\chi}
+\frac{1}{4}\gamma \ln\frac{\tan \chi}{\tan\phi}.
\label{h42b}
\end{equation}

It can be inferred from Eqs. (\ref{h42}) and (\ref{h14}) that
\begin{equation}
\ln(u_1u_2)=\ln u_1+\ln u_2 .
\label{h43}
\end{equation}
The explicit form of Eq. (\ref{h42}) is 
\begin{eqnarray}
\lefteqn{\ln (x+\alpha y+\beta z+\gamma t)=
\frac{1}{4}(1+\alpha+\beta+\gamma)\ln(x+y+z+t)\nonumber}\\
&&+\frac{1}{4}(1-\alpha+\beta-\gamma)\ln(x-y+z-t)
+\frac{1}{4}(1+\alpha-\beta-\gamma)\ln(x+y-z-t)\nonumber\\
&&+\frac{1}{4}(1-\alpha-\beta+\gamma)\ln(x-y-z+t) .
\label{h45}
\end{eqnarray}
The relation (\ref{h45}) can be written with the aid of Eq. (\ref{h21}) as
\begin{equation}
\ln u = e\ln s + e^\prime\ln s^\prime+e^{\prime\prime}\ln s^{\prime\prime}
+e^{\prime\prime\prime}\ln s^{\prime\prime\prime}.
\label{h44}
\end{equation}

The power function $u^n$ can be defined for $s>0, s^\prime>0,
s^{\prime\prime}>0, 
s^{\prime\prime\prime}>0$ and real values of $n$ as
\begin{equation}
u^n=e^{n\ln u} .
\label{h46}
\end{equation}
It can be inferred from Eqs. (\ref{h46}) and (\ref{h43}) that
\begin{equation}
(u_1u_2)^n=u_1^n\:u_2^n .
\label{h47}
\end{equation}
Using the expression (\ref{h45}) for $\ln u$ and the relations (\ref{h22}) and
(\ref{h23}) it can be shown that
\begin{eqnarray}
\lefteqn{(x+\alpha y+\beta z+\gamma t)^n=
\frac{1}{4}(1+\alpha+\beta+\gamma)(x+y+z+t)^n
+\frac{1}{4}(1-\alpha+\beta-\gamma)(x-y+z-t)^n\nonumber}\\
&&+\frac{1}{4}(1+\alpha-\beta-\gamma)(x+y-z-t)^n
+\frac{1}{4}(1-\alpha-\beta+\gamma)(x-y-z+t)^n .
\label{h48}
\end{eqnarray}
For integer $n$, the relation (\ref{h48}) is valid for any $x,y,z,t$. The
relation (\ref{h48}) for $n=-1$ is
\begin{equation}
\frac{1}{x+\alpha y+\beta z+\gamma t}=
\frac{1}{4}\left(\frac{1+\alpha+\beta+\gamma}{x+y+z+t}
+\frac{1-\alpha+\beta-\gamma}{x-y+z-t}
+\frac{1+\alpha-\beta-\gamma}{x+y-z-t}
+\frac{1-\alpha-\beta+\gamma}{x-y-z+t}\right) .
\label{h49}
\end{equation}

The trigonometric functions of the fourcomplex variable
$u$ and the addition theorems for these functions have been written in Eqs.
(\ref{g3})-(\ref{g6}). 
The cosine and sine functions of the hypercomplex variables $\alpha y, 
\beta z$ and $ \gamma t$ can be expressed as
\begin{equation}
\cos\alpha y=\cos y, \: \sin\alpha y=\alpha\sin y, 
\label{h54}
\end{equation}
\begin{equation}
\cos\beta y=\cos y, \: \sin\beta y=\beta\sin y, 
\label{h55}
\end{equation}
\begin{equation}
\cos\gamma y=\cos y, \: \sin\gamma y=\gamma\sin y .
\label{h56}
\end{equation}
The cosine and sine functions of a hyperbolic fourcomplex number $x+\alpha
y+\beta 
z+\gamma t$ can then be
expressed in terms of elementary functions with the aid of the addition
theorems Eqs. (\ref{g5}), (\ref{g6}) and of the expressions in  Eqs. 
(\ref{h54})-(\ref{h56}). 

The hyperbolic functions of the fourcomplex variable
$u$ and the addition theorems for these functions have been written in Eqs.
(\ref{g7})-(\ref{g10}). 
The hyperbolic cosine and sine functions of the hypercomplex variables $\alpha
y,  
\beta z$ and $ \gamma t$ can be expressed as
\begin{equation}
\cosh\alpha y=\cosh y, \: \sinh\alpha y=\alpha\sinh y, 
\label{h61}
\end{equation}
\begin{equation}
\cosh\beta y=\cosh y, \: \sinh\beta y=\beta\sinh y, 
\label{h62}
\end{equation}
\begin{equation}
\cosh\gamma y=\cosh y, \: \sinh\gamma y=\gamma\sinh y .
\label{h63}
\end{equation}
The hyperbolic cosine and sine functions of a hyperbolic fourcomplex number
$x+\alpha y+\beta 
z+\gamma t$ can then be
expressed in terms of elementary functions with the aid of the addition
theorems Eqs. (\ref{g9}), (\ref{g10}) and of the expressions in  Eqs. 
(\ref{h61})-(\ref{h63}).

\subsection{Power series of hyperbolic fourcomplex variables}

A hyperbolic fourcomplex series is an infinite sum of the form
\begin{equation}
a_0+a_1+a_2+\cdots+a_l+\cdots , 
\label{h64}
\end{equation}
where the coefficients $a_l$ are hyperbolic fourcomplex numbers. The
convergence of  
the series (\ref{h64}) can be defined in terms of the convergence of its 4 real
components. The convergence of a hyperbolic fourcomplex series can however be
studied 
using hyperbolic fourcomplex variables. The main criterion for absolute
convergence  
remains the comparison theorem, but this requires a number of inequalities
which will be discussed further.

The modulus of a hyperbolic fourcomplex number $u=x+\alpha y+\beta z+\gamma t$
can be defined as 
\begin{equation}
|u|=(x^2+y^2+z^2+t^2)^{1/2} ,
\label{h65}
\end{equation}
so that according to Eq. (\ref{h12}) $d=|u|$. Since $|x|\leq |u|, |y|\leq |u|,
|z|\leq |u|, |t|\leq |u|$, a property of 
absolute convergence established via a comparison theorem based on the modulus
of the series (\ref{h64}) will ensure the absolute convergence of each real
component of that series.

The modulus of the sum $u_1+u_2$ of the hyperbolic fourcomplex numbers $u_1,
u_2$ fulfils 
the inequality
\begin{equation}
||u_1|-|u_2||\leq |u_1+u_2|\leq |u_1|+|u_2| .
\label{h66}
\end{equation}
For the product the relation is 
\begin{equation}
|u_1u_2|\leq 2|u_1||u_2| ,
\label{h67}
\end{equation}
which replaces the relation of equality extant for regular complex numbers.
The equality in Eq. (\ref{h67}) takes place for $x_1^2=y_1^2=z_1^2=t_1^2$ and
$x_2/x_1=y_2/y_1=z_2/z_1=t_2/t_1$.
In particular
\begin{equation}
|u^2|\leq 2(x^2+y^2+z^2+t^2) .
\label{h68}
\end{equation}
The inequality in Eq. (\ref{h67}) implies that
\begin{equation}
|u^l|\leq 2^{l-1}|u|^l .
\label{h69}
\end{equation}
From Eqs. (\ref{h67}) and (\ref{h69}) it results that
\begin{equation}
|au^l|\leq 2^l |a| |u|^l .
\label{h70}
\end{equation}

A power series of the hyperbolic fourcomplex variable $u$ is a series of the
form 
\begin{equation}
a_0+a_1 u + a_2 u^2+\cdots +a_l u^l+\cdots .
\label{h71}
\end{equation}
Since
\begin{equation}
\left|\sum_{l=0}^\infty a_l u^l\right| \leq  \sum_{l=0}^\infty
2^l|a_l| |u|^l ,
\label{h72}
\end{equation}
a sufficient condition for the absolute convergence of this series is that
\begin{equation}
\lim_{l\rightarrow \infty}\frac{2|a_{l+1}||u|}{|a_l|}<1 .
\label{h73}
\end{equation}
Thus the series is absolutely convergent for 
\begin{equation}
|u|<c_0,
\label{h74}
\end{equation}
where 
\begin{equation}
c_0=\lim_{l\rightarrow\infty} \frac{|a_l|}{2|a_{l+1}|} .
\label{h75}
\end{equation}

The convergence of the series (\ref{h71}) can be also studied with the aid of
the formula (\ref{h48}) which, for integer values of $l$, is valid for any $x,
y, z, t$. If $a_l=a_{lx}+\alpha a_{ly}+\beta a_{lz}+\gamma a_{lt}$, and
\begin{eqnarray}
\label{h76a}
A_l=a_{lx}+a_{ly}+a_{lz}+a_{lt}, \\
A_l^\prime= a_{lx}-a_{ly}+a_{lz}-a_{lt} , \\
A_l^{\prime\prime}=a_{lx}+a_{ly}-a_{lz}-a_{lt},\\
A_l^{\prime\prime\prime}=a_{lx}-a_{ly}-a_{lz}+a_{lt} ,
\label{h76}
\end{eqnarray}
it can be shown with the aid of relations (\ref{h22}) and (\ref{h23}) that
\begin{equation}
a_l e=A_l e, \: a_l e^\prime=A_l^\prime e^\prime, \: 
a_l e^{\prime\prime}=A_l^{\prime\prime}e^{\prime\prime}, \:
a_l e^{\prime\prime\prime}=A_l^{\prime\prime\prime}e^{\prime\prime\prime} ,
\label{h77}
\end{equation}
so that the expression of the series (\ref{h71}) becomes
\begin{equation}
\sum_{l=0}^\infty \left(A_l  s^l e+
A_l^\prime  s^{\prime l}e^\prime+
A_l^{\prime\prime}s^{\prime\prime l}e^{\prime\prime}+
A_l^{\prime\prime\prime}s^{\prime\prime\prime l}e^{\prime\prime\prime}\right) ,
\label{h78}
\end{equation}
where the quantities $s, s^\prime, s^{\prime\prime}, s^{\prime\prime\prime}$
have been defined in Eq. (\ref{h8}).
The sufficient conditions for the absolute convergence of the series 
in Eq. (\ref{h78}) are that
\begin{equation}
\lim_{l\rightarrow \infty}\frac{|A_{l+1}||s|}{|A_l|}<1,
\lim_{l\rightarrow \infty}\frac{|A_{l+1}^\prime||s^\prime|}{|A_l^\prime|}<1,
\lim_{l\rightarrow
\infty}\frac{|A_{l+1}^{\prime\prime}||s^{\prime\prime}|}
{|A_l^{\prime\prime}|}<1,
\lim_{l\rightarrow \infty}\frac{|A_{l+1}^{\prime\prime\prime}|
|s^{\prime\prime\prime}|}{|A_l^{\prime\prime\prime}|}<1,
\label{h79}
\end{equation}
Thus the series in Eq. (\ref{h78}) is absolutely convergent for 
\begin{equation}
|x+y+z+t|<c,\:
|x-y+z-t|<c^\prime,\:
|x+y-z-t|<c^{\prime\prime},\:
|x-y-z+t|<c^{\prime\prime\prime},
\label{h80}
\end{equation}
where 
\begin{equation}
c=\lim_{l\rightarrow\infty} \frac{|A_l|}{|A_{l+1}|} ,\:
c^\prime=\lim_{l\rightarrow\infty} \frac{|A_l^\prime|}{|A_{l+1}^\prime|} ,\:
c^{\prime\prime}=\lim_{l\rightarrow\infty} 
\frac{|A_l^{\prime\prime}|}{|A_{l+1}^{\prime\prime}|} ,\:
c^{\prime\prime\prime}=\lim_{l\rightarrow\infty} 
\frac{|A_l^{\prime\prime\prime}|}{|A_{l+1}^{\prime\prime\prime}|} .
\label{h81}
\end{equation}
The relations (\ref{h80}) show that the region of convergence of the series
(\ref{h78}) is a four-dimensional parallelepiped.
It can be shown that $c_0=(1/2){\rm min}(c,c^\prime,c^{\prime\prime},
c^{\prime\prime\prime})$, where
${\rm min}$ designates the smallest of the numbers $c, c^\prime,
c^{\prime\prime},c^{\prime\prime\prime}$.
Using Eq. (\ref{h25c}), it can be seen that the circular region of
convergence defined in Eqs. (\ref{h74}), (\ref{h75})
is included in the parallelogram defined in Eqs. (\ref{h80}) and (\ref{h81}).

\subsection{Analytic functions of hyperbolic fourcomplex variables}

The fourcomplex function $f(u)$ of the fourcomplex variable $u$ has
been expressed in Eq. (\ref{g16}) in terms of 
the real functions $P(x,y,z,t),Q(x,y,z,t),R(x,y,z,t), S(x,y,z,t)$ of real
variables $x,y,z,t$. The
relations between the partial derivatives of the functions $P, Q, R, S$ are
obtained by setting succesively in   
Eq. (\ref{g17}) $\Delta x\rightarrow 0, \Delta y=\Delta z=\Delta t=0$;
then $ \Delta y\rightarrow 0, \Delta x=\Delta z=\Delta t=0;$  
then $  \Delta z\rightarrow 0,\Delta x=\Delta y=\Delta t=0$; and finally
$ \Delta t\rightarrow 0,\Delta x=\Delta y=\Delta z=0 $. 
The relations are 
\begin{equation}
\frac{\partial P}{\partial x} = \frac{\partial Q}{\partial y} =
\frac{\partial R}{\partial z} = \frac{\partial S}{\partial t},
\label{h89}
\end{equation}
\begin{equation}
\frac{\partial Q}{\partial x} = \frac{\partial P}{\partial y} =
\frac{\partial S}{\partial z} = \frac{\partial R}{\partial t},
\label{h90}
\end{equation}
\begin{equation}
\frac{\partial R}{\partial x} = \frac{\partial S}{\partial y} =
\frac{\partial P}{\partial z} = \frac{\partial Q}{\partial t},
\label{h91}
\end{equation}
\begin{equation}
\frac{\partial S}{\partial x} = \frac{\partial R}{\partial y} =
\frac{\partial Q}{\partial z} = \frac{\partial P}{\partial t}.
\label{h92}
\end{equation}

The relations (\ref{h89})-(\ref{h92}) are analogous to the Riemann relations
for the real and imaginary components of a complex function. It can be shown
from Eqs. (\ref{h89})-(\ref{h92}) that the component $P$ is a solution
of the equations 
\begin{equation}
\frac{\partial^2 P}{\partial x^2}-\frac{\partial^2 P}{\partial y^2}=0,
\:\: 
\frac{\partial^2 P}{\partial x^2}-\frac{\partial^2 P}{\partial z^2}=0,
\:\:
\frac{\partial^2 P}{\partial y^2}-\frac{\partial^2 P}{\partial t^2}=0,
\:\:
\frac{\partial^2 P}{\partial z^2}-\frac{\partial^2 P}{\partial t^2}=0,
\:\:
\label{h93}
\end{equation}
\begin{equation}
\frac{\partial^2 P}{\partial x^2}-\frac{\partial^2 P}{\partial t^2}=0,
\:\:
\frac{\partial^2 P}{\partial y^2}-\frac{\partial^2 P}{\partial z^2}=0,
\label{h94}
\end{equation}
and the components $Q, R, S$ are solutions of similar equations.
As can be seen from Eqs. (\ref{h93})-(\ref{h94}), the components $P, Q, R, S$
of 
an analytic function of hyperbolic fourcomplex variable are solutions of the
wave  
equation with respect to pairs of the variables $x,y,z,t$.
The component $P$ is also a solution of the mixed-derivative
equations 
\begin{equation}
\frac{\partial^2 P}{\partial x\partial y}=\frac{\partial^2 P}{\partial
z\partial t} ,
\:\: 
\frac{\partial^2 P}{\partial x\partial z}=\frac{\partial^2 P}{\partial
y\partial t} ,
\:\: 
\frac{\partial^2 P}{\partial x\partial t}=\frac{\partial^2 P}{\partial
y\partial z} ,
\end{equation}
and the components $Q, R, S$ are solutions of similar equations.

\subsection{Integrals of functions of hyperbolic fourcomplex variables}

The singularities of hyperbolic fourcomplex functions arise from terms of the
form 
$1/(u-u_0)^m$, with $m>0$. Functions containing such terms are singular not
only at $u=u_0$, but also at all points of the two-dimensional hyperplanes
passing through $u_0$ and which are parallel to the nodal hyperplanes. 

The integral of a hyperbolic fourcomplex function between two points $A, B$
along a path 
situated in a region free of singularities is independent of path, which means
that the integral of an analytic function along a loop situated in a region
free from singularities is zero,
\begin{equation}
\oint_\Gamma f(u) du = 0,
\label{h105}
\end{equation}
where it is supposed that a surface $\Sigma$ spanning 
the closed loop $\Gamma$ is not intersected by any of
the two-dimensional hyperplanes associated with the
singularities of the function $f(u)$. Using the expression, Eq. (\ref{g16}),
for $f(u)$ and the fact that $du=dx+\alpha  dy+\beta dz+\gamma dt$, the
explicit form of the integral in Eq. (\ref{h105}) is
\begin{eqnarray}
\lefteqn{\oint _\Gamma f(u) du = \oint_\Gamma
[(Pdx+Qdy+Rdz+Sdt)+\alpha(Qdx+Pdy+Sdz+Rdt)\nonumber}\\
&&+\beta(Rdx+Sdy+Pdz+Qdt)+\gamma(Sdx+Rdy+Qdz+Pdt)] .
\label{h106}
\end{eqnarray}
If the functions $P, Q, R, S$ are regular on a surface $\Sigma$
spanning the loop $\Gamma$,
the integral along the loop $\Gamma$ can be transformed with the aid of the
theorem of Stokes in an integral over the surface $\Sigma$ of terms of the form
$\partial P/\partial y -  \partial Q/\partial x, \:\:
\partial P/\partial z - \partial R/\partial x,\:\:
\partial P/\partial t - \partial S/\partial x, \:\:
\partial Q/\partial z -  \partial R/\partial y, \:\:
\partial Q/\partial t - \partial S/\partial y,\:\:
\partial R/\partial t - \partial S/\partial z$ and of similar terms arising
from the $\alpha, \beta$ and $\gamma$ components, 
which are equal to zero by Eqs. (\ref{h89})-(\ref{h92}), and this proves Eq.
(\ref{h105}). 

The exponential form of the hyperbolic fourcomplex numbers, Eq. (\ref{h40x}),
contains no 
cyclic variable, and therefore the concept of residue is not applicable to the
hyperbolic fourcomplex numbers defined in Eqs. (\ref{h1}).

\subsection{Factorization of hyperbolic fourcomplex polynomials}

A polynomial of degree $m$ of the hyperbolic fourcomplex variable 
$u=x+\alpha y+\beta z+\gamma t$ has the form
\begin{equation}
P_m(u)=u^m+a_1 u^{m-1}+\cdots+a_{m-1} u +a_m ,
\label{h106b}
\end{equation}
where the constants are in general hyperbolic fourcomplex numbers.
If $a_m=a_{mx}+\alpha a_{my}+\beta a_{mz}+\gamma a_{mt}$, and with the
notations of Eqs. (\ref{h8}) and (\ref{h76a})-(\ref{h76}) applied for $l=0, 1,
\cdots, m$ , the 
polynomial $P_m(u)$ can be written as 
\begin{eqnarray}
\lefteqn{P_m= \left[s^m 
+A_1 s^{m-1}+\cdots+A_{m-1} s+ A_m \right] e
+\left[s^{\prime m} 
+A_1^\prime s^{\prime m-1} +\cdots+A_{m-1}^\prime s^\prime+ A_m^\prime
\right]e^\prime\nonumber}\\ 
 & &
+\left[s^{\prime\prime m} 
+A_1^{\prime\prime} s^{\prime\prime m-1} +\cdots+A_{m-1}^{\prime\prime}
s^{\prime\prime}+ A_m^{\prime\prime} \right]e^{\prime\prime} 
+\left[s^{\prime\prime\prime m} 
+A_1^{\prime\prime\prime} s^{\prime\prime\prime m-1}
+\cdots+A_{m-1}^{\prime\prime\prime} s^{\prime\prime\prime}+
A_m^{\prime\prime\prime} \right]e^{\prime\prime\prime}.\nonumber \\ 
&&
\label{h107}
\end{eqnarray}
Each of the polynomials of degree $m$ with real coefficients in Eq.
(\ref{h107}) 
can be written as a product
of linear or quadratic factors with real coefficients, or as a product of
linear factors which, if imaginary, appear always in complex conjugate pairs.
Using the latter form for the simplicity of notations, the polynomial $P_m$
can be written as
\begin{equation}
P_m=\prod_{l=1}^m (s-s_l)e
+\prod_{l=1}^m (s^\prime-s_l^\prime)e^\prime
+\prod_{l=1}^m (s^{\prime\prime}-s_l^{\prime\prime})e^{\prime\prime}
+\prod_{l=1}^m
(s^{\prime\prime\prime}-s_l^{\prime\prime\prime})e^{\prime\prime\prime} ,
\label{h108}
\end{equation}
where the quantities $s_l$ appear always in complex conjugate pairs, and the
same is true for the quantities $s_l^\prime$, for the quantitites $
s_l^{\prime\prime}$, and for the quantities $s_l^{\prime\prime\prime}$. 
Due to the properties in Eqs. (\ref{h22}) and (\ref{h23}),
the polynomial $P_m(u)$ can be written as a product of factors of
the form  
\begin{equation}
P_m(u)=\prod_{l=1}^m \left[(s-s_l)e
+(s^\prime-s_l^\prime)e^\prime
+(s^{\prime\prime}-s_l^{\prime\prime})e^{\prime\prime}
+(s^{\prime\prime\prime}-s_l^{\prime\prime\prime})e^{\prime\prime\prime}
\right].
\label{h109}
\end{equation}
These relations can be written with the aid of Eq. (\ref{h24}) as
\begin{eqnarray}
P_m(u)=\prod_{p=1}^m (u-u_p) ,
\label{h128c}
\end{eqnarray}
where
\begin{eqnarray}
u_p=s_p e+s_p^{\prime}e^{\prime}
+s_p^{\prime\prime}e^{\prime\prime}
+s_p^{\prime\prime\prime}e^{\prime\prime\prime}.
\label{h128d}
\end{eqnarray}
The roots $s_p, s_p^{\prime}, s_p^{\prime\prime}, s_p^{\prime\prime\prime}$
of the corresponding polynomials in Eq. (\ref{h108}) may be ordered
arbitrarily. 
This means that Eq. (\ref{h128d}) gives sets of $m$ roots
$u_1,...,u_m$ of the polynomial $P_m(u)$, 
corresponding to the various ways in which the roots
$s_p, s_p^{\prime}, s_p^{\prime\prime}, s_p^{\prime\prime\prime}$
are ordered according to $p$ in each
group. Thus, while the hypercomplex components in Eq. (\ref{h107}) taken
separately have unique factorizations, the polynomial $P_m(u)$ can be written
in many different ways as a product of linear factors. 

If $P(u)=u^2-1$, 
the factorization in Eq. (\ref{h128c}) is $u^2-1=(u-u_1)(u-u_2)$, where 
$u_1=\pm  e\pm  e^\prime\pm  e^{\prime\prime}\pm  e^{\prime\prime\prime}$, 
$u_2=-u_1$, 
so that there are 8 distinct factorizations of $u^2-1$,
\begin{eqnarray}
\begin{array}{l}
u^2-1=\left(u-e-e^{\prime}-e^{\prime\prime}-e^{\prime\prime\prime}\right)
\left(u+e+e^{\prime}+e^{\prime\prime}+e^{\prime\prime\prime}\right),\\
u^2-1=\left(u-e-e^{\prime}-e^{\prime\prime}+e^{\prime\prime\prime}\right)
\left(u+e+e^{\prime}+e^{\prime\prime}-e^{\prime\prime\prime}\right),\\
u^2-1=\left(u-e-e^{\prime}+e^{\prime\prime}-e^{\prime\prime\prime}\right)
\left(u+e+e^{\prime}-e^{\prime\prime}+e^{\prime\prime\prime}\right),\\
u^2-1=\left(u-e+e^{\prime}-e^{\prime\prime}-e^{\prime\prime\prime}\right)
\left(u+e-e^{\prime}+e^{\prime\prime}+e^{\prime\prime\prime}\right),\\
u^2-1=\left(u-e-e^{\prime}+e^{\prime\prime}+e^{\prime\prime\prime}\right)
\left(u+e+e^{\prime}-e^{\prime\prime}-e^{\prime\prime\prime}\right),\\
u^2-1=\left(u-e+e^{\prime}-e^{\prime\prime}+e^{\prime\prime\prime}\right)
\left(u+e-e^{\prime}+e^{\prime\prime}-e^{\prime\prime\prime}\right),\\
u^2-1=\left(u-e+e^{\prime}+e^{\prime\prime}-e^{\prime\prime\prime}\right)
\left(u+e-e^{\prime}-e^{\prime\prime}+e^{\prime\prime\prime}\right),\\
u^2-1=\left(u-e+e^{\prime}+e^{\prime\prime}+e^{\prime\prime\prime}\right)
\left(u+e-e^{\prime}-e^{\prime\prime}-e^{\prime\prime\prime}\right).
\end{array}
\label{h129}
\end{eqnarray}
It can be checked that 
$\left\{\pm  e\pm  e^\prime\pm  e^{\prime\prime}
\pm  e^{\prime\prime\prime}\right\}^2= 
e+e^{\prime}+e^{\prime\prime}+e^{\prime\prime\prime}=1$.

\subsection{Representation of hyperbolic 
fourcomplex numbers by irreducible matrices}

If $T$ is the unitary matrix,
\begin{equation}
T =\left(
\begin{array}{cccc}
\frac{1}{{2}}&\frac{1}{{2}}  &\frac{1}{{2}}    &\frac{1}{{2}}    \\
\frac{1}{{2}}  &-\frac{1}{{2}}  &\frac{1}{{2}} &-\frac{1}{{2}}   \\
\frac{1}{{2}}&\frac{1}{{2}}  &-\frac{1}{{2}}    &-\frac{1}{{2}}    \\
\frac{1}{{2}}&-\frac{1}{{2}}  &-\frac{1}{{2}}    &\frac{1}{{2}}    \\
\end{array}
\right),
\label{h129x}
\end{equation}
it can be shown 
that the matrix $T U T^{-1}$ has the form 
\begin{equation}
T U T^{-1}=\left(
\begin{array}{cccc}
x+y+z+t &    0    &     0     &  0       \\
0       & x-y+z-t &     0     &  0       \\
0       &    0    &  x+y-z-t  &  0       \\
0       &    0    &     0     &  x-y-z+t \\
\end{array}
\right),
\label{h129y}
\end{equation}
where $U$ is the matrix in Eq. (\ref{h26}) used to represent the hyperbolic
fourcomplex 
number $u$. 
The relations between the variables $x+y+z+t, x-y+z-t, x+y-z-t, x-y-z+t$
for the multiplication
of hyperbolic fourcomplex numbers have been written in Eqs.
(\ref{h15})-(\ref{h18}). 
The matrix $T U T^{-1}$ provides an irreducible representation
\cite{4} of the hyperbolic fourcomplex number $u$ in terms of matrices with
real 
coefficients.

\section{Planar Complex Numbers in Four Dimensions}

\subsection{Operations with planar fourcomplex numbers}

A planar fourcomplex number is determined by its four components $(x,y,z,t)$.
The sum of the planar fourcomplex numbers $(x,y,z,t)$ and
$(x^\prime,y^\prime,z^\prime,t^\prime)$ is the planar fourcomplex
number $(x+x^\prime,y+y^\prime,z+z^\prime,t+t^\prime)$. 
The product of the planar fourcomplex numbers
$(x,y,z,t)$ and $(x^\prime,y^\prime,z^\prime,t^\prime)$ 
is defined in this work to be the planar fourcomplex
number
$(xx^\prime-yt^\prime-zz^\prime-ty^\prime,
xy^\prime+yx^\prime-zt^\prime-tz^\prime,
xz^\prime+yy^\prime+zx^\prime-tt^\prime,
xt^\prime+yz^\prime+zy^\prime+tx^\prime)$.

Planar fourcomplex numbers and their operations can be represented by  writing
the 
planar fourcomplex number $(x,y,z,t)$ as  
$u=x+\alpha y+\beta z+\gamma t$, where $\alpha, \beta$ and $\gamma$ 
are bases for which the multiplication rules are 
\begin{equation}
\alpha^2=\beta, \:\beta^2=-1, \:\gamma^2=-\beta,
\alpha\beta=\beta\alpha=\gamma,\: 
\alpha\gamma=\gamma\alpha=-1, \:\beta\gamma=\gamma\beta=-\alpha .
\label{c2-1}
\end{equation}
Two planar fourcomplex numbers $u=x+\alpha y+\beta z+\gamma t, 
u^\prime=x^\prime+\alpha y^\prime+\beta z^\prime+\gamma t^\prime$ are equal, 
$u=u^\prime$, if and only if $x=x^\prime, y=y^\prime,
z=z^\prime, t=t^\prime$. 
If 
$u=x+\alpha y+\beta z+\gamma t, 
u^\prime=x^\prime+\alpha y^\prime+\beta z^\prime+\gamma t^\prime$
are planar fourcomplex numbers, 
the sum $u+u^\prime$ and the 
product $uu^\prime$ defined above can be obtained by applying the usual
algebraic rules to the sum 
$(x+\alpha y+\beta z+\gamma t)+ 
(x^\prime+\alpha y^\prime+\beta z^\prime+\gamma t^\prime)$
and to the product 
$(x+\alpha y+\beta z+\gamma t)
(x^\prime+\alpha y^\prime+\beta z^\prime+\gamma t^\prime)$,
and grouping of the resulting terms,
\begin{equation}
u+u^\prime=x+x^\prime+\alpha(y+y^\prime)+\beta(z+z^\prime)+\gamma(t+t^\prime),
\label{c2-1a}
\end{equation}
\begin{eqnarray}
\lefteqn{uu^\prime=
xx^\prime-yt^\prime-zz^\prime-ty^\prime+
\alpha(xy^\prime+yx^\prime-zt^\prime-tz^\prime)+
\beta(xz^\prime+yy^\prime+zx^\prime-tt^\prime)\nonumber}\\
&&+\gamma(xt^\prime+yz^\prime+zy^\prime+tx^\prime).
\label{c2-1b}
\end{eqnarray}

If $u,u^\prime,u^{\prime\prime}$ are planar fourcomplex numbers, the
multiplication  
is associative
\begin{equation}
(uu^\prime)u^{\prime\prime}=u(u^\prime u^{\prime\prime})
\label{c2-2}
\end{equation}
and commutative
\begin{equation}
u u^\prime=u^\prime u ,
\label{c2-3}
\end{equation}
as can be checked through direct calculation.
The planar fourcomplex zero is $0+\alpha\cdot 0+\beta\cdot 0+\gamma\cdot 0,$ 
denoted simply 0, 
and the planar fourcomplex unity is $1+\alpha\cdot 0+\beta\cdot 0+\gamma\cdot
0,$  
denoted simply 1.

The inverse of the planar fourcomplex number 
$u=x+\alpha y+\beta z+\gamma t$ is a planar fourcomplex number
$u^\prime=x^\prime+\alpha y^\prime+\beta z^\prime+\gamma t^\prime$
having the property that
\begin{equation}
uu^\prime=1 .
\label{c2-4}
\end{equation}
Written on components, the condition, Eq. (\ref{c2-4}), is
\begin{equation}
\begin{array}{c}
xx^\prime-ty^\prime-zz^\prime-yt^\prime=1,\\
yx^\prime+xy^\prime-tz^\prime-zt^\prime=0,\\
zx^\prime+yy^\prime+xz^\prime-tt^\prime=0,\\
tx^\prime+zy^\prime+yz^\prime+xt^\prime=0.
\end{array}
\label{c2-5}
\end{equation}
The system (\ref{c2-5}) has the solution
\begin{equation}
x^\prime=\frac{x(x^2+z^2)-z(y^2-t^2)+2xyt}
{\rho^4} ,
\label{c2-6a}
\end{equation}
\begin{equation}
y^\prime=-\frac{y(x^2-z^2)+t(y^2+t^2)+2xzt}
{\rho^4} ,
\label{c2-6c}
\end{equation}
\begin{equation}
z^\prime=
\frac{-z(x^2+z^2)+x(y^2-t^2)+2zyt}
{\rho^4} ,
\label{c2-6b}
\end{equation}
\begin{equation}
t^\prime=-\frac{t(x^2-z^2)+y(y^2+t^2)-2xyz}
{\rho^4} ,
\label{c2-6d}
\end{equation}
provided that $\rho\not=0, $ where
\begin{equation}
\rho^4=x^4+z^4+y^4+t^4+2x^2z^2+2y^2t^2+4x^2yt-4xy^2z+4xzt^2-4yz^2t .
\label{c2-6e}
\end{equation}
The quantity $\rho$ will be called amplitude of the planar fourcomplex number
$x+\alpha y+\beta z +\gamma t$.
Since
\begin{equation}
\rho^4=\rho_+^2\rho_-^2 ,
\label{c2-7a}
\end{equation}
where
\begin{equation}
\rho_+^2=\left(x+\frac{y-t}{\sqrt{2}}\right)^2+\left(z+\frac{y+t}
{\sqrt{2}}\right)^2,
\: \rho_-^2= \left(x-\frac{y-t}{\sqrt{2}}\right)^2+\left(z-\frac{y+t}
{\sqrt{2}}\right)^2, 
\label{c2-7b}
\end{equation}
a planar fourcomplex number $u=x+\alpha y+\beta z+\gamma t$ has an inverse,
unless 
\begin{equation}
x+\frac{y-t}{\sqrt{2}}=0,\: z+\frac{y+t}{\sqrt{2}}=0 ,  
\label{c2-8}
\end{equation}
or
\begin{equation}
x-\frac{y-t}{\sqrt{2}}=0,\: z-\frac{y+t}{\sqrt{2}}=0 .  
\label{c2-9}
\end{equation}

Because of conditions (\ref{c2-8})-(\ref{c2-9}) these 2-dimensional
hypersurfaces 
will be called nodal hyperplanes. 
It can be shown that 
if $uu^\prime=0$ then either $u=0$, or $u^\prime=0$, 
or one of the planar fourcomplex numbers is of the form $x+\alpha
(x+z)/\sqrt{2} +\beta z-\gamma (x-z)/\sqrt{2}$ 
and the other of the form $x^\prime-\alpha (x^\prime+z^\prime)/\sqrt{2}
+\beta z^\prime+\gamma (x^\prime-z^\prime)/\sqrt{2}$.

\subsection{Geometric representation of planar fourcomplex numbers}

The planar fourcomplex number $x+\alpha y+\beta z+\gamma t$ can be represented
by  
the point $A$ of coordinates $(x,y,z,t)$. 
If $O$ is the origin of the four-dimensional space $x,y,z,t,$ the distance 
from $A$ to the origin $O$ can be taken as
\begin{equation}
d^2=x^2+y^2+z^2+t^2 .
\label{c2-10}
\end{equation}
The distance $d$ will be called modulus of the planar fourcomplex number
$x+\alpha 
y+\beta z +\gamma t$, $d=|u|$. 
The orientation in the four-dimensional space of the line $OA$ can be specified
with the aid of three angles $\phi, \chi, \psi$
defined with respect to the rotated system of axes
\begin{equation}
\xi=\frac{x}{\sqrt{2}}+\frac{y-t}{2}, \: 
\tau=\frac{x}{\sqrt{2}}-\frac{y-t}{2}, \:
\upsilon=\frac{z}{\sqrt{2}}+\frac{y+t}{2}, \: 
\zeta=-\frac{z}{\sqrt{2}}+\frac{y+t}{2} .
\label{c2-11}
\end{equation}
The variables $\xi, \upsilon, \tau, \zeta$ will be called canonical 
planar fourcomplex variables.
The use of the rotated axes $\xi, \upsilon, \tau, \zeta$ 
for the definition of the angles $\phi, \chi, \psi$ 
is convenient for the expression of the planar fourcomplex numbers
in exponential and trigonometric forms, as it will be discussed further.
The angle $\phi$ is the angle between the projection of $A$ in the plane
$\xi,\upsilon$ and the $O\xi$ axis, $0\leq\phi<2\pi$,  
$\chi$ is the angle between the projection of $A$ in the plane $\tau,\zeta$ and
the  $O\tau$ axis, $0\leq\chi<2\pi$,
and $\psi$ is the angle between the line $OA$ and the plane $\tau O \zeta
$, $0\leq \psi\leq\pi/2$, 
as shown in Fig. 1.
The definition of the variables in this section is different from the
definition used for the circular fourcomplex numbers, because the definition of
the rotated axes in Eq. (\ref{c2-11}) is different from the definition of the
rotated 
circular axes, Eq. (\ref{11}). 
The angles $\phi$ and $\chi$ will be called azimuthal
angles, the angle $\psi$ will be called planar angle .
The fact that $0\leq \psi\leq\pi/2$ means that $\psi$ has
the same sign on both faces of the two-dimensional hyperplane $\upsilon O
\zeta$. The components of the point $A$
in terms of the distance $d$ and the angles $\phi, \chi, \psi$ are thus
\begin{equation}
\frac{x}{\sqrt{2}}+\frac{y-t}{2}=d\cos\phi \sin\psi , 
\label{c2-12a}
\end{equation}
\begin{equation}
\frac{x}{\sqrt{2}}-\frac{y-t}{2}=d\cos\chi \cos\psi , 
\label{c2-12b}
\end{equation}
\begin{equation}
\frac{z}{\sqrt{2}}+\frac{y+t}{2}=d\sin\phi \sin\psi , 
\label{c2-12c}
\end{equation}
\begin{equation}
-\frac{z}{\sqrt{2}}+\frac{y+t}{2}=d\sin\chi \cos\psi .
\label{c2-12d}
\end{equation}
It can be checked that $\rho_+=\sqrt{2}d\sin\psi, \rho_-=\sqrt{2}d\cos\psi$.
The coordinates $x,y,z,t$ in terms of the variables $d, \phi, \chi,
\psi$ are
\begin{equation}
x=\frac{d}{\sqrt{2}}(\cos\phi\sin\psi+\cos\chi\cos\psi),
\label{c2-12e}
\end{equation}
\begin{equation}
y=\frac{d}{\sqrt{2}}[\sin(\phi+\pi/4)\sin\psi+\sin(\chi-\pi/4)\cos\psi],
\label{c2-12g}
\end{equation}
\begin{equation}
z=\frac{d}{\sqrt{2}}(\sin\phi\sin\psi-\sin\chi\cos\psi),
\label{c2-12f}
\end{equation}
\begin{equation}
t=\frac{d}{\sqrt{2}}[-\cos(\phi+\pi/4)\sin\psi+\cos(\chi-\pi/4)\cos\psi].
\label{c2-12h}
\end{equation}
The angles $\phi, \chi, \psi$ can be expressed in terms of the coordinates
$x,y,z,t$ as 
\begin{equation}
\sin\phi = \frac{z+(y+t)/\sqrt{2}}{\rho_+} ,\: 
\cos\phi = \frac{x+(y-t)/\sqrt{2}}{\rho_+} ,
\label{c2-13a}
\end{equation}
\begin{equation}
\sin\chi = \frac{-z+(y+t)/\sqrt{2}}{\rho_-} ,\: 
\cos\chi = \frac{x-(y-t)/\sqrt{2}}{\rho_-} ,
\label{c2-13b}
\end{equation}
\begin{equation}
\tan\psi=\rho_+/\rho_- .
\label{c2-13c}
\end{equation}
The nodal hyperplanes are $\xi O\upsilon$, for which $\tau=0, \zeta=0$, and
$\tau O\zeta$, for which $\xi=0, \upsilon=0$.
For points in 
the nodal hyperplane $\xi O\upsilon$ the planar angle is $\psi=\pi/2$, 
for points in the nodal hyperplane $\tau O\zeta$ the planar angle is $\psi=0$.

It can be shown that if $u_1=x_1+\alpha y_1+\beta z_1+\gamma t_1, 
u_2=x_2+\alpha y_2+\beta z_2+\gamma t_2$ are planar fourcomplex
numbers of amplitudes and angles $\rho_1, \phi_1, \chi_1, \psi_1$ and
respectively $\rho_2, \phi_2, \chi_2, \psi_2$, then the amplitude $\rho$ and
the angles $\phi, \chi, \psi$ of the product planar fourcomplex number $u_1u_2$
are 
\begin{equation}
\rho=\rho_1\rho_2, 
\label{c2-14a}
\end{equation}
\begin{equation}
 \phi=\phi_1+\phi_2, \: \chi=\chi_1+\chi_2, \: \tan\psi=\tan\psi_1\tan\psi_2 . 
\label{c2-14b}
\end{equation}
The relations (\ref{c2-14a})-(\ref{c2-14b}) are consequences of the definitions
(\ref{c2-6e})-(\ref{c2-7b}), (\ref{c2-13a})-(\ref{c2-13c}) and of the
identities 
\begin{eqnarray}
\lefteqn{\left[(x_1x_2-z_1z_2-y_1t_2-t_1y_2)
+\frac{(x_1y_2+y_1x_2-z_1t_2-t_1z_2)-(x_1t_2+t_1x_2+z_1y_2+y_1z_2)}{\sqrt{2}}
\right]^2
\nonumber}\\
&&+\left[(x_1z_2+z_1x_2+y_1y_2-t_1t_2)
+\frac{(x_1y_2+y_1x_2-z_1t_2-t_1z_2)+(x_1t_2+t_1x_2+z_1y_2+y_1z_2)}{\sqrt{2}}
\right]^2
\nonumber\\
&&=\left[\left(x_1+\frac{y_1-t_1}{\sqrt{2}}\right)^2
+\left(z_1+\frac{y_1+t_1}{\sqrt{2}}\right)^2\right]
\left[\left(x_2+\frac{y_2-t_2}{\sqrt{2}}\right)^2
+\left(z_2+\frac{y_2+t_2}{\sqrt{2}}\right)^2\right],
\label{c2-15}
\end{eqnarray}
\begin{eqnarray}
\lefteqn{\left[(x_1x_2-z_1z_2-y_1t_2-t_1y_2)
-\frac{(x_1y_2+y_1x_2-z_1t_2-t_1z_2)-(x_1t_2+t_1x_2+z_1y_2+y_1z_2)}{\sqrt{2}}
\right]^2
\nonumber}\\
&&+\left[(x_1z_2+z_1x_2+y_1y_2-t_1t_2)
-\frac{(x_1y_2+y_1x_2-z_1t_2-t_1z_2)+(x_1t_2+t_1x_2+z_1y_2+y_1z_2)}{\sqrt{2}}
\right]^2
\nonumber\\
&&=\left[\left(x_1-\frac{y_1-t_1}{\sqrt{2}}\right)^2
+\left(z_1-\frac{y_1+t_1}{\sqrt{2}}\right)^2\right]
\left[\left(x_2-\frac{y_2-t_2}{\sqrt{2}}\right)^2
+\left(z_2-\frac{y_2+t_2}{\sqrt{2}}\right)^2\right],
\label{c2-16}
\end{eqnarray}
\begin{eqnarray}
\lefteqn{(x_1x_2-z_1z_2-y_1t_2-t_1y_2)
+\frac{(x_1y_2+y_1x_2-z_1t_2-t_1z_2)-(x_1t_2+t_1x_2+z_1y_2+y_1z_2)}{\sqrt{2}}
\nonumber}\\
&&=\left(x_1+\frac{y_1-t_1}{\sqrt{2}}\right)
\left(x_2+\frac{y_2-t_2}{\sqrt{2}}\right)
-\left(z_1+\frac{y_1+t_1}{\sqrt{2}}\right)
\left(z_2+\frac{y_2+t_2}{\sqrt{2}}\right) ,
\label{c2-17}
\end{eqnarray}
\begin{eqnarray}
\lefteqn{(x_1z_2+z_1x_2+y_1y_2-t_1t_2)
+\frac{(x_1y_2+y_1x_2-z_1t_2-t_1z_2)+(x_1t_2+t_1x_2+z_1y_2+y_1z_2)}{\sqrt{2}}
\nonumber}\\
&&=\left(z_1+\frac{y_1+t_1}{\sqrt{2}}\right)
\left(x_2+\frac{y_2-t_2}{\sqrt{2}}\right)
+\left(x_1+\frac{y_1-t_1}{\sqrt{2}}\right)
\left(z_2+\frac{y_2+t_2}{\sqrt{2}}\right) ,
\label{c2-18}
\end{eqnarray}
\begin{eqnarray}
\lefteqn{(x_1x_2-z_1z_2-y_1t_2-t_1y_2)
-\frac{(x_1y_2+y_1x_2-z_1t_2-t_1z_2)-(x_1t_2+t_1x_2+z_1y_2+y_1z_2)}{\sqrt{2}}
\nonumber}\\
&&=\left(x_1-\frac{y_1-t_1}{\sqrt{2}}\right)
\left(x_2-\frac{y_2-t_2}{\sqrt{2}}\right)
-\left(-z_1+\frac{y_1+t_1}{\sqrt{2}}\right)
\left(-z_2+\frac{y_2+t_2}{\sqrt{2}}\right) ,
\label{c2-19}
\end{eqnarray}
\begin{eqnarray}
\lefteqn{-(x_1z_2+z_1x_2+y_1y_2-t_1t_2)
+\frac{(x_1y_2+y_1x_2-z_1t_2-t_1z_2)+(x_1t_2+t_1x_2+z_1y_2+y_1z_2)}{\sqrt{2}}
\nonumber}\\
&&=\left(-z_1+\frac{y_1+t_1}{\sqrt{2}}\right)
\left(x_2-\frac{y_2-t_2}{\sqrt{2}}\right)
+\left(x_1-\frac{y_1-t_1}{\sqrt{2}}\right)
\left(-z_2+\frac{y_2+t_2}{\sqrt{2}}\right) .
\label{c2-20}
\end{eqnarray}
The identities (\ref{c2-15}) and (\ref{c2-16}) can also be written as
\begin{equation}
\rho_+^2=\rho_{1+}\rho_{2+} ,
\label{c2-21a}
\end{equation}
\begin{equation}
\rho_-^2=\rho_{1-}\rho_{2-} ,
\label{c2-21b}
\end{equation}
where
\begin{equation}
\rho_{j+}^2=\left(x_j+\frac{y_j-t_j}{\sqrt{2}}\right)^2
+\left(z_j+\frac{y_j+t_j}{\sqrt{2}}\right)^2,
\: \rho_{j-}^2=\left(x_j-\frac{y_j-t_j}{\sqrt{2}}\right)^2
+\left(z-\frac{y_j+t_j}{\sqrt{2}}\right)^2 , 
\label{c2-22}
\end{equation}
for $j=1,2$.

The fact that the amplitude of the product is equal to the product of the 
amplitudes, as written in Eq. (\ref{c2-14a}), can 
be demonstrated also by using a representation of the multiplication of the 
planar fourcomplex numbers by matrices, in which the planar fourcomplex number
$u=x+\alpha 
y+\beta z+\gamma t$ is represented by the matrix
\begin{equation}
A=\left(\begin{array}{cccc}
x &y &z &t\\
-t&x &y &z\\
-z&-t&x &y\\
-y&-z&-t&x 
\end{array}\right) .
\label{c2-23}
\end{equation}
The product $u=x+\alpha y+\beta z+\gamma t$ of the planar fourcomplex numbers
$u_1=x_1+\alpha y_1+\beta z_1+\gamma t_1, u_2=x_2+\alpha y_2+\beta z_2+\gamma
t_2$, can be represented by the matrix multiplication 
\begin{equation}
A=A_1A_2.
\label{c2-24}
\end{equation}
It can be checked that the determinant ${\rm det}(A)$ of the matrix $A$ is
\begin{equation}
{\rm det}A = \rho^4 .
\label{c2-25}
\end{equation}
The identity (\ref{c2-14a}) is then a consequence of the fact the determinant 
of the product of matrices is equal to the product of the determinants 
of the factor matrices. 

\subsection{The planar fourdimensional cosexponential functions}

The exponential function of a fourcomplex variable $u$ and the addition
theorem for the exponential function have been written in Eqs. (\ref{g1}) and
(\ref{g2}).
If $u=x+\alpha y+\beta z+\gamma t$, then $\exp u$ can be calculated as 
$\exp u=\exp x \cdot \exp (\alpha y) \cdot \exp (\beta z)\cdot \exp (\gamma
t)$. 
According to Eq. (\ref{c2-1}),  
\begin{eqnarray}
\begin{array}{l}
\alpha^{8m}=1, \alpha^{8m+1}=\alpha, \alpha^{8m+2}=\beta, \alpha^{8m+3}=\gamma,
\\ 
\alpha^{8m+4}=-1, \alpha^{8m+5}=-\alpha, \alpha^{8m+6}=-\beta,
\alpha^{8m+7}=-\gamma, \\ 
\beta^{4m}=1, \beta^{4m+1}=\beta, \beta^{4m+2}=-1, \beta^{4m+3}=-\beta, \\
\gamma^{8m}=1, \gamma^{8m+1}=\gamma, \gamma^{8m+2}=-\beta,
\gamma^{8m+3}=\alpha, \\ 
\gamma^{8m+4}=-1, \gamma^{8m+5}=-\gamma, \gamma^{8m+6}=\beta,
\gamma^{8m+7}=-\alpha, \\ 
\end{array}
\label{c2-28}
\end{eqnarray}
where $n$ is a natural number,
so that $\exp (\alpha y), \: \exp(\beta z)$ and $\exp(\gamma t)$ can be written
as 
\begin{equation}
\exp (\beta z) = \cos z +\beta \sin z ,
\label{c2-29}
\end{equation}
and
\begin{equation}
\exp (\alpha y) = f_{40}(y)+\alpha f_{41}(y)+\beta f_{42}(y) +\gamma f_{43}(y)
,   
\label{c2-30a}
\end{equation}
\begin{equation}
\exp (\gamma t) = f_{40}(t)+\gamma f_{41}(t)-\beta f_{42}(t) +\alpha f_{43}(t)
, 
\label{c2-30b}
\end{equation}
where the four-dimensional cosexponential functions $f_{40}, f_{41}, f_{42},
f_{43}$ are 
defined by the series 
\begin{equation}
f_{40}(x)=1-x^4/4!+x^8/8!-\cdots ,
\label{c2-30c}
\end{equation}
\begin{equation}
f_{41}(x)=x-x^5/5!+x^9/9!-\cdots,
\label{c2-30d}
\end{equation}
\begin{equation}
f_{42}(x)=x^2/2!-x^6/6!+x^{10}/10!-\cdots ,
\label{c2-30e}
\end{equation}
\begin{equation}
f_{43}(x)=x^3/3!-x^7/7!+x^{11}/11!-\cdots .
\label{c2-30f}
\end{equation}
The functions $f_{40}, f_{42}$ are even, the functions $f_{41}, f_{43}$ are
odd, 
\begin{equation}
f_{40}(-u)=f_{40}(u),\:f_{42}(-u)=f_{42}(u),\:f_{41}(-u)=-f_{41}(u),
\:f_{43}(-u)=-f_{43}(u).
\label{c2-30feo}
\end{equation}

Addition theorems for the four-dimensional cosexponential functions can be
obtained from the relation $\exp \alpha(x+y)=\exp\alpha x\cdot\exp\alpha y $,
by 
substituting the expression of the exponentials as given in Eq. (\ref{c2-30a}),
\begin{equation}
f_{40}(x+y)=f_{40}(x)f_{40}(y)-f_{41}(x)f_{43}(y)
-f_{42}(x)f_{42}(y)-f_{43}(x)f_{41}(y), 
\label{c2-30g}
\end{equation}
\begin{equation}
f_{41}(x+y)=f_{40}(x)f_{41}(y)+f_{41}(x)f_{40}(y)
-f_{42}(x)f_{43}(y)-f_{43}(x)f_{42}(y) ,
\label{c2-30h}
\end{equation}
\begin{equation}
f_{42}(x+y)=f_{40}(x)f_{42}(y)+f_{41}(x)f_{41}(y)
+f_{42}(x)f_{40}(y)-f_{43}(x)f_{43}(y) ,
\label{c2-30i}
\end{equation}
\begin{equation}
f_{43}(x+y)=f_{40}(x)f_{43}(y)+f_{41}(x)f_{42}(y)
+f_{42}(x)f_{41}(y)+f_{43}(x)f_{40}(y) .
\label{c2-30j}
\end{equation}
For $x=y$ the relations (\ref{c2-30g})-(\ref{c2-30j}) take the form
\begin{equation}
f_{40}(2x)=f_{40}^2(x)-f_{42}^2(x)-2f_{41}(x)f_{43}(x) ,
\label{c2-30gg}
\end{equation}
\begin{equation}
f_{41}(2x)=2f_{40}(x)f_{41}(x)-2f_{42}(x)f_{43}(x) ,
\label{c2-30hh}
\end{equation}
\begin{equation}
f_{42}(2x)=f_{41}^2(x)-f_{43}^2(x)+2f_{40}(x)f_{42}(x) ,
\label{c2-30ii}
\end{equation}
\begin{equation}
f_{43}(2x)=2f_{40}(x)f_{43}(x)+2f_{41}(x)f_{42}(x) .
\label{c2-30jj}
\end{equation}
For $x=-y$ the relations (\ref{c2-30g})-(\ref{c2-30j}) and (\ref{c2-30feo})
yield 
\begin{equation}
f_{40}^2(x)-f_{42}^2(x)+2f_{41}(x)f_{43}(x)=1 ,
\label{c2-30eog}
\end{equation}
\begin{equation}
f_{41}^2(x)-f_{43}^2(x)-2f_{40}(x)f_{42}(x)=0 .
\label{c2-30eoi}
\end{equation}
From Eqs. (\ref{c2-29})-(\ref{c2-30b}) it can be shown that, for $m$ integer,
\begin{equation}
(\cos z +\beta \sin z)^m=\cos mz +\beta \sin mz ,
\label{c2-2929}
\end{equation}
and
\begin{equation}
[f_{40}(y)+\alpha f_{41}(y)+\beta f_{42}(y) +\gamma f_{43}(y)]^m= 
f_{40}(my)+\alpha f_{41}(my)+\beta f_{42}(my) +\gamma f_{43}(my),  
\label{c2-30a30a}
\end{equation}
\begin{equation}
[f_{40}(t)+\gamma f_{41}(t)-\beta f_{42}(t) +\alpha f_{43}(t)]^m=
f_{40}(mt)+\gamma f_{41}(mt)-\beta f_{42}(mt) +\alpha f_{43}(mt) .
\label{c2-30b30b}
\end{equation}

Since
\begin{equation}
(\alpha-\gamma)^{2m}=2^m,\: (\alpha-\gamma)^{2m+1}=2^m(\alpha-\gamma) ,
\label{c2-30k}
\end{equation}
it can be shown from the definition of the exponential function, Eq.
(\ref{g1}) that
\begin{equation}
\exp(\alpha-\gamma)x=\cosh\sqrt{2}x+\frac{\alpha-\gamma}{\sqrt{2}}
\sinh\sqrt{2}x .
\label{c2-30l}
\end{equation}
Substituting in the relation $\exp(\alpha-\gamma)x=\exp\alpha x\exp(-\gamma x)$
the 
expression of the exponentials from Eqs. (\ref{c2-30a}), (\ref{c2-30b}) and
(\ref{c2-30l}) yields  
\begin{equation}
f_{40}^2+f_{41}^2+f_{42}^2+f_{43}^2=\cosh\sqrt{2}x ,
\label{c2-30m}
\end{equation}
\begin{equation}
f_{40}f_{41}-f_{40}f_{43}+f_{41}f_{42}
+f_{42}f_{43}=\frac{1}{\sqrt{2}}\sinh\sqrt{2}x, 
\label{c2-30n}
\end{equation}
where $f_{40}, f_{41}, f_{42}, f_{43}$ are functions of $x$.
From relations (\ref{c2-30m}) and (\ref{c2-30n}) it can be inferred that
\begin{equation}
\left(f_{40}+\frac{f_{41}-f_{43}}{\sqrt{2}}\right)^2+
\left(f_{42}+\frac{f_{41}+f_{43}}{\sqrt{2}}\right)^2=\exp\sqrt{2}x ,
\label{c2-30o}
\end{equation}
\begin{equation}
\left(f_{40}-\frac{f_{41}-f_{43}}{\sqrt{2}}\right)^2+
\left(f_{42}-\frac{f_{41}+f_{43}}{\sqrt{2}}\right)^2=\exp(-\sqrt{2}x) ,
\label{c2-30p}
\end{equation}
which means that
\begin{equation}
\left[\left(f_{40}+\frac{f_{41}-f_{43}}{\sqrt{2}}\right)^2+
\left(f_{42}+\frac{f_{41}+f_{43}}{\sqrt{2}}\right)^2\right]
\left[\left(f_{40}-\frac{f_{41}-f_{43}}{\sqrt{2}}\right)^2+
\left(f_{42}-\frac{f_{41}+f_{43}}{\sqrt{2}}\right)^2\right]=1 .
\label{c2-30q}
\end{equation}
An equivalent form of Eq. (\ref{c2-30q}) is
\begin{equation}
f_{40}^4+f_{41}^4+f_{42}^4+f_{43}^4+2(f_{40}^2f_{42}^2+f_{41}^2f_{43}^2)
+4(f_{40}^2f_{41}f_{43}+
f_{40}f_{42}f_{43}^2-f_{40}f_{41}^2f_{42}-f_{41}f_{42}^2f_{43})=1.
\label{c2-30r}
\end{equation}
The form of this relation is similar to the expression in Eq. (\ref{c2-6e}).
Similarly, since
\begin{equation}
(\alpha+\gamma)^{2m}=(-1)^m2^m,\:
(\alpha+\gamma)^{2m+1}=(-1)^m2^m(\alpha+\gamma) , 
\label{c2-30k30k}
\end{equation}
it can be shown from the definition of the exponential function, Eq.
(\ref{g1}) 
that
\begin{equation}
\exp(\alpha+\gamma)x=\cos\sqrt{2}x+\frac{\alpha+\gamma}{\sqrt{2}}
\sin\sqrt{2}x .
\label{c2-30l30l}
\end{equation}
Substituting in the relation $\exp(\alpha+\gamma)x=\exp\alpha x\exp\gamma x$
the 
expression of the exponentials from Eqs. (\ref{c2-30a}), (\ref{c2-30b}) and
(\ref{c2-30l30l}) yields  
\begin{equation}
f_{40}^2-f_{41}^2+f_{42}^2-f_{43}^2=\cos\sqrt{2}x ,
\label{c2-30m30m}
\end{equation}
\begin{equation}
f_{40}f_{41}+f_{40}f_{43}-f_{41}f_{42}+f_{42}f_{43}
=\frac{1}{\sqrt{2}}\sin\sqrt{2}x, 
\label{c2-30n30n}
\end{equation}
where $f_{40}, f_{41}, f_{42}, f_{43}$ are functions of $x$.

Expressions of the four-dimensional cosexponential functions
(\ref{c2-30c})-(\ref{c2-30f}) can be obtained using the fact that
$[(1+i)/\sqrt{2}]^4=-1$, so that 
\begin{equation}
f_{40}(x)=\frac{1}{2}\left(\cosh\frac{1+i}{\sqrt{2}}x+
\cos\frac{1+i}{\sqrt{2}}x\right),
\label{c2-30s}
\end{equation}
\begin{equation}
f_{41}(x)=\frac{1}{\sqrt{2}(1+i)}\left(\sinh\frac{1+i}{\sqrt{2}}x+
\sin\frac{1+i}{\sqrt{2}}x\right),
\label{c2-30t}
\end{equation}
\begin{equation}
f_{42}(x)=\frac{1}{2i}\left(\cosh\frac{1+i}{\sqrt{2}}x-
\cos\frac{1+i}{\sqrt{2}}x\right),
\label{c2-30u}
\end{equation}
\begin{equation}
f_{43}(x)=\frac{1}{\sqrt{2}(-1+i)}\left(\sinh\frac{1+i}{\sqrt{2}}x-
\sin\frac{1+i}{\sqrt{2}}x\right).
\label{c2-30w}
\end{equation}
Using the addition theorems for the functions in the right-hand sides of Eqs.
(\ref{c2-30s})-(\ref{c2-30w}), the expressions of the four-dimensional
cosexponential 
functions become
\begin{equation}
f_{40}(x)=\cos\frac{x}{\sqrt{2}}\cosh\frac{x}{\sqrt{2}} ,
\label{c2-30c30c}
\end{equation}
\begin{equation}
f_{41}(x)=\frac{1}{\sqrt{2}}
\left(\sin\frac{x}{\sqrt{2}}\cosh\frac{x}{\sqrt{2}}+
\sinh\frac{x}{\sqrt{2}}\cos\frac{x}{\sqrt{2}}\right),
\label{c2-30d30d}
\end{equation}
\begin{equation}
f_{42}(x)=\sin\frac{x}{\sqrt{2}}\sinh\frac{x}{\sqrt{2}}, 
\label{c2-30e30e}
\end{equation}
\begin{equation}
f_{43}(x)=\frac{1}{\sqrt{2}}
\left(\sin\frac{x}{\sqrt{2}}\cosh\frac{x}{\sqrt{2}}-
\sinh\frac{x}{\sqrt{2}}\cos\frac{x}{\sqrt{2}}\right) .
\label{c2-30f30f}
\end{equation}
It is remarkable that the series in Eqs. (\ref{c2-30c30c})-(\ref{c2-30f30f}),
in 
which the  
terms are either of the form $x^{4m}$, or $x^{4m+1}$, or $x^{4m+2}$,
or $x^{4m+3}$ can be expressed
in terms of elementary functions whose power series are not subject to such 
restrictions. 
The graphs of the four-dimensional cosexponential functions are shown in 
Fig. 3.

It can be checked that the cosexponential functions are solutions of the
fourth-order differential equation
\begin{equation}
\frac{{\rm d}^4\zeta}{{\rm d}u^4}=-\zeta ,
\label{c2-30x}
\end{equation}
whose solutions are of the form
$\zeta(u)=Af_{40}(u)+Bf_{41}(u)+Cf_{42}(u)+Df_{43}(u).$ 
It can also be checked that the derivatives of the cosexponential functions are
related by
\begin{equation}
\frac{df_{40}}{du}=-f_{43}, \:
\frac{df_{41}}{du}=f_{40}, \:
\frac{df_{42}}{du}=f_{41} ,
\frac{df_{43}}{du}=f_{42} .
\label{c2-30x30x}
\end{equation}

\subsection{The exponential and trigonometric forms of planar
fourcomplex numbers}

Any planar fourcomplex number $u=x+\alpha y+\beta z+\gamma t$ can be writen in
the 
form 
\begin{equation}
x+\alpha y+\beta z+\gamma t=e^{x_1+\alpha y_1+\beta z_1+\gamma t_1} .
\label{c2-31}
\end{equation}
The expressions of $x_1, y_1, z_1, t_1$ as functions of 
$x, y, z, t$ can be obtained by
developing $e^{\alpha y_1}, e^{\beta z_1}$ and $e^{\gamma t_1}$ with the aid of
Eqs. (\ref{c2-29})-(\ref{c2-30b}), by multiplying these expressions and
separating 
the hypercomplex components, and then substituting the expressions of the
four-dimensional cosexponential functions, Eqs.
(\ref{c2-30c30c})-(\ref{c2-30f30f}),  
\begin{equation}
x=e^{x_1}\left(\cos z_1\cos \frac{y_1+t_1}{\sqrt{2}}\cosh
\frac{y_1-t_1}{\sqrt{2}} 
-\sin z_1\sin \frac{y_1+t_1}{\sqrt{2}}\sinh \frac{y_1-t_1}{\sqrt{2}}\right) ,
\label{c2-32}
\end{equation}
\begin{equation}
y=e^{x_1}\left[\sin \left(z_1+\frac{\pi}{4}\right)\cos
\frac{y_1+t_1}{\sqrt{2}}\sinh \frac{y_1-t_1}{\sqrt{2}} 
+\cos \left(z_1+\frac{\pi}{4}\right)\sin \frac{y_1+t_1}{\sqrt{2}}\cosh
\frac{y_1-t_1}{\sqrt{2}}\right] , 
\label{c2-34}
\end{equation}
\begin{equation}
z=e^{x_1}\left(\cos z_1\sin \frac{y_1+t_1}{\sqrt{2}}\sinh
\frac{y_1-t_1}{\sqrt{2}} 
+\sin z_1\cos \frac{y_1+t_1}{\sqrt{2}}\cosh \frac{y_1-t_1}{\sqrt{2}}\right) ,
\label{c2-33}
\end{equation}
\begin{equation}
t=e^{x_1}\left[-\cos \left(z_1+\frac{\pi}{4}\right)\cos
\frac{y_1+t_1}{\sqrt{2}}\sinh \frac{y_1-t_1}{\sqrt{2}} 
+\sin \left(z_1+\frac{\pi}{4}\right)\sin \frac{y_1+t_1}{\sqrt{2}}\cosh
\frac{y_1-t_1}{\sqrt{2}}\right] , 
\label{c2-35}
\end{equation}
The relations (\ref{c2-32})-(\ref{c2-35}) can be rewritten as
\begin{equation}
x+\frac{y-t}{\sqrt{2}}=e^{x_1}\cos\left(z_1
+\frac{y_1+t_1}{\sqrt{2}}\right)e^{(y_1-t_1)/\sqrt{2}},
\label{c2-36}
\end{equation}
\begin{equation}
z+\frac{y+t}{\sqrt{2}}
=e^{x_1}\sin\left(z_1+\frac{y_1+t_1}{\sqrt{2}}\right)e^{(y_1-t_1)/\sqrt{2}},
\label{c2-37}
\end{equation}
\begin{equation}
x-\frac{y-t}{\sqrt{2}}
=e^{x_1}\cos\left(z_1-\frac{y_1+t_1}{\sqrt{2}}\right)e^{-(y_1-t_1)/\sqrt{2}},
\label{c2-38}
\end{equation}
\begin{equation}
z-\frac{y+t}{\sqrt{2}}
=e^{x_1}\sin\left(z_1-\frac{y_1+t_1}{\sqrt{2}}\right)e^{-(y_1-t_1)/\sqrt{2}}.
\label{c2-39}
\end{equation}
By multiplying the sum of the squares of the first two and of the last
two relations (\ref{c2-36})-(\ref{c2-39}) it results that
\begin{equation}
e^{4x_1}=\rho_+^2\rho_-^2 ,
\label{c2-40}
\end{equation}
or
\begin{equation}
e^{x_1}=\rho .
\label{c2-41}
\end{equation}
By summing the squares of all relations (\ref{c2-36})-(\ref{c2-39}) it results
that 
\begin{equation}
d^2=\rho^2\cosh\left[\sqrt{2}(y_1-t_1)\right].
\label{c2-42}
\end{equation}
Then the quantities $y_1, z_1, t_1$ can be expressed in terms of the angles
$\phi, \chi, \psi$ defined in Eqs. (\ref{c2-12a})-(\ref{c2-12d}) as
\begin{equation}
z_1+\frac{y_1+t_1}{\sqrt{2}}=\phi ,
\label{c2-43}
\end{equation}
\begin{equation}
-z_1+\frac{y_1+t_1}{\sqrt{2}}=\chi ,
\label{c2-44}
\end{equation}
\begin{equation}
\frac{e^{(y_1-t_1)/\sqrt{2}}}
{\sqrt{2}\left[\cosh\sqrt{2}(y_1-t_1)\right]^{1/2}}=\sin\psi,\:\:
\frac{e^{-(y_1-t_1)/\sqrt{2}}}
{\sqrt{2}\left[\cosh\sqrt{2}(y_1-t_1)\right]^{1/2}}=\cos\psi.
\label{c2-45}
\end{equation}
From Eq. (\ref{c2-45}) it results that
\begin{equation}
y_1-t_1=\frac{1}{\sqrt{2}}\ln\tan\psi,
\label{c2-45b}
\end{equation}
so that 
\begin{equation}
y_1=\frac{\phi+\chi}{2\sqrt{2}}+\frac{1}{2\sqrt{2}}\ln\tan\psi,\;
z_1=\frac{\phi-\chi}{2},\;
t_1=\frac{\phi+\chi}{2\sqrt{2}}-\frac{1}{2\sqrt{2}}\ln\tan\psi.
\label{c2-45c}
\end{equation}
Substituting the expressions of the quantities $x_1, y_1, z_1, t_1$ in Eq.
(\ref{c2-31}) yields 
\begin{equation}
u=\rho\exp\left[\frac{1}{2\sqrt{2}}(\alpha-\gamma)\ln\tan\psi
+\frac{1}{2}\left(\beta+\frac{\alpha+\gamma}{\sqrt{2}}\right)\phi
-\frac{1}{2}\left(\beta-\frac{\alpha+\gamma}{\sqrt{2}}\right)\chi
\right],
\label{c2-46}
\end{equation}
which will be called the exponential form of the planar fourcomplex number $u$.
It can be checked that
\begin{equation}
\exp\left[\frac{1}{2}\left(\beta+\frac{\alpha+\gamma}
{\sqrt{2}}\right)\phi\right]
=\frac{1}{2}-\frac{\alpha-\gamma}{2\sqrt{2}}
+\left(\frac{1}{2}+\frac{\alpha-\gamma}{2\sqrt{2}}\right)\cos\phi+
\left(\frac{\beta}{2}+\frac{\alpha+\gamma}{2\sqrt{2}}\right)\sin\phi ,
\label{c2-51b}
\end{equation}
\begin{equation}
\exp\left[-\frac{1}{2}\left(\beta-\frac{\alpha+\gamma}
{\sqrt{2}}\right)\chi\right]
=\frac{1}{2}+\frac{\alpha-\gamma}{2\sqrt{2}}
+\left(\frac{1}{2}-\frac{\alpha-\gamma}{2\sqrt{2}}\right)\cos\chi-
\left(\frac{\beta}{2}-\frac{\alpha+\gamma}{2\sqrt{2}}\right)\sin\chi ,
\label{c2-51c}
\end{equation}
which shows that $e^{[\beta+(\alpha+\gamma)/\sqrt{2}]\phi/2}$ and 
$e^{-[\beta-(\alpha+\gamma)/\sqrt{2}]\chi/2}$ are periodic
functions of $\phi$ and respectively $\chi$, with period $2\pi$.

The exponential of the logarithmic term in Eq. (\ref{c2-46}) 
can be expanded with
the aid of the relation (\ref{c2-30l}) as
\begin{equation}
\exp\left[\frac{1}{2\sqrt{2}}(\alpha-\gamma)\ln\tan\psi\right]=
\frac{1}{(\sin2\psi)^{1/2}}\left[\cos\left(\psi-\frac{\pi}{4}\right)
+\frac{\alpha-\gamma}{\sqrt{2}}\sin\left(\psi-\frac{\pi}{4}\right)\right].
\label{c2-47}
\end{equation}
Since according to Eq. (\ref{c2-13c}) $\tan\psi=\rho_+/\rho_-$, then 
\begin{equation}
\sin\psi\cos\psi=\frac{\rho_+\rho_-}{\rho_+^2+\rho_-^2} ,
\label{c2-48}
\end{equation}
and it can be checked that
\begin{equation}
\rho_+^2+\rho_-^2=2d^2,
\label{c2-49}
\end{equation}
where $d$ has been defined in Eq. (\ref{c2-10}). Thus
\begin{equation}
\rho^2=d^2\sin2\psi,
\label{c2-50}
\end{equation}
so that the planar fourcomplex number $u$ can be written as
\begin{equation}
u=d\left[\cos\left(\psi-\frac{\pi}{4}\right)
+\frac{\alpha-\gamma}{\sqrt{2}}\sin\left(\psi-\frac{\pi}{4}\right)\right]
\exp\left[\frac{1}{2}\left(\beta+\frac{\alpha+\gamma}{\sqrt{2}}\right)\phi
-\frac{1}{2}\left(\beta-\frac{\alpha+\gamma}{\sqrt{2}}\right)\chi\right],
\label{c2-51}
\end{equation}
which will be called the trigonometric form of the planar fourcomplex number
$u$. 

If $u_1, u_2$ are planar fourcomplex numbers of moduli and angles $d_1, \phi_1,
\chi_1, \psi_1$ and respectively $d_2, \phi_2, \chi_2, \psi_2$, the product of
the factors depending on the planar angles can be calculated to be
\begin{eqnarray}
\lefteqn{[\cos(\psi_1-\pi/4)+\frac{\alpha-\gamma}{\sqrt{2}}\sin(\psi_1-\pi/4)] 
[\cos(\psi_2-\pi/4)+\frac{\alpha-\gamma}{\sqrt{2}}
\sin(\psi_2-\pi/4)]}\nonumber\\
&&=[\cos(\psi_1-\psi_2)-\frac{\alpha-\gamma}{\sqrt{2}}\cos(\psi_1+\psi_2)] .
\label{c2-52-53}
\end{eqnarray}
The right-hand side of Eq. (\ref{c2-52-53}) can be written as
\begin{eqnarray}
\lefteqn{\cos(\psi_1-\psi_2)-\frac{\alpha-\gamma}
{\sqrt{2}}\cos(\psi_1+\psi_2)}\nonumber\\
&&=[2(\cos^2\psi_1\cos^2\psi_2+\sin^2\psi_1\sin^2\psi_2)]^{1/2}
[\cos(\psi-\pi/4)+\frac{\alpha-\gamma}{\sqrt{2}}\sin(\psi-\pi/4)] ,
\label{c2-54}
\end{eqnarray}
where the angle $\psi$, determined by the condition that
\begin{equation}
\tan(\psi-\pi/4)=-\cos(\psi_1+\psi_2)/\cos(\psi_1-\psi_2)
\label{c2-55}
\end{equation}
is given by $\tan\psi=\tan\psi_1\tan\psi_2$ ,
which is consistent with Eq. (\ref{c2-14b}).
The modulus $d$ of the product $u_1u_2$ is then
\begin{equation}
d=\sqrt{2}d_1d_2
\left(\cos^2\psi_1\cos^2\psi_2+\sin^2\psi_1\sin^2\psi_2\right)^{1/2} .
\label{c2-56}
\end{equation}

\subsection{Elementary functions of planar fourcomplex variables}

The logarithm $u_1$ of the planar fourcomplex number $u$, $u_1=\ln u$, can be
defined 
as the solution of the equation
\begin{equation}
u=e^{u_1} ,
\label{c2-57}
\end{equation}
written explicitly previously in Eq. (\ref{c2-31}), for $u_1$ as a function of
$u$. From Eq. (\ref{c2-46}) it results that 
\begin{equation}
\ln u=\ln \rho+\frac{1}{2\sqrt{2}}(\alpha-\gamma)\ln\tan\psi+
\frac{1}{2}\left(\beta+\frac{\alpha+\gamma}{\sqrt{2}}\right)\phi
-\frac{1}{2}\left(\beta-\frac{\alpha+\gamma}{\sqrt{2}}\right)\chi ,
\label{c2-58}
\end{equation}
which is multivalued because of the presence of the terms proportional to
$\phi$ and $\chi$.
It can be inferred from Eqs. (\ref{c2-14a}) and (\ref{c2-14b}) that
\begin{equation}
\ln(uu^\prime)=\ln u+\ln u^\prime ,
\label{c2-59}
\end{equation}
up to multiples of $\pi[\beta+(\alpha+\gamma)/\sqrt{2}]$ and 
$\pi[\beta-(\alpha+\gamma)/\sqrt{2}]$.

The power function $u^m$ can be defined for real values of $n$ as
\begin{equation}
u^m=e^{m\ln u} .
\label{c2-60}
\end{equation}
The power function is multivalued unless $n$ is an integer. 
For integer $n$, it can be inferred from Eq. (\ref{c2-59}) that
\begin{equation}
(uu^\prime)^m=u^m\:u^{\prime m} .
\label{c2-61}
\end{equation}
If, for example, $m=2$, it can be checked with the aid of Eq. (\ref{c2-51})
that 
Eq. (\ref{c2-60}) gives indeed $(x+\alpha y+\beta z+\gamma t)^2=
x^2-z^2-2yt+2\alpha(xy-zt)+\beta(y^2-t^2+2xz)+2\gamma(xt+yz)$.

The trigonometric functions of the fourcomplex variable
$u$ and the addition theorems for these functions have been written in Eqs.
(\ref{g3})-(\ref{g6}). 
The cosine and sine functions of the hypercomplex variables $\alpha y, 
\beta z$ and $ \gamma t$ can be expressed as
\begin{equation}
\cos\alpha y=f_{40}(y)-\beta f_{42}(y), \: \sin\alpha y=\alpha f_{41}(y)-\gamma
f_{43}(y),  
\label{c2-67}
\end{equation}
\begin{equation}
\cos\beta z=\cosh z, \: \sin\beta z=\beta\sinh z, 
\label{c2-66}
\end{equation}
\begin{equation}
\cos\gamma t=f_{40}(t)+\beta f_{42}(t), \: \sin\gamma t=\gamma f_{41}(t)-\alpha
f_{43}(t) . 
\label{c2-68}
\end{equation}
The cosine and sine functions of a planar fourcomplex number $x+\alpha y+\beta
z+\gamma t$ can then be
expressed in terms of elementary functions with the aid of the addition
theorems Eqs. (\ref{g5}), (\ref{g6}) and of the expressions in  Eqs. 
(\ref{c2-67})-(\ref{c2-68}). 

The hyperbolic functions of the fourcomplex variable
$u$ and the addition theorems for these functions have been written in Eqs.
(\ref{g7})-(\ref{g10}). 
The hyperbolic cosine and sine functions of the hypercomplex variables $\alpha
y,  
\beta z$ and $ \gamma t$ can be expressed as
\begin{equation}
\cosh\alpha y=f_{40}(y)+\beta f_{42}(y), \: \sinh\alpha y=\alpha
f_{41}(y)+\gamma f_{43}(y),  
\label{c2-74}
\end{equation}
\begin{equation}
\cosh\beta z=\cos z, \: \sinh\beta z=\beta\sin z, 
\label{c2-73}
\end{equation}
\begin{equation}
\cosh\gamma t=f_{40}(t)-\beta f_{42}(t), \: \sinh\gamma t=\gamma
f_{41}(t)+\alpha f_{43}(t) . 
\label{c2-75}
\end{equation}
The hyperbolic cosine and sine functions of a planar fourcomplex number
$x+\alpha y+\beta 
z+\gamma t$ can then be
expressed in terms of elementary functions with the aid of the addition
theorems Eqs. (\ref{g9}), (\ref{g10}) and of the expressions in  Eqs. 
(\ref{c2-74})-(\ref{c2-75}).

\subsection{Power series of planar fourcomplex variables}

A planar fourcomplex series is an infinite sum of the form
\begin{equation}
a_0+a_1+a_2+\cdots+a_l+\cdots , 
\label{c2-76}
\end{equation}
where the coefficients $a_l$ are planar fourcomplex numbers. The convergence of
the series (\ref{c2-76}) can be defined in terms of the convergence of its 4
real 
components. The convergence of a planar fourcomplex series can however be
studied 
using planar fourcomplex variables. The main criterion for absolute convergence
remains the comparison theorem, but this requires a number of inequalities
which will be discussed further.

The modulus of a planar fourcomplex number $u=x+\alpha y+\beta z+\gamma t$ can
be defined as 
\begin{equation}
|u|=(x^2+y^2+z^2+t^2)^{1/2} ,
\label{c2-77}
\end{equation}
so that, according to Eq. (\ref{c2-10}), $d=|u|$. Since $|x|\leq |u|, |y|\leq
|u|, 
|z|\leq |u|, |t|\leq |u|$, a property of 
absolute convergence established via a comparison theorem based on the modulus
of the series (\ref{c2-76}) will ensure the absolute convergence of each real
component of that series.

The modulus of the sum $u_1+u_2$ of the planar fourcomplex numbers $u_1, u_2$
fulfils 
the inequality
\begin{equation}
||u_1|-|u_2||\leq |u_1+u_2|\leq |u_1|+|u_2| .
\label{c2-78}
\end{equation}
For the product the relation is 
\begin{equation}
|u_1u_2|\leq \sqrt{2}|u_1||u_2| ,
\label{c2-79}
\end{equation}
as can be shown from Eq. (\ref{c2-56}). The relation (\ref{c2-79}) replaces the
relation of equality extant for regular complex numbers. 
The equality in Eq. (\ref{c2-79}) takes place for 
$\cos^2(\psi_1-\psi_2)=1, \:\cos^2(\psi_1+\psi_2)=1,$ which means that
$x_1+(y_1-t_1)/\sqrt{2}=0, \:z_1+(y_1+t_1)/\sqrt{2}=0,
\:x_2+(y_2-t_2)/\sqrt{2}=0, \:z_2+(y_2+t_2)/\sqrt{2}=0,$ or
$x_1-(y_1-t_1)/\sqrt{2}=0, \:z_1-(y_1+t_1)/\sqrt{2}=0,
\:x_2-(y_2-t_2)/\sqrt{2}=0, \:z_2-(y_2+t_2)/\sqrt{2}=0$.
The modulus of the product, which has the property that $0\leq |u_1u_2|$,
becomes equal to zero for
$\cos^2(\psi_1-\psi_2)=0, \:\cos^2(\psi_1+\psi_2)=0,$ which means that
$x_1+(y_1-t_1)/\sqrt{2}=0, \:z_1+(y_1+t_1)/\sqrt{2}=0,
\:x_2-(y_2-t_2)/\sqrt{2}=0, \:z_2-(y_2+t_2)/\sqrt{2}=0$,
or
$x_1-(y_1-t_1)/\sqrt{2}=0, \:z_1-(y_1+t_1)/\sqrt{2}=0,
\:x_2+(y_2-t_2)/\sqrt{2}=0, \:z_2+(y_2+t_2)/\sqrt{2}=0$.
as discussed after Eq. (\ref{c2-9}).

It can be shown that
\begin{equation}
x^2+y^2+z^2+t^2\leq|u^2|\leq \sqrt{2}(x^2+y^2+z^2+t^2) .
\label{c2-80}
\end{equation}
The left relation in Eq. (\ref{c2-80}) becomes an equality for $\sin^2
2\psi=1$, 
when $\rho_+=\rho_-$, which means that $x(y-t)+z(y+t)=0$.
The right relation in Eq. (\ref{c2-80}) becomes an equality for $\sin^2
2\psi=0$, 
when $x+(y-t)/\sqrt{2}=0, \:z+(y+t)/\sqrt{2}=0,$
or $x-(y-t)/\sqrt{2}=0, \:z-(y+t)/\sqrt{2}=0.$
The inequality in Eq. (\ref{c2-79}) implies that
\begin{equation}
|u^l|\leq 2^{(l-1)/2}|u|^l .
\label{c2-81}
\end{equation}
From Eqs. (\ref{c2-79}) and (\ref{c2-81}) it results that
\begin{equation}
|au^l|\leq 2^{l/2} |a| |u|^l .
\label{c2-82}
\end{equation}

A power series of the planar fourcomplex variable $u$ is a series of the form
\begin{equation}
a_0+a_1 u + a_2 u^2+\cdots +a_l u^l+\cdots .
\label{c2-83}
\end{equation}
Since
\begin{equation}
\left|\sum_{l=0}^\infty a_l u^l\right| \leq  \sum_{l=0}^\infty
2^{l/2}|a_l| |u|^l ,
\label{c2-84}
\end{equation}
a sufficient condition for the absolute convergence of this series is that
\begin{equation}
\lim_{l\rightarrow \infty}\frac{\sqrt{2}|a_{l+1}||u|}{|a_l|}<1 .
\label{c2-85}
\end{equation}
Thus the series is absolutely convergent for 
\begin{equation}
|u|<c,
\label{c2-86}
\end{equation}
where 
\begin{equation}
c=\lim_{l\rightarrow\infty} \frac{|a_l|}{\sqrt{2}|a_{l+1}|} .
\label{c2-87}
\end{equation}

The convergence of the series (\ref{c2-83}) can be also studied with the aid of
the transformation 
\begin{equation}
x+\alpha y+\beta z+\gamma t=\sqrt{2}(e_1\xi+\tilde e_1 \upsilon+e_2\tau
+\tilde e_2\zeta) , 
\label{c2-87b}
\end{equation}
where $\xi,\upsilon, \tau, \zeta$ have been defined in Eq. (\ref{c2-11}),
and
\begin{equation}
e_1=\frac{1}{2}+\frac{\alpha-\gamma}{2\sqrt{2}},\:\:
\tilde e_1=\frac{\beta}{2}+\frac{\alpha+\gamma}{2\sqrt{2}},\:\:
e_2=\frac{1}{2}-\frac{\alpha-\gamma}{2\sqrt{2}},\:\:
\tilde e_2=-\frac{\beta}{2}+\frac{\alpha+\gamma}{2\sqrt{2}}.
\label{c2-87c}
\end{equation}
It can be checked that
\begin{eqnarray}
\lefteqn{e_1^2=e_1, \:\:\tilde e_1^2=-e_1,\:\: e_1\tilde e_1=\tilde e_1,\:\:
e_2^2=e_2, \:\:\tilde e_2^2=-e_2,\:\: e_2\tilde e_2=\tilde e_2,\:\:\nonumber}\\
&&e_1e_2=0,\:\: \tilde e_1\tilde e_2=0, \:\:e_1\tilde e_2=0, \:\:
e_2\tilde e_1=0.  
\label{c2-87d}
\end{eqnarray}
The moduli of the bases in Eq. (\ref{c2-87c}) are
\begin{equation}
|e_1|=\frac{1}{\sqrt{2}},\;|\tilde e_1|=\frac{1}{\sqrt{2}},\;
|e_2|=\frac{1}{\sqrt{2}},\;|\tilde e_2|=\frac{1}{\sqrt{2}},
\label{c2-87e}
\end{equation}
and it can be checked that
\begin{equation}
|x+\alpha y+\beta z+\gamma t|^2=\xi^2+\upsilon^2+\tau^2+\zeta^2.
\label{c2-87f}
\end{equation}
The ensemble $e_1, \tilde e_1, e_2, \tilde e_2$ will be called the canonical
planar fourcomplex base, and Eq. (\ref{c2-87b}) gives the canonical form of the
planar fourcomplex number.

If $u=u^\prime u^{\prime\prime}$, the components $\xi,\upsilon, \tau, \zeta$
are related, according to Eqs. (\ref{c2-17})-(\ref{c2-20}) by 
\begin{equation}
\xi=\sqrt{2}(\xi^\prime \xi^{\prime\prime}
-\upsilon^\prime \upsilon^{\prime\prime}), \:\:
\upsilon=\sqrt{2}(\xi^\prime \upsilon^{\prime\prime}
+\upsilon^\prime \xi^{\prime\prime}), \:\:
\tau=\sqrt{2}(\tau^\prime \tau^{\prime\prime}
-\zeta^\prime \zeta^{\prime\prime}), \:\:
\zeta=\sqrt{2}(\tau^\prime \zeta^{\prime\prime}
+\zeta^\prime \xi^{\prime\prime}), \:\:
\label{c2-87g}
\end{equation}
which show that, upon multiplication, the components $\xi,\upsilon$ and $\tau,
\zeta$ obey, up to a normalization constant, the same
rules as the real and imaginary components of usual, two-dimensional complex
numbers.

If the coefficients in Eq. (\ref{c2-83}) are 
\begin{equation}
a_l= a_{l0}+\alpha a_{l1}+\beta a_{l2}+\gamma a_{l3}, 
\label{c2-n88a}
\end{equation}
and
\begin{equation}
A_{l1}=a_{l0}+\frac{a_{l1}-a_{l3}}{\sqrt{2}},\;
\tilde A_{l1}=a_{l2}+\frac{a_{l1}+a_{l3}}{\sqrt{2}},\;
A_{l2}=a_{l0}-\frac{a_{l1}-a_{l3}}{\sqrt{2}},\;
\tilde A_{l2}=-a_{l2}+\frac{a_{l1}+a_{l3}}{\sqrt{2}},
\label{c2-n88b}
\end{equation}
the series (\ref{c2-83}) can be written as
\begin{equation}
\sum_{l=0}^\infty 2^{l/2}\left[
(e_1 A_{l1}+\tilde e_1\tilde A_{l1})(e_1 \xi+\tilde e_1\upsilon)^l 
+(e_2 A_{l2}+\tilde e_2\tilde A_{l2})(e_2 \tau+\tilde e_2\zeta)^l 
\right].
\label{c2-n89a}
\end{equation}
Thus, the series in Eqs. (\ref{c2-83}) and (\ref{c2-n89a}) are
absolutely convergent for   
\begin{equation}
\rho_+<c_1, \;\rho_-<c_2,
\label{c2-n90}
\end{equation}
where 
\begin{equation}
c_1=\lim_{l\rightarrow\infty} \frac
{\left[A_{l1}^2+\tilde A_{l1}^2\right]^{1/2}}
{\sqrt{2}\left[A_{l+1,1}^2+\tilde A_{l+1,1}^2\right]^{1/2}},\;\;
c_2=\lim_{l\rightarrow\infty} \frac
{\left[A_{l2}^2+\tilde A_{l2}^2\right]^{1/2}}
{\sqrt{2}\left[A_{l+1,2}^2+\tilde A_{l+1,2}^2\right]^{1/2}}.
\label{c2-n91}
\end{equation}

It can be shown that $c=(1/\sqrt{2}){\rm min}(c_1,c_2)$, where ${\rm min}$
designates the smallest of the numbers $c_1,c_2$. Using the expression of $|u|$
in Eq. (\ref{c2-87f}),  it can be seen that the spherical region of convergence
defined in Eqs. (\ref{c2-86}), (\ref{c2-87}) is included in the cylindrical
region of convergence defined in Eqs. (\ref{c2-n90}) and (\ref{c2-n91}).

\subsection{Analytic functions of planar fourcomplex variables}

The fourcomplex function $f(u)$ of the fourcomplex variable $u$ has
been expressed in Eq. (\ref{g16}) in terms of 
the real functions $P(x,y,z,t),Q(x,y,z,t),R(x,y,z,t), S(x,y,z,t)$ of real
variables $x,y,z,t$.
The
relations between the partial derivatives of the functions $P, Q, R, S$ are
obtained by setting succesively in   
Eq. (\ref{g17}) $\Delta x\rightarrow 0, \Delta y=\Delta z=\Delta t=0$;
then $ \Delta y\rightarrow 0, \Delta x=\Delta z=\Delta t=0;$  
then $  \Delta z\rightarrow 0,\Delta x=\Delta y=\Delta t=0$; and finally
$ \Delta t\rightarrow 0,\Delta x=\Delta y=\Delta z=0 $. 
The relations are 
\begin{equation}
\frac{\partial P}{\partial x} = \frac{\partial Q}{\partial y} =
\frac{\partial R}{\partial z} = \frac{\partial S}{\partial t},
\label{c2-95}
\end{equation}
\begin{equation}
\frac{\partial Q}{\partial x} = \frac{\partial R}{\partial y} =
\frac{\partial S}{\partial z} = -\frac{\partial P}{\partial t},
\label{c2-97}
\end{equation}
\begin{equation}
\frac{\partial R}{\partial x} = \frac{\partial S}{\partial y} =
-\frac{\partial P}{\partial z} = -\frac{\partial Q}{\partial t},
\label{c2-96}
\end{equation}
\begin{equation}
\frac{\partial S}{\partial x} =-\frac{\partial P}{\partial y} =
-\frac{\partial Q}{\partial z} =-\frac{\partial R}{\partial t}.
\label{c2-98}
\end{equation}
The relations (\ref{c2-95})-(\ref{c2-98}) are analogous to the Riemann
relations for the real and imaginary components of a complex function. It can
be shown from Eqs. (\ref{c2-95})-(\ref{c2-98}) that the component $P$ is a
solution of the equations
\begin{equation}
\frac{\partial^2 P}{\partial x^2}+\frac{\partial^2 P}{\partial z^2}=0,
\:\: 
\frac{\partial^2 P}{\partial y^2}+\frac{\partial^2 P}{\partial t^2}=0,
\:\:
\label{c2-99}
\end{equation}
and the components $Q, R, S$ are solutions of similar equations.
As can be seen from Eqs. (\ref{c2-99}), the components $P, Q, R, S$ of
an analytic function of planar fourcomplex variable are harmonic 
with respect to the pairs of variables $x,y$ and $ z,t$.
The component $P$ is also a solution of the mixed-derivative
equations
\begin{equation}
\frac{\partial^2 P}{\partial x^2}=-\frac{\partial^2 P}{\partial y\partial t},
\:\: 
\frac{\partial^2 P}{\partial y^2}=\frac{\partial^2 P}{\partial x\partial z},
\:\:
\frac{\partial^2 P}{\partial z^2}=\frac{\partial^2 P}{\partial y\partial t},
\:\:
\frac{\partial^2 P}{\partial t^2}=-\frac{\partial^2 P}{\partial x\partial z},
\:\:
\label{c2-106b}
\end{equation}
and the components $Q, R, S$ are solutions of similar equations.
The component $P$ is also a solution of the mixed-derivative
equations 
\begin{equation}
\frac{\partial^2 P}{\partial x\partial y}=-\frac{\partial^2 P}{\partial
z\partial t} ,
\:\: 
\frac{\partial^2 P}{\partial x\partial t}=\frac{\partial^2 P}{\partial
y\partial z} ,
\label{c2-107}
\end{equation}
and the components $Q, R, S$ are solutions of similar equations.

\subsection{Integrals of functions of planar fourcomplex variables}

The singularities of planar fourcomplex functions arise from terms of the form
$1/(u-u_0)^m$, with $m>0$. Functions containing such terms are singular not
only at $u=u_0$, but also at all points of the two-dimensional hyperplanes
passing through $u_0$ and which are parallel to the nodal hyperplanes. 

The integral of a planar fourcomplex function between two points $A, B$ along a
path situated in a region free of singularities is independent of path, which
means that the integral of an analytic function along a loop situated in a
region free from singularities is zero,
\begin{equation}
\oint_\Gamma f(u) du = 0,
\label{c2-111}
\end{equation}
where it is supposed that a surface $\Sigma$ spanning 
the closed loop $\Gamma$ is not intersected by any of
the two-dimensional hyperplanes associated with the
singularities of the function $f(u)$. Using the expression, Eq. (\ref{g16})
for $f(u)$ and the fact that $du=dx+\alpha  dy+\beta dz+\gamma dt$, the
explicit form of the integral in Eq. (\ref{c2-111}) is
\begin{eqnarray}
\lefteqn{\oint _\Gamma f(u) du = \oint_\Gamma
[(Pdx-Sdy-Rdz-Qdt)+\alpha(Qdx+Pdy-Sdz-Rdt)\nonumber}\\
&&+\beta(Rdx+Qdy+Pdz-Sdt)+\gamma(Sdx+Rdy+Qdz+Pdt)] .
\label{c2-112}
\end{eqnarray}
If the functions $P, Q, R, S$ are regular on a surface $\Sigma$
spanning the loop $\Gamma$,
the integral along the loop $\Gamma$ can be transformed with the aid of the
theorem of Stokes in an integral over the surface $\Sigma$ of terms of the form
$\partial P/\partial y + \partial S/\partial x,\:\:
\partial P/\partial z +  \partial R/\partial x, \:\:
\partial P/\partial t + \partial Q/\partial x, \:\:
\partial R/\partial y -  \partial S/\partial z, \:\:
\partial S/\partial t - \partial Q/\partial y, \:\:
\partial R/\partial t - \partial Q/\partial z$ 
and of similar terms arising
from the $\alpha, \beta$ and $\gamma$ components, 
which are equal to zero by Eqs. (\ref{c2-95})-(\ref{c2-98}), and this proves
Eq.  (\ref{c2-111}).

The integral of the function $(u-u_0)^m$ on a closed loop $\Gamma$ is equal to
zero for $m$ a positive or negative integer not equal to -1,
\begin{equation}
\oint_\Gamma (u-u_0)^m du = 0, \:\: m \:\:{\rm integer},\: m\not=-1 .
\label{c2-112b}
\end{equation}
This is due to the fact that $\int (u-u_0)^m du=(u-u_0)^{m+1}/(m+1), $ and to
the fact that the function $(u-u_0)^{m+1}$ is singlevalued for $m$ an integer.

The integral $\oint du/(u-u_0)$ can be calculated using the exponential form 
(\ref{c2-46}),
\begin{eqnarray}
u-u_0=\rho\exp\left[\frac{1}{2\sqrt{2}}(\alpha-\gamma)\ln\tan\psi
+\frac{1}{2}\left(\beta+\frac{\alpha+\gamma}{\sqrt{2}}\right)\phi
-\frac{1}{2}\left(\beta-\frac{\alpha+\gamma}{\sqrt{2}}\right)\chi\right],
\label{c2-113}
\end{eqnarray}
so that 
\begin{equation}
\frac{du}{u-u_0}=\frac{d\rho}{\rho}+
\frac{1}{2\sqrt{2}}(\alpha-\gamma)d\ln\tan\psi
+\frac{1}{2}\left(\beta+\frac{\alpha+\gamma}{\sqrt{2}}\right)d\phi
-\frac{1}{2}\left(\beta-\frac{\alpha+\gamma}{\sqrt{2}}\right)d\chi .
\label{c2-114}
\end{equation}
Since $\rho$ and $\psi$ are singlevalued variables, it follows that
$\oint_\Gamma d\rho/\rho =0, \oint_\Gamma d\ln\tan\psi=0$. On the other hand,
$\phi$ and $\chi$ are cyclic variables, so that they may give a contribution to
the integral around the closed loop $\Gamma$.
Thus, if $C_+$ is a circle of radius $r$
parallel to the $\xi O\upsilon$ plane, whose
projection of the center of this circle on the $\xi O\upsilon$ plane
coincides with the projection of the point $u_0$ on this plane, the points
of the circle $C_+$ are described according to Eqs.
(\ref{c2-11})-(\ref{c2-12d}) by the equations
\begin{eqnarray}
\lefteqn{\xi=\xi_0+r \sin\psi\cos\phi , \:
\upsilon=\upsilon_0+r \sin\psi\sin\phi , \:
\tau=\tau_0+r\cos\psi \cos\chi , \nonumber}\\
&&\zeta=\zeta_0+r \cos\psi\sin\chi , 
\label{c2-115}
\end{eqnarray}
for constant values of $\chi$ and $\psi, \:\psi\not=0, \pi/2$, where
$u_0=x_0+\alpha y_0+\beta 
z_0+\gamma t_0$,  and $\xi_0, \upsilon_0, \tau_0, \zeta_0$ are calculated from
$x_0, y_0, z_0, t_0$ according to Eqs. (\ref{c2-11}).
Then
\begin{equation}
\oint_{C_+}\frac{du}{u-u_0}
=\pi\left(\beta+\frac{\alpha+\gamma}{\sqrt{2}}\right). 
\label{c2-116}
\end{equation}
If $C_-$ is a circle of radius $r$
parallel to the $\tau O\zeta$ plane,
whose projection of the center of this circle on the $\tau O\zeta$ plane
coincides with the projection of the point $u_0$ on this plane, the points
of the circle $C_-$ are described by the same Eqs. (\ref{c2-115}) 
but for constant values of $\phi$ and $\psi, \:\psi\not=0, \pi/2$.
Then
\begin{equation}
\oint_{C_-}\frac{du}{u-u_0}
=-\pi\left(\beta-\frac{\alpha+\gamma}{\sqrt{2}}\right) .
\label{c2-117}
\end{equation}
The expression of $\oint_\Gamma du/(u-u_0)$ can be written as a single equation
with the aid of the functional int($M,C$) defined in Eq. (\ref{118}) as
\begin{equation}
\oint_\Gamma\frac{du}{u-u_0}=
\pi\left(\beta+\frac{\alpha+\gamma}{\sqrt{2}}\right) \;{\rm
int}(u_{0\xi\upsilon},\Gamma_{\xi\upsilon}) 
-\pi \left(\beta-\frac{\alpha+\gamma}{\sqrt{2}}\right)\;{\rm
int}(u_{0\tau\zeta},\Gamma_{\tau\zeta}), 
\label{c2-119}
\end{equation}
where $u_{0\xi\upsilon}, u_{0\tau\zeta}$ and $\Gamma_{\xi\upsilon},
\Gamma_{\tau\zeta}$ are respectively the projections of the point $u_0$ and of
the loop $\Gamma$ on the planes $\xi \upsilon$ and $\tau \zeta$.

If $f(u)$ is an analytic planar fourcomplex function which can be expanded in a
series as written in Eq. (\ref{g14}), and the expansion holds on the curve
$\Gamma$ and on a surface spanning $\Gamma$, then from Eqs. (\ref{c2-112b}) and
(\ref{c2-119}) it follows that
\begin{equation}
\oint_\Gamma \frac{f(u)du}{u-u_0}=
\pi\left[\left(\beta+\frac{\alpha+\gamma}{\sqrt{2}}\right) \;{\rm
int}(u_{0\xi\upsilon},\Gamma_{\xi\upsilon}) 
- \left(\beta-\frac{\alpha+\gamma}{\sqrt{2}}\right)\;{\rm
int}(u_{0\tau\zeta},\Gamma_{\tau\zeta})\right]\;f(u_0) , 
\label{c2-120}
\end{equation}
where $\Gamma_{\xi\upsilon}, \Gamma_{\tau\zeta}$ are the projections of 
the curve $\Gamma$ on the planes $\xi \upsilon$ and respectively $\tau \zeta$,
as shown in Fig. 2. As remarked previously,
the definition of the variables in this section is different from the
former definition for the circular hypercomplex numbers. 

Substituting in the right-hand side of 
Eq. (\ref{c2-120}) the expression of $f(u)$ in terms of the real 
components $P, Q, R, S$, Eq. (\ref{g16}), yields
\begin{eqnarray}
\lefteqn{\oint_\Gamma \frac{f(u)du}{u-u_0}\nonumber}\\
&&=\pi \left[
\left(\beta+\frac{\alpha+\gamma}{\sqrt{2}}\right)P
-\left(1+\frac{\alpha-\gamma}{\sqrt{2}}\right)R
-\left(\gamma-\frac{1-\beta}{\sqrt{2}}\right)Q
-\left(\alpha-\frac{1+\beta}{\sqrt{2}}\right)S
\right] 
{\rm int}\left(u_{0\xi\upsilon},\Gamma_{\xi\upsilon}\right)\nonumber\\
&&-\pi \left[
\left(\beta-\frac{\alpha+\gamma}{\sqrt{2}}\right)P
-\left(1-\frac{\alpha-\gamma}{\sqrt{2}}\right)R
-\left(\gamma-\frac{1+\beta}{\sqrt{2}}\right)Q
-\left(\alpha-\frac{1-\beta}{\sqrt{2}}\right)S
\right] 
{\rm int}(u_{0\tau\zeta},\Gamma_{\tau\zeta}) ,\nonumber\\
&&
\label{c2-121}
\end{eqnarray}
where $P, Q, R, S$ are the values of the components of $f$ at $u=u_0$.

If $f(u)$ can be expanded as written in Eq. (\ref{g14}) on 
$\Gamma$ and on a surface spanning $\Gamma$, then from Eqs. (\ref{c2-112b}) and
(\ref{c2-119}) it also results that
\begin{equation}
\oint_\Gamma \frac{f(u)du}{(u-u_0)^{m+1}}=
\frac{\pi}{m!}\left[\left(\beta+\frac{\alpha+\gamma}{\sqrt{2}}\right) \;{\rm
int}(u_{0\xi\upsilon},\Gamma_{\xi\upsilon}) 
- \left(\beta-\frac{\alpha+\gamma}{\sqrt{2}}\right)\;{\rm
int}(u_{0\tau\zeta},\Gamma_{\tau\zeta})\right]\; 
f^{(m)}(u_0) ,
\label{c2-122}
\end{equation}
where it has been used the fact that the derivative $f^{(m)}(u_0)$ of order $n$
of $f(u)$ at $u=u_0$ is related to the expansion coefficient in Eq. (\ref{g14})
according to Eq. (\ref{g15}).

If a function $f(u)$ is expanded in positive and negative powers of $u-u_j$,
where $u_j$ are planar fourcomplex constants, $j$ being an index, the integral
of $f$ on a closed loop $\Gamma$ is determined by the terms in the expansion of
$f$ which are of the form $a_j/(u-u_j)$,
\begin{equation}
f(u)=\cdots+\sum_j\frac{a_j}{u-u_j}+\cdots
\label{c2-123}
\end{equation}
Then the integral of $f$ on a closed loop $\Gamma$ is
\begin{equation}
\oint_\Gamma f(u) du = 
\pi\left(\beta+\frac{\alpha+\gamma}{\sqrt{2}}\right) \sum_j{\rm
int}(u_{j\xi\upsilon},\Gamma_{\xi\upsilon})a_j 
- \pi\left(\beta-\frac{\alpha+\gamma}{\sqrt{2}}\right)\sum_j{\rm
int}(u_{j\tau\zeta},\Gamma_{\tau\zeta})a_j .
\label{c2-124}
\end{equation}

\subsection{Factorization of planar fourcomplex polynomials}

A polynomial of degree $m$ of the planar fourcomplex variable 
$u=x+\alpha y+\beta z+\gamma t$ has the form
\begin{equation}
P_m(u)=u^m+a_1 u^{m-1}+\cdots+a_{m-1} u +a_m ,
\label{c2-125}
\end{equation}
where the constants are in general planar fourcomplex numbers.

It can be shown that any planar fourcomplex polynomial has a planar fourcomplex
root, whence it follows that a polynomial of degree $m$ can be written as a
product of $m$ linear factors of the form $u-u_j$, where the planar fourcomplex
numbers $u_j$ are the roots of the polynomials, although the factorization may
not be unique,
\begin{equation}
P_m(u)=\prod_{j=1}^m (u-u_j) .
\label{c2-126}
\end{equation}

The fact that any planar fourcomplex polynomial has a root can be shown by
considering the transformation of a fourdimensional sphere with the center at
the origin by the function $u^m$. The points of the hypersphere of radius $d$
are of the form written in Eq. (\ref{c2-51}), with $d$ constant and $\phi,
\chi, \psi$ arbitrary. The point $u^m$ is
\begin{eqnarray}
\lefteqn{u^m=d^m\left[\cos\left(\psi-\frac{\pi}{4}\right)
+\frac{\alpha-\gamma}{\sqrt{2}}
\sin\left(\psi-\frac{\pi}{4}\right)\right]^m\nonumber}\\
&&\exp\left[\frac{1}{2}\left(\beta+\frac{\alpha+\gamma}{\sqrt{2}}\right)m\phi
-\frac{1}{2}\left(\beta-\frac{\alpha+\gamma}{\sqrt{2}}\right)m\chi\right].
\label{c2-127}
\end{eqnarray}
It can be shown with the aid of Eq. (\ref{c2-56}) that
\begin{equation}
\left|u
\exp\left[\frac{1}{2}\left(\beta+\frac{\alpha+\gamma}{\sqrt{2}}\right)\phi
-\frac{1}{2}\left(\beta-\frac{\alpha+\gamma}{\sqrt{2}}\right)\chi\right]\right|
=|u|,
\label{c2-127bb}
\end{equation}
so that
\begin{eqnarray}
\lefteqn{\left|
\left[\cos(\psi-\pi/4)+\frac{\alpha-\gamma}\sin(\psi-\pi/4)\right]^m
\exp\left[\frac{1}{2}\left(\beta+\frac{\alpha+\gamma}{\sqrt{2}}\right)m\phi
-\frac{1}{2}\left(\beta-\frac{\alpha+\gamma}{\sqrt{2}}\right)m\chi\right]
\right|\nonumber}\\
&&=\left|\left(\cos(\psi-\pi/4)
+\frac{\alpha-\gamma}{\sqrt{2}}\sin(\psi-\pi/4)\right)^m\right| .
\label{c2-127b}
\end{eqnarray}
The right-hand side of Eq. (\ref{c2-127b}) is
\begin{equation}
\left|\left(\cos\epsilon
+\frac{\alpha-\gamma}{\sqrt{2}}\sin\epsilon\right)^m\right|^2
=\sum_{k=0}^m C_{2m}^{2k}\cos^{2m-2k}\epsilon\sin^{2k}\epsilon ,
\label{c2-128}
\end{equation}
where $\epsilon=\psi-\pi/4$, 
and since $C_{2m}^{2k}\geq C_m^k$, it can be concluded that
\begin{equation}
\left|\left(\cos\epsilon
+\frac{\alpha-\gamma}{\sqrt{2}}\sin\epsilon\right)^m\right|^2\geq 1 .
\label{c2-129}
\end{equation}
Then
\begin{equation}
d^m\leq |u^m|\leq 2^{(m-1)/2} d^m ,
\label{c2-129b}
\end{equation}
which shows that the image of a four-dimensional sphere via the transformation
operated by the function $u^m$ is a finite hypersurface.

If $u^\prime=u^m$, and
\begin{equation}
u^\prime=d^\prime
\left[\cos(\psi^\prime-\pi/4)
+\frac{\alpha-\gamma}{\sqrt{2}}\sin(\psi^\prime-\pi/4)\right]
\exp\left[
\frac{1}{2}\left(\beta+\frac{\alpha+\gamma}{\sqrt{2}}\right)\phi^\prime
-\frac{1}{2}\left(\beta-\frac{\alpha+\gamma}{\sqrt{2}}\right)\chi^\prime
\right],
\label{c2-130}
\end{equation}
then 
\begin{equation}
\phi^\prime=m\phi, \: \chi^\prime=m\chi, \: \tan\psi^\prime=\tan^m\psi .
\label{c2-131}
\end{equation}
Since for any values of the angles $\phi^\prime, \chi^\prime, \psi^\prime$
there is a set of solutions $\phi, \chi, \psi$ of Eqs. (\ref{c2-131}), and
since the image of the hypersphere is a finite hypersurface, it follows that
the image of the four-dimensional sphere via the function $u^m$ is also a
closed hypersurface. A continuous hypersurface is called closed when any ray
issued from the origin intersects that surface at least once in the finite part
of the space.

A transformation of the four-dimensional space by the polynomial $P_m(u)$
will be considered further. By this transformation, a hypersphere of radius $d$
having the center at the origin is changed into a certain finite closed
surface, as discussed previously. 
The transformation of the four-dimensional space by the polynomial $P_m(u)$
associates to the point $u=0$ the point $f(0)=a_m$, and the image of a
hypersphere of very large radius $d$ can be represented with good approximation
by the image of that hypersphere by the function $u^m$. 
The origin of the axes is an inner
point of the latter image. If the radius of the hypersphere is now reduced
continuously from the initial very large values to zero, the image hypersphere
encloses initially the origin, but the image shrinks to $a_m$ when the radius
approaches the value zero.  Thus, the
origin is initially inside the image hypersurface, and it lies outside the
image hypersurface when the radius of the hypersphere tends to zero. Then since
the image hypersurface is closed, the image surface must intersect at some
stage the origin of the axes, which means that there is a point $u_1$ such that
$f(u_1)=0$. The factorization in Eq. (\ref{c2-126}) can then be obtained by
iterations.

The roots of the polynomial $P_m$ can be obtained by the following method.
If the constants in Eq. (\ref{c2-125}) are $a_l=a_{l0}+\alpha a_{l1}
+\beta a_{l2}+\gamma a_{l3}$, and with the 
notations of Eq. (\ref{c2-n88b}), the polynomial $P_m(u)$ can be written as
\begin{eqnarray}
\lefteqn{P_m=\sum_{l=0}^{m} 2^{(m-l)/2}
(e_1 A_{l1}+\tilde e_1\tilde A_{l1})
(e_1 \xi+\tilde e_1 \upsilon)^{m-l}\nonumber}\\
&&+\sum_{l=0}^{m} 2^{(m-l)/2}
(e_2 A_{l2}+\tilde e_2\tilde A_{l2})(e_2 \tau+\tilde e_2 \zeta)^{m-l} ,
\label{c2-126a}
\end{eqnarray}
where the constants $A_{lk}, \tilde A_{lk}, k=1,2$ are real numbers.
Each of the polynomials of degree $m$ in $e_1 \xi+\tilde e_1\upsilon, 
e_2 \tau+\tilde e_2\zeta$
in Eq. (\ref{c2-126a}) 
can always be written as a product of linear factors of the form
$e_1 (\xi-\xi_p)+\tilde e_1(\upsilon- \upsilon_p)$ and respectively
$e_2 (\tau-\tau_p)+\tilde e_2(\zeta- \zeta_p)$, where the
constants $\xi_p, \upsilon_p, \tau_p, \zeta_p$ are real,
\begin{eqnarray}
\sum_{l=0}^{m} 2^{(m-l)/2}
(e_1 A_{l1}+\tilde e_1\tilde A_{l1})(e_1 \xi+\tilde e_1 \upsilon)^{m-l}
=\prod_{p=1}^{m}2^{m/2}\left\{e_1 (\xi-\xi_p)+\tilde e_1(\upsilon- \upsilon_p)
\right\},
\label{c2-126bb}
\end{eqnarray}
\begin{eqnarray}
\sum_{l=0}^{m} 2^{(m-l)/2}
(e_2 A_{l2}+\tilde e_2\tilde A_{l2})(e_2 \tau+\tilde e_2 \zeta)^{m-l}
=\prod_{p=1}^{m}2^{m/2}\left\{e_2 (\tau-\tau_p)+\tilde e_2(\zeta- \zeta_p)
\right\}.
\label{c2-126bc}
\end{eqnarray}

Due to the relations  (\ref{c2-87d}),
the polynomial $P_m(u)$ can be written as a product of factors of
the form 
\begin{eqnarray}
P_m(u)=\prod_{p=1}^m 2^{m/2}
\left\{e_1 (\xi-\xi_p)+\tilde e_1(\upsilon- \upsilon_p)
+e_2 (\tau-\tau_p)+\tilde e_2(\zeta- \zeta_p)\right\}.
\label{c2-128b}
\end{eqnarray}
This relation can be written with the aid of Eq. (\ref{c2-87b}) in the form 
(\ref{c2-126}), where
\begin{eqnarray}
u_p=\sqrt{2}(e_1 \xi_p+\tilde e_1 \upsilon_p
+e_2 \tau_p+\tilde e_2 \zeta_p) .
\label{c2-129bx}
\end{eqnarray}
The roots $e_1 \xi_p+\tilde e_1 \upsilon_p$ and $e_2 \tau_p+\tilde e_2 \zeta_p$
defined in Eqs. (\ref{c2-126bb}) and respectively (\ref{c2-126bc}) may be
ordered arbitrarily. This means that Eq.  (\ref{c2-129bx}) gives sets of $m$
roots $u_1,...,u_m$ of the polynomial $P_m(u)$, corresponding to the various
ways in which the roots $e_1 \xi_p+\tilde e_1 \upsilon_p$ and $e_2
\tau_p+\tilde e_2 \zeta_p$ are ordered according to $p$ for each polynomial.
Thus, while the hypercomplex components in Eqs. (\ref{c2-126bb}),
(\ref{c2-126bc}) taken separately have unique factorizations, the polynomial
$P_m(u)$ can be written in many different ways as a product of linear factors.
The result of the planar fourcomplex integration, Eq. (\ref{c2-124}), is
however unique.

If, for example, $P(u)=u^2+1$, the possible factorizations are
$P=(u-\tilde e_1-\tilde e_2)(u+\tilde e_1+\tilde e_2)$ and
$P=(u-\tilde e_1+\tilde e_2)(u+\tilde e_1-\tilde e_2)$ which can also be
written 
as $u^2+1=(u-\beta)(u+\beta)$ or as
$u^2+1=\left\{u-(\alpha+\gamma)/\sqrt{2}\right\}
\left\{(u+(\alpha+\gamma)/\sqrt{2}\right\}$. The result of the planar
fourcomplex 
integration, Eq. (\ref{c2-124}), is however unique. 
It can be checked
that $(\pm \tilde e_1\pm\tilde e_2)^2=
-e_1-e_2=-1$.

\subsection{Representation of planar fourcomplex numbers 
by irreducible matrices}

If $T$ is the unitary matrix,
\begin{equation}
T =\left(
\begin{array}{cccc}
\frac{1}{\sqrt{2}}&\frac{1}{2}       &0                  &-\frac{1}{2}  \\
0                 &\frac{1}{2}       &\frac{1}{\sqrt{2}} &\frac{1}{2} \\
\frac{1}{\sqrt{2}}&-\frac{1}{2}      &0                  &\frac{1}{2}  \\
0                 &\frac{1}{2}       &-\frac{1}{\sqrt{2}}&\frac{1}{2} \\
\end{array}
\right),
\label{c2-129x}
\end{equation}
it can be shown 
that the matrix $T U T^{-1}$ has the form 
\begin{equation}
T U T^{-1}=\left(
\begin{array}{cc}
V_1      &     0    \\
0        &     V_2  \\
\end{array}
\right),
\label{c2-129y}
\end{equation}
where $U$ is the matrix in Eq. (\ref{c2-23}) used to represent the planar
fourcomplex number $u$. In Eq. (\ref{c2-129y}), $V_1, V_2$ are the matrices
\begin{equation}
V_1=\left(
\begin{array}{cc}
x+\frac{y-t}{\sqrt{2}}    &  z+\frac{y+t}{\sqrt{2}}   \\
-z-\frac{y+t}{\sqrt{2}}     &  x+\frac{y-t}{\sqrt{2}}   \\
\end{array}\right),\;\;
V_2=\left(
\begin{array}{cc}
x-\frac{y-t}{\sqrt{2}}    &-z+\frac{y+t}{\sqrt{2}}   \\
z-\frac{y+t}{\sqrt{2}}  &   x-\frac{y-t}{\sqrt{2}}   \\
\end{array}\right).
\label{c2-130x}
\end{equation}
In Eq. (\ref{c2-129y}), the symbols 0 denote the matrix
\begin{equation}
\left(
\begin{array}{cc}
0   &  0   \\
0   &  0   \\
\end{array}\right).
\label{c2-131x}
\end{equation}
The relations between the variables 
$x+(y-t)/\sqrt{2}, z+(y+t)/\sqrt{2}$,
$x-(y-t)/\sqrt{2},-z+(y+t)/\sqrt{2}$ for the multiplication
of planar fourcomplex numbers have been written in Eqs.
(\ref{c2-17})-(\ref{c2-20}). The matrix $T U T^{-1}$ provides an irreducible
representation
\cite{4} of the planar fourcomplex number $u$ in terms of matrices with real
coefficients.

\section{Polar omplex Numbers in Four Dimensions}

\subsection{Operations with polar fourcomplex numbers}

A polar fourcomplex number is determined by its four components $(x,y,z,t)$.
The sum of the polar fourcomplex numbers $(x,y,z,t)$ and
$(x^\prime,y^\prime,z^\prime,t^\prime)$ is the polar fourcomplex number
$(x+x^\prime,y+y^\prime,z+z^\prime,t+t^\prime)$.  The product of the polar
fourcomplex numbers $(x,y,z,t)$ and $(x^\prime,y^\prime,z^\prime,t^\prime)$ is
defined in this work to be the polar fourcomplex number
$(xx^\prime+yt^\prime+zz^\prime+ty^\prime,
xy^\prime+yx^\prime+zt^\prime+tz^\prime,
xz^\prime+yy^\prime+zx^\prime+tt^\prime,
xt^\prime+yz^\prime+zy^\prime+tx^\prime)$.  Polar fourcomplex numbers and their
operations can be represented by  writing the polar fourcomplex number
$(x,y,z,t)$ as $u=x+\alpha y+\beta z+\gamma t$, where $\alpha, \beta$ and
$\gamma$ are bases for which the multiplication rules are
\begin{equation}
\alpha^2=\beta, \:\beta^2=1, \:\gamma^2=\beta,
\alpha\beta=\beta\alpha=\gamma,\: 
\alpha\gamma=\gamma\alpha=-1, \:\beta\gamma=\gamma\beta=\alpha .
\label{ch1}
\end{equation}
Two polar fourcomplex numbers $u=x+\alpha y+\beta z+\gamma t, 
u^\prime=x^\prime+\alpha y^\prime+\beta z^\prime+\gamma t^\prime$ are equal, 
$u=u^\prime$, if and only if $x=x^\prime, y=y^\prime,
z=z^\prime, t=t^\prime$. 
If 
$u=x+\alpha y+\beta z+\gamma t, 
u^\prime=x^\prime+\alpha y^\prime+\beta z^\prime+\gamma t^\prime$
are polar fourcomplex numbers, 
the sum $u+u^\prime$ and the 
product $uu^\prime$ defined above can be obtained by applying the usual
algebraic rules to the sum 
$(x+\alpha y+\beta z+\gamma t)+ 
(x^\prime+\alpha y^\prime+\beta z^\prime+\gamma t^\prime)$
and to the product 
$(x+\alpha y+\beta z+\gamma t)
(x^\prime+\alpha y^\prime+\beta z^\prime+\gamma t^\prime)$,
and grouping of the resulting terms,
\begin{equation}
u+u^\prime=x+x^\prime+\alpha(y+y^\prime)+\beta(z+z^\prime)+\gamma(t+t^\prime),
\label{ch1a}
\end{equation}
\begin{eqnarray}
\lefteqn{uu^\prime=
xx^\prime+yt^\prime+zz^\prime+ty^\prime+
\alpha(xy^\prime+yx^\prime+zt^\prime+tz^\prime)+
\beta(xz^\prime+yy^\prime+zx^\prime+tt^\prime)\nonumber}\\
&&+\gamma(xt^\prime+yz^\prime+zy^\prime+tx^\prime).
\label{ch1b}
\end{eqnarray}

If $u,u^\prime,u^{\prime\prime}$ are polar fourcomplex numbers, the
multiplication is associative
\begin{equation}
(uu^\prime)u^{\prime\prime}=u(u^\prime u^{\prime\prime})
\label{ch2}
\end{equation}
and commutative
\begin{equation}
u u^\prime=u^\prime u ,
\label{ch3}
\end{equation}
as can be checked through direct calculation.
The polar fourcomplex zero is $0+\alpha\cdot 0+\beta\cdot 0+\gamma\cdot 0,$
denoted simply 0, and the polar fourcomplex unity is $1+\alpha\cdot
0+\beta\cdot 0+\gamma\cdot 0,$ denoted simply 1.

The inverse of the polar fourcomplex number 
$u=x+\alpha y+\beta z+\gamma t$ is a polar fourcomplex number
$u^\prime=x^\prime+\alpha y^\prime+\beta z^\prime+\gamma t^\prime$
having the property that
\begin{equation}
uu^\prime=1 .
\label{ch4}
\end{equation}
Written on components, the condition, Eq. (\ref{ch4}), is
\begin{equation}
\begin{array}{c}
xx^\prime+ty^\prime+zz^\prime+yt^\prime=1,\\
yx^\prime+xy^\prime+tz^\prime+zt^\prime=0,\\
zx^\prime+yy^\prime+xz^\prime+tt^\prime=0,\\
tx^\prime+zy^\prime+yz^\prime+xt^\prime=0.
\end{array}
\label{ch5}
\end{equation}
The system (\ref{ch5}) has the solution
\begin{equation}
x^\prime=\frac{x(x^2-z^2)+z(y^2+t^2)-2xyt}
{\nu} ,
\label{ch6a}
\end{equation}
\begin{equation}
y^\prime=\frac{-y(x^2+z^2)+t(y^2-t^2)+2xzt}
{\nu} ,
\label{ch6c}
\end{equation}
\begin{equation}
z^\prime=
\frac{-z(x^2-z^2)+x(y^2+t^2)-2yzt}
{\nu} ,
\label{ch6b}
\end{equation}
\begin{equation}
t^\prime=\frac{-t(x^2+z^2)-y(y^2-t^2)+2xyz}
{\nu} ,
\label{ch6d}
\end{equation}
provided that $\nu\not=0, $ where
\begin{equation}
\nu=x^4+z^4-y^4-t^4-2x^2z^2+2y^2t^2-4x^2yt-4yz^2t+4xy^2z+4xzt^2 .
\label{ch6e}
\end{equation}
The quantity $\nu$ can be written as
\begin{equation}
\nu=v_+ v_-\mu_+^2 ,
\label{ch7}
\end{equation}
where
\begin{equation}
v_+=x+y+z+t,\:v_-=x-y+z-t,
\label{ch8a}
\end{equation}
and
\begin{equation}
\mu_+^2=(x-z)^2+(y-t)^2 .
\label{ch8b}
\end{equation}
Then a polar fourcomplex number $q=x+\alpha y+\beta z+\gamma t$ has an inverse,
unless 
\begin{equation}
v_+=0 ,\:\:{\rm or}\:\: v_-=0,
\:\:{\rm or}\:\:\mu_+=0 . 
\label{ch9}
\end{equation}
The condition $v_+=0$ represents the 3-dimensional hyperplane $x+y+z+t=0$, the
condition $v_-=0$ represents the 3-dimensional hyperplane
$x-y+z-t=0$, and the condition $\mu_+=0$ represents the 2-dimensional
hyperplane $x=z, y=t$.
For arbitrary values of the variables $x,y,z,t$, the quantity $\nu$ can be
positive or negative. If $\nu\geq 0$, the quantity $\rho=\nu^{1/4}$
will be called amplitude of the polar fourcomplex number
$x+\alpha y+\beta z +\gamma t$.
Because of conditions (\ref{ch9}), these hyperplanes
will be called nodal hyperplanes. 

It can be shown that if $uu^\prime=0$ then either $u=0$, or $u^\prime=0$, or
the polar fourcomplex numbers $u, u^\prime$ belong to different members of the
pairs of orthogonal hypersurfaces listed further,
\begin{equation}
x+y+z+t=0 \:\:{\rm and} \:\: x^\prime=y^\prime=z^\prime=t^\prime,
\label{ch10a}
\end{equation}
\begin{equation}
x-y+z-t=0 \:\:{\rm and}\:\: x^\prime=-y^\prime=z^\prime=-t^\prime.
\label{ch10b}
\end{equation}
Divisors of zero also exist if the polar fourcomplex numbers $u,u^\prime$
belong to different members of the pair of two-dimensional hypersurfaces,
\begin{equation}
x-z=0,\:y-t=0\:\: {\rm and}\:\:  x^\prime+z^\prime=0 , \: y^\prime+t^\prime=0.
\label{ch11}
\end{equation}

\subsection{Geometric representation of polar fourcomplex numbers}

The polar fourcomplex number $x+\alpha y+\beta z+\gamma t$ can be represented
by the point $A$ of coordinates $(x,y,z,t)$.  If $O$ is the origin of the
four-dimensional space $x,y,z,t,$ the distance from $A$ to the origin $O$ can
be taken as
\begin{equation}
d^2=x^2+y^2+z^2+t^2 .
\label{ch12}
\end{equation}
The distance $d$ will be called modulus of the polar fourcomplex number
$x+\alpha y+\beta z +\gamma t$, $d=|u|$.

If $u=x+\alpha y+\beta z +\gamma t, u_1=x_1+\alpha y_1+\beta z_1 +\gamma t_1,
u_2=x_2+\alpha y_2+\beta z_2 +\gamma t_2$, and $u=u_1u_2$, and if
\begin{equation}
s_{j+}=x_j+y_j+z_j+t_j, \: 
s_{j-}=x_j-y_j+z_j-t_j,  
\label{ch13}
\end{equation}
for $j=1,2$, it can be shown that
\begin{equation}
v_+=s_{1+}s_{2+} ,\:\:
v_-=s_{1-}s_{2-}. \:\:
\label{ch14}
\end{equation}
The relations (\ref{ch14}) are a consequence of the identities
\begin{eqnarray}
\lefteqn{(x_1x_2+z_1z_2+t_1y_2+y_1t_2)+(x_1y_2+y_1x_2+z_1t_2+t_1z_2)
\nonumber}\\
&&+(x_1z_2+z_1x_2+y_1y_2+t_1t_2)+(x_1t_2+t_1x_2+z_1y_2+y_1z_2)\nonumber\\
&&=(x_1+y_1+z_1+t_1)(x_2+y_2+z_2+t_2),
\label{ch15}
\end{eqnarray}
\begin{eqnarray}
\lefteqn{(x_1x_2+z_1z_2+t_1y_2+y_1t_2)-(x_1y_2+y_1x_2+z_1t_2+t_1z_2)
\nonumber}\\
&&+(x_1z_2+z_1x_2+y_1y_2+t_1t_2)-(x_1t_2+t_1x_2+z_1y_2+y_1z_2)\nonumber\\
&&=(x_1+z_1-y_1-t_1)(x_2+z_2-y_2-t_2).
\label{ch16}
\end{eqnarray}

The differences 
\begin{equation}
v_1=x-z,\: \tilde v_1=y-t
\label{ch16a}
\end{equation}
can be written with the aid of the radius $\mu_+$, Eq. (\ref{ch8b}), and of the
azimuthal angle $\phi$, where $0\leq\phi<2\pi$ as
\begin{equation}
v_1=\mu_+\cos\phi,\:\:\tilde v_1=\mu_+\sin\phi .
\label{ch17}
\end{equation}
The variables $v_+, v_-, v_1, \tilde v_1$ will be called canonical 
polar fourcomplex variables.
The distance $d$, Eq. (\ref{ch12}), can then be written as
\begin{equation}
d^2
=\frac{1}{4}v_+^2+\frac{1}{4}v_-^2
+\frac{1}{2}\mu_+^2.  
\label{ch24d}
\end{equation}
It can be shown that if $u_1=x_1+\alpha y_1+\beta z_1+\gamma t_1, 
u_2=x_2+\alpha y_2+\beta z_2+\gamma t_2$ are polar fourcomplex
numbers of polar radii and angles $\rho_{1-}, \phi_1$ and
respectively $\rho_{2-}, \phi_2$, then the polar radius $\rho$ and
the angle $\phi$ of the product polar fourcomplex number $u_1u_2$
are 
\begin{equation}
\mu_+=\rho_{1-}\rho_{2-}, 
\label{ch17a}
\end{equation}
\begin{equation}
\phi=\phi_1+\phi_2.
\label{ch17b}
\end{equation}
The relation (\ref{ch17a}) is a consequence of the identity
\begin{eqnarray}
\lefteqn{\left[(x_1x_2+z_1z_2+y_1t_2+t_1y_2)
-(x_1z_2+z_1x_2+y_1y_2+t_1t_2)\right]^2\nonumber}\\
&&+\left[(x_1y_2+y_1x_2+z_1t_2+t_1z_2)
-(x_1t_2+t_1x_2+z_1y_2+y_1z_2)\right]^2\nonumber\\
&&=\left[(x_1-z_1)^2+(y_1-t_1)^2\right]\left[(x_2-z_2)^2+(y_2-t_2)^2\right],
\label{ch18}
\end{eqnarray}
and the relation (\ref{ch17b}) is a consequence of the identities
\begin{eqnarray}
\lefteqn{(x_1x_2+z_1z_2+y_1t_2+t_1y_2)
-(x_1z_2+z_1x_2+y_1y_2+t_1t_2)\nonumber}\\
&&=(x_1-z_1)(x_2-z_2)-(y_1-t_1)(y_2-t_2) ,
\label{ch19a}
\end{eqnarray}
\begin{eqnarray}
\lefteqn{(x_1y_2+y_1x_2+z_1t_2+t_1z_2)
-(x_1t_2+t_1x_2+z_1y_2+y_1z_2)\nonumber}\\
&&=(y_1-t_1)(x_2-z_2)+(x_1-z_1)(y_2-t_2) .
\label{ch19b}
\end{eqnarray}

A consequence of Eqs. (\ref{ch14}) and (\ref{ch17a}) is that if $u=u_1u_2$,
and $\nu_j=s_j s_j^{\prime\prime} \rho_{j-}$ , where $j=1,2$, then
\begin{equation}
\nu=\nu_1\nu_2 .
\label{ch20}
\end{equation}

The angles $\theta_+, \theta_-$ between the line $OA$ and the $v_+$ and
respectively $v_-$ axes are
\begin{equation}
\tan\theta_+=\frac{\sqrt{2}\mu_+}{v_+},
\tan\theta_-=\frac{\sqrt{2}\mu_+}{v_-},
\label{ch20b}
\end{equation}
where $0\leq\theta_+\leq\pi,\;0\leq\theta_-\leq\pi $.
The variable $\mu_+$ can be expressed with the aid of Eq. (\ref{ch24d}) as
\begin{equation}
\mu_+^2=2d^2\left(1+\frac{1}{\tan^2\theta_+}
+\frac{1}{\tan^2\theta_-}\right)^{-1}.
\label{ch20c}
\end{equation}
The coordinates $x,y,z,t$ can then be expressed in terms of the distance
$d$, of the polar angles $\theta_+,\theta_-$ and of the azimuthal angle $\phi$
as 
\begin{equation}
x=\frac{\mu_+(\tan\theta_++\tan\theta_-)}
{2\sqrt{2}\tan\theta_+\tan\theta_-}+\frac{1}{2}\mu_+\cos\phi,
\label{ch20i}
\end{equation}
\begin{equation}
y=\frac{\mu_+(-\tan\theta_++\tan\theta_-)}
{2\sqrt{2}\tan\theta_+\tan\theta_-}+\frac{1}{2}\mu_+\sin\phi,
\label{ch20k}
\end{equation}
\begin{equation}
z=\frac{\mu_+(\tan\theta_++\tan\theta_-)}
{2\sqrt{2}\tan\theta_+\tan\theta_-}-\frac{1}{2}\mu_+\cos\phi,
\label{ch20j}
\end{equation}
\begin{equation}
t=\frac{\mu_+(-\tan\theta_++\tan\theta_-)}
{2\sqrt{2}\tan\theta_+\tan\theta_-}-\frac{1}{2}\mu_+\sin\phi.
\label{ch20l}
\end{equation}
If $u=u_1u_2$, then Eqs. (\ref{ch14}) and (\ref{ch17a}) imply that
\begin{equation}
\tan\theta_+=\frac{1}{\sqrt{2}}\tan\theta_{1+}\tan\theta_{2+},\;
\tan\theta_-=\frac{1}{\sqrt{2}}\tan\theta_{1-}\tan\theta_{2-},
\label{ch20d}
\end{equation}
where 
\begin{equation}
\tan\theta_{j+}=\frac{\sqrt{2}\rho_{i+}}{s_j},
\tan\theta_{j-}=\frac{\sqrt{2}\rho_{i-}}{s_j^{\prime\prime}}.
\label{ch20e}
\end{equation}

An alternative choice of the angular variables is 
\begin{equation}
\mu_+=\sqrt{2}d\cos\theta,\;
v_+=2d\sin\theta\cos\lambda,\;v_-=2d\sin\theta\sin\lambda ,
\label{ch20f}
\end{equation}
where $0\leq\theta\leq\pi/2, 0\leq\lambda<2\pi$.
If $u=u_1u_2$, then
\begin{equation}
\tan\lambda=\tan\lambda_1\tan\lambda_2,\;
d\cos\theta=\sqrt{2}d_1d_2\cos\theta_1\cos\theta_2,
\label{ch20g}
\end{equation}
where 
\begin{equation}
\rho_{j-}=\sqrt{2}d_j\cos\theta_j,\;
s_j=2d_j\sin\theta_j\cos\lambda_j,\;
s_j^{\prime\prime}=2d_j\sin\theta_j\sin\lambda_j, 
\label{ch20h}
\end{equation}
for $j=1,2$. The coordinates $x,y,z,t$ can then be expressed in terms of the
distance $d$, of the polar angles $\theta_,\lambda$ and of the azimuthal angle
$\phi$ as
\begin{equation}
x=\frac{d}{\sqrt{2}}\sin\theta
\sin(\lambda+\pi/4)+\frac{d}{\sqrt{2}}\cos\theta\cos\phi,
\label{ch20m}
\end{equation}
\begin{equation}
y=\frac{d}{\sqrt{2}}\sin\theta
\cos(\lambda+\pi/4)+\frac{d}{\sqrt{2}}\cos\theta\sin\phi,
\label{ch20o}
\end{equation}
\begin{equation}
z=\frac{d}{\sqrt{2}}\sin\theta\sin(\lambda+\pi/4)
-\frac{d}{\sqrt{2}}\cos\theta\cos\phi,
\label{ch20n}
\end{equation}
\begin{equation}
t=\frac{d}{\sqrt{2}}\sin\theta\cos(\lambda+\pi/4)
-\frac{d}{\sqrt{2}}\cos\theta\sin\phi.
\label{ch20p}
\end{equation}

The polar fourcomplex numbers
\begin{equation}
e_+=\frac{1+\alpha+\beta+\gamma}{4},\:
e_-=\frac{1-\alpha+\beta-\gamma}{4},
\label{ch21}
\end{equation}
have the property that
\begin{equation}
e_+^2=e_+, \: 
e_-^2=e_-, \:e_+e_-=0.
\label{ch22}
\end{equation}
The polar fourcomplex numbers 
\begin{equation}
e_1=\frac{1-\beta}{2}, \tilde e_1=\frac{\alpha-\gamma}{2}
\label{ch22b}
\end{equation}
have the property that
\begin{equation}
e_1^2=e_1,\:\:
\tilde e_1^2=-e_1,\:\:
e_1\tilde e_1=\tilde e_1 .
\label{ch23}
\end{equation}
The polar fourcomplex numbers $ e_+, e_-$ are orthogonal to
$e_1, \tilde e_1$, 
\begin{equation}
e_+\:e_1=0,\:\:e_+\:\tilde e_1=0,\:\:
e_-\:e_1=0,
\:\:e_-\:\tilde e_1=0.\:\: 
\label{ch23b}
\end{equation}
The polar fourcomplex number $q=x+\alpha y+\beta z+\gamma t$ can then be
written as 
\begin{equation}
x+\alpha y+\beta z+\gamma t
=v_+e_++v_-e_-+v_1e_1
+\tilde v_1\tilde e_1.
\label{ch24}
\end{equation}
The ensemble $e_+, e_-, e_1, \tilde e_1$ will be called the canonical
polar fourcomplex base, and Eq. (\ref{ch24}) gives the canonical form of the
polar fourcomplex number.
Thus, the product of the polar fourcomplex numbers $u, u^\prime$ can be
expressed as 
\begin{equation}
uu^\prime=v_+v_+^\prime e_++v_-v_-^\prime e_-
+(v_1v_1^\prime-\tilde v_1\tilde v_1^\prime)e_1+
(v_1\tilde v_1^\prime+v_1^\prime\tilde v_1)\tilde e_1,
\label{ch24b}
\end{equation}
where $v_+^\prime=x^\prime+y^\prime+z^\prime+t^\prime, 
v_-^\prime=x^\prime-y^\prime+z^\prime-t^\prime,
v_1^\prime=x^\prime-y^\prime, \tilde v_1^\prime=z^\prime-t^\prime$.
The moduli of the bases used in Eq. (\ref{ch24}) are
\begin{equation}
|e_+|=\frac{1}{2},\;|e_-|=\frac{1}{2},\;
\left|e_1\right|=\frac{1}{\sqrt{2}},\;
\left|\tilde e_1\right|=\frac{1}{\sqrt{2}}.
\label{ch24c}
\end{equation}

The fact that the amplitude of the product is equal to the product of the
amplitudes, as written in Eq. (\ref{ch20}), can be demonstrated also by using a
representation of the multiplication of the polar fourcomplex numbers by
matrices, in which the polar fourcomplex number $u=x+\alpha y+\beta z+\gamma t$
is represented by the matrix
\begin{equation}
A=\left(\begin{array}{cccc}
x &y &z &t\\
t&x &y &z\\
z&t&x &y\\
y&z&t&x 
\end{array}\right) .
\label{ch25a}
\end{equation}
The product $u=x+\alpha y+\beta z+\gamma t$ of the polar fourcomplex numbers
$u_1=x_1+\alpha y_1+\beta z_1+\gamma t_1, u_2=x_2+\alpha y_2+\beta z_2+\gamma
t_2$, can be represented by the matrix multiplication 
\begin{equation}
A=A_1A_2.
\label{ch25b}
\end{equation}
It can be checked that the determinant ${\rm det}(A)$ of the matrix $A$ is
\begin{equation}
{\rm det}A = \nu .
\label{ch25c}
\end{equation}
The identity (\ref{ch20}) is then a consequence of the fact the determinant 
of the product of matrices is equal to the product of the determinants 
of the factor matrices.

\subsection{The polar fourdimensional cosexponential functions}

The exponential function of a fourcomplex variable $u$ and the addition
theorem for the exponential function have been written in Eqs. (\ref{g1}) and
(\ref{g2}).
If $u=x+\alpha y+\beta z+\gamma t$, then $\exp u$ can be calculated as $\exp
u=\exp x \cdot \exp (\alpha y) \cdot \exp (\beta z)\cdot \exp (\gamma t)$.
According to Eq. (\ref{ch1}),
\begin{eqnarray}
\begin{array}{l}
\alpha^{4m}=1, \alpha^{4m+1}=\alpha, 
\alpha^{4m+2}=\beta, \alpha^{4m+3}=\gamma,  \\
\beta^{2m}=1, \beta^{2m+1}=\beta, \\
\gamma^{4m}=1, \gamma^{4m+1}=\gamma, 
\gamma^{4m+2}=\beta, \gamma^{4m+3}=\alpha, \\
\end{array}
\label{ch28}
\end{eqnarray}
where $m$ is a natural number,
so that $\exp (\alpha y), \: \exp(\beta z)$ and $\exp(\gamma z)$ can be written
as 
\begin{equation}
\exp (\beta z) = \cosh z +\beta \sinh z ,
\label{ch29}
\end{equation}
and
\begin{equation}
\exp (\alpha y) = g_{40}(y)+\alpha g_{41}(y)
+\beta g_{42}(y) +\gamma g_{43}(y) ,  
\label{ch30a}
\end{equation}
\begin{equation}
\exp (\gamma t) = g_{40}(t)+\gamma g_{41}(t)
+\beta g_{42}(t) +\alpha g_{43}(t) ,
\label{ch30b}
\end{equation}
where the four-dimensional cosexponential 
functions $g_{40}, g_{41}, g_{42}, g_{43}$ are
defined by the series 
\begin{equation}
g_{40}(x)=1+x^4/4!+x^8/8!+\cdots ,
\label{ch30c}
\end{equation}
\begin{equation}
g_{41}(x)=x+x^5/5!+x^9/9!+\cdots,
\label{ch30d}
\end{equation}
\begin{equation}
g_{42}(x)=x^2/2!+x^6/6!+x^{10}/10!+\cdots ,
\label{ch30e}
\end{equation}
\begin{equation}
g_{43}(x)=x^3/3!+x^7/7!+x^{11}/11!+\cdots .
\label{ch30f}
\end{equation}
The functions $g_{40}, g_{42}$ are even, the functions $g_{41}, g_{43}$ are
odd, 
\begin{equation}
g_{40}(-u)=g_{40}(u),\:g_{42}(-u)=g_{42}(u),
\:g_{41}(-u)=-g_{41}(u),\:g_{43}(-u)=-g_{43}(u).
\label{ch30feo}
\end{equation}
It can be seen from Eqs. (\ref{ch30c})-(\ref{ch30f}) that
\begin{equation}
g_{40}+g_{41}+g_{42}+g_{43}=e^x , \:g_{40}-g_{41}+g_{42}-g_{43}=e^{-x} ,
\label{ch30fep}
\end{equation}
and
\begin{equation}
g_{40}-g_{42}=\cos x , \: g_{41}-g_{43}=\sin x ,
\label{ch30feq}
\end{equation}
so that
\begin{equation}
(g_{40}+g_{41}+g_{42}+g_{43})(g_{40}-g_{41}+g_{42}-g_{43})
\left[(g_{40}-g_{42})^2+(g_{41}-g_{43})^2\right]=1 ,
\label{ch30fer}
\end{equation}
which can be also written as
\begin{equation}
g_{40}^4-g_{41}^4+g_{42}^4-g_{43}^4-2(g_{40}^2g_{42}^2-g_{41}^2g_{43}^2)
-4(g_{40}^2g_{41}g_{43}+g_{42}^2g_{41}g_{43}
-g_{41}^2g_{40}g_{42}-g_{43}^2g_{40}g_{42})=1. 
\label{ch30fes}
\end{equation}
The combination of terms in Eq. (\ref{ch30fes}) is similar to that in Eq.
(\ref{ch6e}). 

Addition theorems for the four-dimensional cosexponential functions can be
obtained from the relation $\exp \alpha(x+y)=\exp\alpha x\cdot\exp\alpha y $,
by substituting the expression of the exponentials as given in Eq.
(\ref{ch30a}),
\begin{equation}
g_{40}(x+y)=g_{40}(x)g_{40}(y)
+g_{41}(x)g_{43}(y)+g_{42}(x)g_{42}(y)+g_{43}(x)g_{41}(y) ,
\label{ch30g}
\end{equation}
\begin{equation}
g_{41}(x+y)=g_{40}(x)g_{41}(y)+g_{41}(x)g_{40}(y)
+g_{42}(x)g_{43}(y)+g_{43}(x)g_{42}(y) ,
\label{ch30h}
\end{equation}
\begin{equation}
g_{42}(x+y)=g_{40}(x)g_{42}(y)+g_{41}(x)g_{41}(y)
+g_{42}(x)g_{40}(y)+g_{43}(x)g_{43}(y) ,
\label{ch30i}
\end{equation}
\begin{equation}
g_{43}(x+y)=g_{40}(x)g_{43}(y)+g_{41}(x)g_{42}(y)
+g_{42}(x)g_{41}(y)+g_{43}(x)g_{40}(y) .
\label{ch30j}
\end{equation}
For $x=y$ the relations (\ref{ch30g})-(\ref{ch30j}) take the form
\begin{equation}
g_{40}(2x)=g_{40}^2(x)+g_{42}^2(x)+2g_{41}(x)g_{43}(x) ,
\label{ch30gg}
\end{equation}
\begin{equation}
g_{41}(2x)=2g_{40}(x)g_{41}(x)+2g_{42}(x)g_{43}(x) ,
\label{ch30hh}
\end{equation}
\begin{equation}
g_{42}(2x)=g_{41}^2(x)+g_{43}^2(x)+2g_{40}(x)g_{42}(x) ,
\label{ch30ii}
\end{equation}
\begin{equation}
g_{43}(2x)=2g_{40}(x)g_{43}(x)+2g_{41}(x)g_{42}(x) .
\label{ch30jj}
\end{equation}
For $x=-y$ the relations (\ref{ch30g})-(\ref{ch30j}) and (\ref{ch30feo}) yield
\begin{equation}
g_{40}^2(x)+g_{42}^2(x)-2g_{41}(x)g_{43}(x)=1 ,
\label{ch30eog}
\end{equation}
\begin{equation}
g_{41}^2(x)+g_{43}^2(x)-2g_{40}(x)g_{42}(x)=0 .
\label{ch30eoi}
\end{equation}
From Eqs. (\ref{ch29})-(\ref{ch30b}) it can be shown that, for $m$ integer,
\begin{equation}
(\cosh z +\beta \sinh z)^m=\cosh mz +\beta \sinh mz ,
\label{ch2929}
\end{equation}
and
\begin{equation}
[g_{40}(y)+\alpha g_{41}(y)+\beta g_{42}(y) +\gamma g_{43}(y)]^m= 
g_{40}(my)+\alpha g_{41}(my)+\beta g_{42}(my) +\gamma g_{43}(my),  
\label{ch30a30a}
\end{equation}
\begin{equation}
[g_{40}(t)+\gamma g_{41}(t)+\beta g_{42}(t) +\alpha g_{43}(t)]^m=
g_{40}(mt)+\gamma g_{41}(mt)+\beta g_{42}(mt) +\alpha g_{43}(mt) .
\label{ch30b30b}
\end{equation}
Since
\begin{equation}
(\alpha+\gamma)^{2m}=2^{2m-1}(1+\beta),
\: (\alpha+\gamma)^{2m+1}=2^{2m}(\alpha+\gamma) ,
\label{ch30k}
\end{equation}
it can be shown from the definition of the exponential function, Eq.
(\ref{g1}) that
\begin{equation}
\exp(\alpha+\gamma)x=e_1+\frac{1+\beta}{2}\cosh 2x+\frac{\alpha+\gamma}{2}
\sinh 2x .
\label{ch30l}
\end{equation}
Substituting in the relation $\exp(\alpha+\gamma)x=\exp\alpha x\exp\gamma x$
the expression of the exponentials from Eqs. (\ref{ch30a}), (\ref{ch30b}) and
(\ref{ch30l}) yields
\begin{equation}
g_{40}^2+g_{41}^2+g_{42}^2+g_{43}^2=\frac{1+\cosh 2x}{2} ,
\label{ch30m}
\end{equation}
\begin{equation}
g_{40} g_{42} +g_{41} g_{43}=\frac{-1+\cosh 2x}{4} ,
\label{ch30mm}
\end{equation}
\begin{equation}
g_{40}g_{41}+g_{40}g_{43}+g_{41}g_{42}+g_{42}g_{43}=\frac{1}{2}\sinh 2x ,
\label{ch30n}
\end{equation}
where $g_{40}, g_{41}, g_{42}, g_{43}$ are functions of $x$.

Similarly, since
\begin{equation}
(\alpha-\gamma)^{2m}=(-1)^m 2^{2m-1}(1-\beta),\:
(\alpha-\gamma)^{2m+1}=(-1)^m2^{2m}(\alpha-\gamma) , 
\label{ch30k30k}
\end{equation}
it can be shown from the definition of the exponential function, Eq. (\ref{g1})
that
\begin{equation}
\exp(\alpha-\gamma)x=\frac{1+\beta}{2}+e_1\cos 2x+\tilde e_1
\sin 2x .
\label{ch30l30l}
\end{equation}
Substituting in the relation $\exp(\alpha-\gamma)x=\exp\alpha x\exp(-\gamma x)$
the expression of the exponentials from Eqs. (\ref{ch30a}), (\ref{ch30b}) and
(\ref{ch30l30l}) yields
\begin{equation}
g_{40}^2-g_{41}^2+g_{42}^2-g_{43}^2=\frac{1+\cos 2x}{2} ,
\label{ch30m30m}
\end{equation}
\begin{equation}
g_{40} g_{42} -g_{41} g_{43}=\frac{1-\cos 2x}{2} ,
\label{ch30mm30mm}
\end{equation}
\begin{equation}
g_{40}g_{41}-g_{40}g_{43}-g_{41}g_{42}+g_{42}g_{43}=\frac{1}{2}\sin 2x ,
\label{ch30n30n}
\end{equation}
where $g_{40}, g_{41}, g_{42}, g_{43}$ are functions of $x$.

The expressions of the four-dimensional cosexponential
functions are
\begin{equation}
g_{40}(x)=\frac{1}{2}(\cosh x + \cos x) ,
\label{ch30c30c}
\end{equation}
\begin{equation}
g_{41}(x)=\frac{1}{2}(\sinh x+\sin x),
\label{ch30d30d}
\end{equation}
\begin{equation}
g_{42}(x)=\frac{1}{2}(\cosh x - \cos x) ,
\label{ch30e30e}
\end{equation}
\begin{equation}
g_{43}(x)=\frac{1}{2}(\sinh x-\sin x).
\label{ch30f30f}
\end{equation}
The graphs of these four-dimensional cosexponential functions are shown in 
Fig. 4. 

It can be checked that the cosexponential functions are solutions of the
fourth-order differential equation
\begin{equation}
\frac{{\rm d}^4\zeta}{{\rm d}u^4}=\zeta ,
\label{ch30x}
\end{equation}
whose solutions are of the form
$\zeta(u)=Ag_{40}(u)+Bg_{41}(u)+Cg_{42(u)}+Dg_{43}(u).$ It can also be checked
that the derivatives of the cosexponential functions are related by
\begin{equation}
\frac{dg_{40}}{dw}=g_{43}, \:
\frac{dg_{41}}{dw}=g_{40}, \:
\frac{dg_{42}}{dw}=g_{41} ,
\frac{dg_{43}}{dw}=g_{42} .
\label{ch30x30x}
\end{equation}

\subsection{The exponential and trigonometric forms 
of a polar fourcomplex number}

The polar fourcomplex numbers $u=x+\alpha y+\beta z+\gamma t$ for which
$v_+=x+y+z+t>0, \:  v_-=x-y+z-t>0$ can be written in the form 
\begin{equation}
x+\alpha y+\beta z+\gamma t=e^{x_1+\alpha y_1+\beta z_1+\gamma t_1} .
\label{ch31-34}
\end{equation}

The expressions of $x_1, y_1, z_1, t_1$ as functions of $x, y, z, t$ can be
obtained by developing $e^{\alpha y_1}, e^{\beta z_1}$ and $e^{\gamma t_1}$
with the aid of Eqs. (\ref{ch29})-(\ref{ch30b}), by multiplying these
expressions and separating the hypercomplex components, and then substituting
the expressions of the four-dimensional cosexponential functions, Eqs.
(\ref{ch30c30c})-(\ref{ch30f30f}),
\begin{equation}
x+y+z+t=e^{x_1+y_1+z_1+t_1},
\label{ch35}
\end{equation}
\begin{equation}
x-z=e^{x_1-z_1}\cos(y_1-t_1),
\label{ch36}
\end{equation}
\begin{equation}
x-y+z-t=e^{x_1-y_1+z_1-t_1},
\label{ch37}
\end{equation}
\begin{equation}
y-t=e^{x_1-z_1}\sin(y_1-t_1).
\label{ch38}
\end{equation}
It can be shown from Eqs. (\ref{ch35})-(\ref{ch38}) that
\begin{equation}
x_1=\frac{1}{2} \ln(\mu_+\mu_-) , \:
y_1=\frac{1}{2}(\phi+\omega) ,\:
z_1=\frac{1}{2}\ln\frac{\mu_-}{\mu_+},\:
t_1=\frac{1}{2}(\phi-\omega),
\label{ch39}
\end{equation}
where
\begin{equation}
\mu_-^2=(x+z)^2-(y+t)^2=v_+v_-,\: \:v_+>0, \:\:v_->0.
\label{ch39a}
\end{equation}
The quantitites $\phi$ and $\omega$ are determined by
\begin{equation}
\cos\phi=(x-z)/\mu_+,\:\sin\phi=(y-t)/\mu_+ 
\label{ch39b}
\end{equation}
and
\begin{equation}
\cosh\omega=(x+z)/\mu_-,\:
\sinh\omega=(y+t)/\mu_- .
\label{ch39c}
\end{equation}
The explicit form of $\omega$ is
\begin{equation}
\omega=\frac{1}{2}\ln\frac{x+y+z+t}{x-y+z-t}.
\label{ch39cc}
\end{equation}
If $u=u_1u_2$, and $\mu_{j-}^2=(x_j+z_j)^2-(y_j+t_j)^2, j=1,2$, it can be
checked with the aid of Eqs. (\ref{ch14}) that
\begin{equation}
\mu_-=\mu_{1-}\mu_{2-}. 
\label{ch39d}
\end{equation}
Moreover, if
$\cosh\omega_j=(x_j+z_j)/\mu_{j-},\:\sinh\omega_j=(y_j+t_j)/\mu_{j-}, j=1,2$,
it can be checked that 
\begin{equation}
\omega=\omega_1+\omega_2.
\label{ch39e}
\end{equation}
The relation (\ref{ch39e}) is a consequence of the identities
\begin{eqnarray}
\lefteqn{(x_1x_2+z_1z_2+y_1t_2+t_1y_2)
+(x_1z_2+z_1x_2+y_1y_2+t_1t_2)\nonumber}\\
&&=(x_1+z_1)(x_2+z_2)+(y_1+t_1)(y_2+t_2) ,
\label{ch19c}
\end{eqnarray}
\begin{eqnarray}
\lefteqn{(x_1y_2+y_1x_2+z_1t_2+t_1z_2)
+(x_1t_2+t_1x_2+z_1y_2+y_1z_2)\nonumber}\\
&&=(y_1+t_1)(x_2+z_2)+(x_1+z_1)(y_2+t_2) .
\label{ch19d}
\end{eqnarray}

According to Eq. (\ref{ch39b}), $\phi$ is a cyclic variable, $0\leq\phi<2\pi$.
As it has been assumed that $v_+>0, v_->0$, it follows that $x+z>0$
and $x+z>|y+t|$. The range of the variable $\omega$ is $-\infty<\omega<\infty$.
The exponential form of the polar fourcomplex number $u$ is then
\begin{equation}
u=\rho\exp\left[\frac{1}{2}\beta \ln\frac{\mu_-}{\mu_+}
+\frac{1}{2}\alpha (\omega+\phi) 
+\frac{1}{2}\gamma (\omega-\phi)\right] ,
\label{ch40}
\end{equation}
where
\begin{equation}
\rho=(\mu_+\mu_-)^{1/2}.
\label{ch40a}
\end{equation}
If $u=u_1u_2$, and $\rho_j=(\mu_{j+}\mu_{j-})^{1/2}, j=1,2$, then 
from Eqs. (\ref{ch17a}) and (\ref{ch39d}) it results that
\begin{equation}
\rho=\rho_1\rho_2.
\label{ch40bb}
\end{equation}

It can be checked with the aid of Eq. (\ref{ch23}) that
\begin{equation}
\exp\left(\tilde e_1\phi\right)
=\frac{1+\beta}{2}+e_1\cos\phi+\tilde e_1\sin\phi,
\label{ch40x}
\end{equation}
which shows that $e^{(\alpha-\gamma)\phi/2}$ is a periodic function of $\phi$,
with period $2\pi$.
The modulus has the property that
\begin{equation}
\left|u\exp\left(\tilde e_1\phi\right)\right|=|u|.
\label{ch40bc}
\end{equation}

By introducing in Eq. (\ref{ch40}) the polar angles $\theta_+, \theta_-$
defined in Eqs. (\ref{ch20b}), the exponential form of the fourcomplex number
$u$ becomes
\begin{equation}
u=\rho\exp\left[\frac{1}{4}(\alpha+\beta+\gamma) 
\ln\frac{\sqrt{2}}{\tan\theta_+}
-\frac{1}{4}(\alpha-\beta+\gamma) \ln\frac{\sqrt{2}}{\tan\theta_-}
+\tilde e_1\phi\right] ,
\label{ch40c}
\end{equation}
where $0<\theta_+<\pi/2, 0<\theta_-<\pi/2$.
The relation between the amplitude $\rho$, Eq. (\ref{ch40a}), and the
distance $d$, Eq. (\ref{ch12}), is, according to Eqs. (\ref{ch20b}) and
(\ref{ch20c}), 
\begin{equation}
\rho=\frac{2^{3/4}d}{\left(\tan\theta_+\tan\theta_-\right)^{1/4}}
\left(1+\frac{1}{\tan^2\theta_+}+\frac{1}{\tan^2\theta_-}\right)^{-1/2}.
\label{ch40d}
\end{equation}
Using the properties of the vectors $e_+, e_-$ written in Eq.
(\ref{ch22}), the first part of the exponential, Eq. (\ref{ch40c}) can be
developed as
\begin{eqnarray}
\lefteqn{\exp\left[\frac{1}{4}(\alpha+\beta+\gamma) 
\ln\frac{\sqrt{2}}{\tan\theta_+}
-\frac{1}{4}(\alpha-\beta+\gamma) 
\ln\frac{\sqrt{2}}{\tan\theta_-}\right] \nonumber}\\
&&=\left(\frac{1}{2}\tan\theta_+\tan\theta_-\right)^{1/4}
\left(e_1+e_+\frac{\sqrt{2}}{\tan\theta_+}
+e_-\frac{\sqrt{2}}{\tan\theta_-}\right).
\label{ch40e}
\end{eqnarray}
The fourcomplex number $u$, Eq. (\ref{ch40c}), can then be written as
\begin{equation}
u=d\sqrt{2}
\left(1+\frac{1}{\tan^2\theta_+}+\frac{1}{\tan^2\theta_-}\right)^{-1/2}
\left(e_1+e_+\frac{\sqrt{2}}{\tan\theta_+}
+e_-\frac{\sqrt{2}}{\tan\theta_-}\right)
\exp\left(\tilde e_1\phi\right),
\label{ch40f}
\end{equation}
which is the trigonometric form of the fourcomplex number $u$.

The polar angles $\theta_+, \theta_-$, Eq. (\ref{ch20b}), can be expressed in
terms of the variables $\theta, \lambda$, Eq. (\ref{ch20f}), as
\begin{equation}
\tan\theta_+=\frac{1}{\tan\theta\cos\lambda},\;
\tan\theta_-=\frac{1}{\tan\theta\sin\lambda},
\label{ch40g}
\end{equation}
so that
\begin{equation}
1+\frac{1}{\tan^2\theta_+}+\frac{1}{\tan^2\theta_-}=\frac{1}{\cos^2\theta}.
\label{ch40h}
\end{equation}
The exponential form of the fourcomplex number $u$, written in terms of the
amplitude $\rho$ and of the angles $\theta, \lambda, \phi$ is
\begin{equation}
u=\rho\exp\left[\frac{1}{4}(\alpha+\beta+\gamma)
\ln(\sqrt{2}\tan\theta\cos\lambda) 
-\frac{1}{4}(\alpha-\beta+\gamma) \ln(\sqrt{2}\tan\theta\sin\lambda)
+\tilde e_1\phi
\right],
\label{ch40i}
\end{equation}
where $0\leq\lambda<\pi/2$.
The trigonometric form of the fourcomplex number $u$, written in terms of the
amplitude $\rho$ and of the angles $\theta, \lambda, \phi$ is
\begin{equation}
u=d\sqrt{2}
\left(e_1\cos\theta+e_+\sqrt{2}\sin\theta\cos\lambda
+e_-\sqrt{2}\sin\theta\sin\lambda\right)
\exp\left(\tilde e_1\phi\right).
\label{ch40j}
\end{equation}

If $u=u_1u_2$, it can be shown with the aid of the trigonometric form, Eq.
(\ref{ch40f}), that the modulus of the product as a function of the polar
angles is 
\begin{eqnarray}
\lefteqn{d^2=4d_1^2d_2^2\left(\frac{1}{2}
+\frac{1}{\tan^2\theta_{1+}\tan^2\theta_{2+}}
+\frac{1}{\tan^2\theta_{1-}\tan^2\theta_{2-}}\right)\nonumber}\\
&&\left(1+\frac{1}{\tan^2\theta_{1+}}+\frac{1}{\tan^2\theta_{1-}}\right)^{-1}
\left(1+\frac{1}{\tan^2\theta_{2+}}+\frac{1}{\tan^2\theta_{2-}}\right)^{-1}.
\label{ch40k}
\end{eqnarray}
The modulus $d$ of the product $u_1u_2$ can be expressed alternatively in terms
of the angles $\theta, \lambda, \phi$ with the aid of the
trigonometric form, Eq. (\ref{ch40j}), as
\begin{eqnarray}
\lefteqn{d^2=4d_1^2d_2^2\left(\frac{1}{2}\cos^2\theta_1\cos^2\theta_2
+\sin^2\theta_1\sin^2\theta_2\cos^2\lambda_1\cos^2\lambda_2
+\sin^2\theta_1\sin^2\theta_2\sin^2\lambda_1\sin^2\lambda_2\right).\nonumber}\\
&&
\label{ch40l}
\end{eqnarray}

\subsection{Elementary functions of polar fourcomplex variables}

The logarithm $u_1$ of the polar fourcomplex number $u$, $u_1=\ln u$, can be
defined as the solution of the equation
\begin{equation}
u=e^{u_1} ,
\label{ch41}
\end{equation}
written explicitly previously in Eq. (\ref{ch31-34}), for $u_1$ as a function
of $u$. From Eq. (\ref{ch40}) it results that 
\begin{equation}
\ln u=\ln\rho+
\frac{1}{2}\beta \ln\frac{\mu_-}{\mu_+}
+\frac{1}{2}\alpha (\omega+\phi) 
+\frac{1}{2}\gamma (\omega-\phi).
\label{ch42}
\end{equation}
If the fourcomplex number $u$ is written in terms of the
amplitude $\rho$ and of the angles $\theta_+, \theta_-, \phi$, the logarithm
is 
\begin{equation}
\ln u=\ln \rho+
\frac{1}{4}(\alpha+\beta+\gamma) \ln\frac{\sqrt{2}}{\tan\theta_+}
-\frac{1}{4}(\alpha-\beta+\gamma) \ln\frac{\sqrt{2}}{\tan\theta_-}
+\tilde e_1\phi,
\label{ch42b}
\end{equation}
where $0<\theta_+<\pi/2, 0<\theta_+<\pi/2$.
If the fourcomplex number $u$ is written in terms of the
amplitude $\rho$ and of the angles $\theta, \lambda, \phi$, the logarithm is
\begin{equation}
\ln u=\ln \rho+
\frac{1}{4}(\alpha+\beta+\gamma)
\ln(\sqrt{2}\tan\theta\cos\lambda) 
-\frac{1}{4}(\alpha-\beta+\gamma) \ln(\sqrt{2}\tan\theta\sin\lambda)
+\tilde e_1\phi,
\label{ch42c}
\end{equation}
where $0<\theta<\pi/2, 0\leq\lambda<\pi/2$.
The logarithm is multivalued because of the term proportional to $\phi$.
It can be inferred from Eq. (\ref{ch42b}) that
\begin{equation}
\ln(u_1u_2)=\ln u_1+\ln u_2 ,
\label{ch43}
\end{equation}
up to multiples of $\pi(\alpha-\gamma)$.
If the expressions of $\rho, \mu_+, \mu_-$ and $\omega$ in terms of
$x,y,z,t$ are introduced in
Eq. (\ref{ch42}), the logarithm of the polar fourcomplex number becomes
\begin{equation}
\ln u=\frac{1+\alpha+\beta+\gamma}{4}\ln(x+y+z+t)
+\frac{1-\alpha+\beta-\gamma}{4}\ln(x-y+z-t)
+e_1\ln\mu_+
+\tilde e_1\phi .
\label{ch44}
\end{equation}

The power function $u^m$ can be defined for $v_+>0, 
v_->0$ and real values of $m$ as
\begin{equation}
u^m=e^{m\ln u} .
\label{ch44b}
\end{equation}
The power function is multivalued unless $m$ is an integer. 
For integer $m$, it can be inferred from Eqs. (\ref{ch42}) and (\ref{ch44b})
that
\begin{equation}
(u_1u_2)^m=u_1^m\:u_2^m .
\label{ch45}
\end{equation}

Using the expression (\ref{ch44}) for $\ln u$ and the relations
(\ref{ch22})-(\ref{ch23b}) it can be shown that 
\begin{eqnarray}
u^m=e_+v_+^m+e_-v_-^m
+\mu_+^m\left(e_1\cos m\phi
+\tilde e_1\sin m\phi\right).
\label{ch46}
\end{eqnarray}
For integer $m$, the relation (\ref{ch46}) is valid for any $x,y,z,t$. 
For natural $m$ this relation can be written as
\begin{eqnarray}
u^m=e_+v_+^m+e_-v_-^m
+\left[e_1(x-z)
+\tilde e_1(y-t)\right]^m,
\label{ch46b}
\end{eqnarray}
as can be shown with the aid of the relation
\begin{eqnarray}
e_1\cos m\phi+\tilde e_1\sin m\phi
=\left(e_1\cos \phi+\tilde e_1\sin \phi\right)^m,
\label{ch46c}
\end{eqnarray}
valid for natural $m$.
For $m=-1$ the relation (\ref{ch46}) becomes
\begin{eqnarray}
\lefteqn{\frac{1}{x+\alpha y+\beta z+\gamma t}=
\frac{1}{4}\left(\frac{1+\alpha+\beta+\gamma}{x+y+z+t}
+\frac{1-\alpha+\beta-\gamma}{x-y+z-t}\right)\nonumber}\\
&&+\frac{1}{2}\frac{(1-\beta)(x-z)-(\alpha-\gamma)(y-t)}{(x-z)^2+(y-t)^2}.
\label{ch47}
\end{eqnarray}
If $m=2$, it can be checked that the right-hand side of
Eq. (\ref{ch46b}) is equal to $(x+\alpha y+\beta z+\gamma t)^2=
x^2+z^2+2yt+2\alpha(xy+zt)+\beta(y^2+t^2+2xz)+2\gamma(xt+yz)$.

The trigonometric functions of the fourcomplex variable
$u$ and the addition theorems for these functions have been written in Eqs.
(\ref{g3})-(\ref{g6}). 
The cosine and sine functions of the hypercomplex variables $\alpha y, 
\beta z$ and $ \gamma t$ can be expressed as
\begin{equation}
\cos\alpha y=g_{40}-\beta g_{42}, \: \sin\alpha y=\alpha g_{41}-\gamma g_{43}, 
\label{ch52}
\end{equation}
\begin{equation}
\cos\beta y=\cos y, \: \sin\beta y=\beta\sin y, 
\label{ch51}
\end{equation}
\begin{equation}
\cos\gamma y=g_{40}-\beta g_{42}, \: \sin\gamma y=\gamma g_{41}-\alpha g_{43} .
\label{ch53}
\end{equation}
The cosine and sine functions of a polar fourcomplex number $x+\alpha y+\beta
z+\gamma t$ can then be
expressed in terms of elementary functions with the aid of the addition
theorems Eqs. (\ref{g5}), (\ref{g6}) and of the expressions in  Eqs. 
(\ref{ch52})-(\ref{ch53}). 

The hyperbolic functions of the fourcomplex variable
$u$ and the addition theorems for these functions have been written in Eqs.
(\ref{g7})-(\ref{g10}). 
The hyperbolic cosine and sine functions of the hypercomplex variables $\alpha
y,
\beta z$ and $ \gamma t$ can be expressed as
\begin{equation}
\cosh\alpha y=g_{40}+\beta g_{42}, \: 
\sinh\alpha y=\alpha g_{41}+\gamma g_{43}, 
\label{ch59}
\end{equation}
\begin{equation}
\cosh\beta y=\cosh y, \: \sinh\beta y=\beta\sinh y, 
\label{ch58}
\end{equation}
\begin{equation}
\cosh\gamma y=g_{40}+\beta g_{42}, 
\: \sinh\gamma y=\gamma g_{41}+\alpha g_{43} .
\label{ch60-63}
\end{equation}
The hyperbolic cosine and sine
functions of a polar fourcomplex number $x+\alpha y+\beta
z+\gamma t$ can then be
expressed in terms of elementary functions with the aid of the addition
theorems Eqs. (\ref{g9}), (\ref{g10}) and of the expressions in  Eqs. 
(\ref{ch59})-(\ref{ch60-63}). 

\subsection{Power series of polar fourcomplex variables}

A polar fourcomplex series is an infinite sum of the form
\begin{equation}
a_0+a_1+a_2+\cdots+a_l+\cdots , 
\label{ch64}
\end{equation}
where the coefficients $a_l$ are polar fourcomplex numbers. The convergence of
the series (\ref{ch64}) can be defined in terms of the convergence of its 4
real components. The convergence of a polar fourcomplex series can however be
studied using polar fourcomplex variables. The main criterion for absolute
convergence remains the comparison theorem, but this requires a number of
inequalities which will be discussed further.

The modulus of a polar fourcomplex number $u=x+\alpha y+\beta z+\gamma t$ can
be defined as
\begin{equation}
|u|=(x^2+y^2+z^2+t^2)^{1/2} ,
\label{ch65}
\end{equation}
so that according to Eq. (\ref{ch12}) $d=|u|$. Since $|x|\leq |u|, |y|\leq |u|,
|z|\leq |u|, |t|\leq |u|$, a property of 
absolute convergence established via a comparison theorem based on the modulus
of the series (\ref{ch64}) will ensure the absolute convergence of each real
component of that series.

The modulus of the sum $u_1+u_2$ of the polar fourcomplex numbers $u_1, u_2$
fulfils the inequality
\begin{equation}
||u_1|-|u_2||\leq |u_1+u_2|\leq |u_1|+|u_2| .
\label{ch66}
\end{equation}
For the product the relation is 
\begin{equation}
|u_1u_2|\leq 2|u_1||u_2| ,
\label{ch67}
\end{equation}
which replaces the relation of equality extant for regular complex numbers.
The equality in Eq. (\ref{ch67}) takes place for 
$x_1=y_1=z_1=t_1, x_2=y_2=z_2=t_2$, or
$x_1=-y_1=z_1=-t_1, x_2=-y_2=z_2=-t_2$.
In particular
\begin{equation}
|u^2|\leq 2(x^2+y^2+z^2+t^2) .
\label{ch68}
\end{equation}
The inequality in Eq. (\ref{ch67}) implies that
\begin{equation}
|u^l|\leq 2^{l-1}|u|^l .
\label{ch69}
\end{equation}
From Eqs. (\ref{ch67}) and (\ref{ch69}) it results that
\begin{equation}
|au^l|\leq 2^l |a| |u|^l .
\label{ch70}
\end{equation}

A power series of the polar fourcomplex variable $u$ is a series of the form
\begin{equation}
a_0+a_1 u + a_2 u^2+\cdots +a_l u^l+\cdots .
\label{ch71}
\end{equation}
Since
\begin{equation}
\left|\sum_{l=0}^\infty a_l u^l\right| \leq  \sum_{l=0}^\infty
2^l|a_l| |u|^l ,
\label{ch72}
\end{equation}
a sufficient condition for the absolute convergence of this series is that
\begin{equation}
\lim_{l\rightarrow \infty}\frac{2|a_{l+1}||u|}{|a_l|}<1 .
\label{ch73}
\end{equation}
Thus the series is absolutely convergent for 
\begin{equation}
|u|<c_0,
\label{ch74}
\end{equation}
where 
\begin{equation}
c_0=\lim_{l\rightarrow\infty} \frac{|a_l|}{2|a_{l+1}|} .
\label{ch75}
\end{equation}

The convergence of the series (\ref{ch71}) can be also studied with the aid of
the formula (\ref{ch46b}) which, for integer values of $l$, is valid for any
$x, y, z, t$. If $a_l=a_{lx}+\alpha a_{ly}+\beta a_{lz}+\gamma a_{lt}$, and
\begin{eqnarray}
\begin{array}{l}
A_{l+}=a_{lx}+a_{ly}+a_{lz}+a_{lt}, \\
A_{l-}=a_{lx}-a_{ly}+a_{lz}-a_{lt},\\
A_{l1}=a_{lx}-a_{lz},\\
\tilde A_{l1}=a_{ly}-a_{lt},
\end{array}
\label{ch76}
\end{eqnarray}
it can be shown with the aid of relations (\ref{ch22})-(\ref{ch23b}) 
and (\ref{ch46b}) that the expression of the series (\ref{ch71}) is
\begin{equation}
\sum_{l=0}^\infty \left[A_{l+}  v_+^l e_++
A_{l-}v_-e_-+
\left(e_1A_{l1}+\tilde e_1\tilde A_{l1}\right)
\left(e_1v_1+\tilde e_1\tilde v_1\right)^l\right],
\label{ch77-78}
\end{equation}
where the quantities $v_+,  v_-$
have been defined in Eq. (\ref{ch8a}), and the quantities $v_1,\tilde v_1$ have
been defined in Eq. (\ref{ch16a}).

The sufficient conditions for the absolute convergence of the series 
in Eq. (\ref{ch77-78}) are that
\begin{equation}
\lim_{l\rightarrow \infty}\frac{|A_{l+1,+}||v_+|}{|A_{l+}|}<1,
\lim_{l\rightarrow \infty}\frac{|A_{l+1,-}||v_-|}{|A_{l-}|}<1,
\lim_{l\rightarrow \infty}\frac{A_{l+1}\mu_+}{A_l}<1,
\label{ch79}
\end{equation}
where the real and positive quantity $A_{l-}>0$ is given by
\begin{equation}
A_l^2=A_{l1}^2+\tilde A_{l1}^2.
\label{ch79a}
\end{equation}

Thus the series in Eq. (\ref{ch77-78}) is absolutely convergent for 
\begin{equation}
|x+y+z+t|<c_+,\:
|x-y+z-t|<c_-,\:
\mu_+<c_1
\label{ch80}
\end{equation}
where 
\begin{equation}
c_+=\lim_{l\rightarrow\infty} \frac{|A_{l+}|}{|A_{l+1,+}|} ,\:
c_-=\lim_{l\rightarrow\infty} \frac{|A_{l-}|}{|A_{l+1,-}|} ,\:
c_1=\lim_{l\rightarrow\infty} \frac{A_{l-}}{A_{l+1,-}} .
\label{ch81-87}
\end{equation}
The relations (\ref{ch80}) show that the region of convergence of the series
(\ref{ch77-78}) is a four-dimensional cylinder.
It can be shown that $c_0=(1/2)\;{\rm min}(c_+,c_-,c_1)$, 
where ${\rm min}$ designates the smallest of
the numbers in the argument of this function. Using the expression of $|u|$ in
Eq. (\ref{ch24d}), it can be seen that the spherical region of
convergence defined in Eqs. (\ref{ch74}), (\ref{ch75}) is included in the
cylindrical region of convergence defined in Eqs. (\ref{ch81-87}).

\subsection{Analytic functions of polar fourcomplex variables}

The fourcomplex function $f(u)$ of the fourcomplex variable $u$ has
been expressed in Eq. (\ref{g16}) in terms of 
the real functions $P(x,y,z,t),Q(x,y,z,t),R(x,y,z,t), S(x,y,z,t)$ of real
variables $x,y,z,t$.
The
relations between the partial derivatives of the functions $P, Q, R, S$ are
obtained by setting succesively in   
Eq. (\ref{g17}) $\Delta x\rightarrow 0, \Delta y=\Delta z=\Delta t=0$;
then $ \Delta y\rightarrow 0, \Delta x=\Delta z=\Delta t=0;$  
then $  \Delta z\rightarrow 0,\Delta x=\Delta y=\Delta t=0$; and finally
$ \Delta t\rightarrow 0,\Delta x=\Delta y=\Delta z=0 $. 
The relations are 
\begin{equation}
\frac{\partial P}{\partial x} = \frac{\partial Q}{\partial y} =
\frac{\partial R}{\partial z} = \frac{\partial S}{\partial t},
\label{ch95}
\end{equation}
\begin{equation}
\frac{\partial Q}{\partial x} = \frac{\partial R}{\partial y} =
\frac{\partial S}{\partial z} = \frac{\partial P}{\partial t},
\label{ch97}
\end{equation}
\begin{equation}
\frac{\partial R}{\partial x} = \frac{\partial S}{\partial y} =
\frac{\partial P}{\partial z} = \frac{\partial Q}{\partial t},
\label{ch96}
\end{equation}
\begin{equation}
\frac{\partial S}{\partial x} =\frac{\partial P}{\partial y} =
\frac{\partial Q}{\partial z} =\frac{\partial R}{\partial t}.
\label{ch98}
\end{equation}

The relations (\ref{ch95})-(\ref{ch98}) are analogous to the Riemann relations
for the real and imaginary components of a complex function. It can be shown
from Eqs. (\ref{ch95})-(\ref{ch98}) that the component $P$ is a solution
of the equations 
\begin{equation}
\frac{\partial^2 P}{\partial x^2}-\frac{\partial^2 P}{\partial z^2}=0,
\:\: 
\frac{\partial^2 P}{\partial y^2}-\frac{\partial^2 P}{\partial t^2}=0,
\:\:
\label{ch99}
\end{equation}
and the components $Q, R, S$ are solutions of similar equations.
As can be seen from Eqs. (\ref{ch99})-(\ref{ch99}), the components $P, Q, R, S$
of an analytic function of polar fourcomplex variable are harmonic with respect
to the pairs of variables $x,y$ and $ z,t$.  The component $P$ is also a
solution of the mixed-derivative equations
\begin{equation}
\frac{\partial^2 P}{\partial x^2}=\frac{\partial^2 P}{\partial y\partial t},
\:\: 
\frac{\partial^2 P}{\partial y^2}=\frac{\partial^2 P}{\partial x\partial z},
\:\:
\frac{\partial^2 P}{\partial z^2}=\frac{\partial^2 P}{\partial y\partial t},
\:\:
\frac{\partial^2 P}{\partial t^2}=\frac{\partial^2 P}{\partial x\partial z},
\:\:
\label{ch106b}
\end{equation}
and the components $Q, R, S$ are solutions of similar equations.
The component $P$ is also a solution of the mixed-derivative
equations 
\begin{equation}
\frac{\partial^2 P}{\partial x\partial y}=\frac{\partial^2 P}{\partial
z\partial t} ,
\:\: 
\frac{\partial^2 P}{\partial x\partial t}=\frac{\partial^2 P}{\partial
y\partial z} ,
\label{ch107}
\end{equation}
and the components $Q, R, S$ are solutions of similar equations.

\subsection{Integrals of functions of polar fourcomplex variables}

The singularities of polar fourcomplex functions arise from terms of the form
$1/(u-u_0)^m$, with $m>0$. Functions containing such terms are singular not
only at $u=u_0$, but also at all points of the two-dimensional hyperplanes
passing through $u_0$ and which are parallel to the nodal hyperplanes. 

The integral of a polar fourcomplex function between two points $A, B$ along a
path situated in a region free of singularities is independent of path, which
means that the integral of an analytic function along a loop situated in a
region free from singularities is zero,
\begin{equation}
\oint_\Gamma f(u) du = 0,
\label{ch111}
\end{equation}
where it is supposed that a surface $\Sigma$ spanning 
the closed loop $\Gamma$ is not intersected by any of
the hyperplanes associated with the
singularities of the function $f(u)$. Using the expression, Eq. (\ref{g16})
for $f(u)$ and the fact that $du=dx+\alpha  dy+\beta dz+\gamma dt$, the
explicit form of the integral in Eq. (\ref{ch111}) is
\begin{eqnarray}
\lefteqn{\oint _\Gamma f(u) du = \oint_\Gamma
[(Pdx+Sdy+Rdz+Qdt)+\alpha(Qdx+Pdy+Sdz+Rdt)\nonumber}\\
&&+\beta(Rdx+Qdy+Pdz+Sdt)+\gamma(Sdx+Rdy+Qdz+Pdt)] .
\label{ch112}
\end{eqnarray}
If the functions $P, Q, R, S$ are regular on a surface $\Sigma$
spanning the loop $\Gamma$,
the integral along the loop $\Gamma$ can be transformed with the aid of the
theorem of Stokes in an integral over the surface $\Sigma$ of terms of the form
$\partial P/\partial y - \partial S/\partial x,\:\:
\partial P/\partial z -  \partial R/\partial x, \:\:
\partial P/\partial t - \partial Q/\partial x, \:\:
\partial R/\partial y -  \partial S/\partial z, \:\:
\partial S/\partial t - \partial Q/\partial y, \:\:
\partial R/\partial t - \partial Q/\partial z$
and of similar terms arising
from the $\alpha, \beta$ and $\gamma$ components, 
which are equal to zero by Eqs. (\ref{ch95})-(\ref{ch98}), and this proves Eq.
(\ref{ch111}). 

The integral of the function $(u-u_0)^m$ on a closed loop $\Gamma$ is equal to
zero for $m$ a positive or negative integer not equal to -1,
\begin{equation}
\oint_\Gamma (u-u_0)^m du = 0, \:\: m \:\:{\rm integer},\: m\not=-1 .
\label{ch112b}
\end{equation}
This is due to the fact that $\int (u-u_0)^m du=(u-u_0)^{m+1}/(m+1), $ and to
the fact that the function $(u-u_0)^{m+1}$ is singlevalued for $m$ an integer.

The integral $\oint_\Gamma du/(u-u_0)$ can be calculated using the exponential
form (\ref{ch40c}),
\begin{equation}
u-u_0=
\rho\exp\left[\frac{1}{4}(\alpha+\beta+\gamma) \ln\frac{\sqrt{2}}{\tan\theta_+}
-\frac{1}{4}(\alpha-\beta+\gamma) \ln\frac{\sqrt{2}}{\tan\theta_-}
+\tilde e_1\phi\right] ,
\label{ch113}
\end{equation}
so that 
\begin{equation}
\frac{du}{u-u_0}=\frac{d\rho}{\rho}
+\frac{1}{4}(\alpha+\beta+\gamma) d\ln\frac{\sqrt{2}}{\tan\theta_+}
-\frac{1}{4}(\alpha-\beta+\gamma) d\ln\frac{\sqrt{2}}{\tan\theta_-}
+\tilde e_1d\phi.
\label{ch114}
\end{equation}
Since $\rho$, $\tan\theta_+$ and $\tan\theta_-$ are singlevalued variables,
it follows that 
$\oint_\Gamma d\rho/\rho =0, \oint_\Gamma d\ln\sqrt{2}/\tan\theta_+=0$,
and $\oint_\Gamma d\ln\sqrt{2}/\tan\theta_+=0$. On the other hand,
$\phi$ is a cyclic variables, so that it may give a contribution to
the integral around the closed loop $\Gamma$.
The result of the integrations will be given in the rotated system of
coordinates 
\begin{equation}
\xi=\frac{1}{\sqrt{2}}(x-z),\:
\upsilon=\frac{1}{\sqrt{2}}(y-t),\:
\tau=\frac{1}{2}(x+y+z+t),\:
\upsilon=\frac{1}{2}(x-y+z-t) .
\label{ch115}
\end{equation}
Thus, if $C_\parallel$ is a circle of radius $r$
parallel to the $\xi O\upsilon$ plane, and the
projection of the center of this circle on the $\xi O\upsilon$ plane
coincides with the projection of the point $u_0$ on this plane, the points
of the circle $C_\parallel$ are described according to Eqs.
(\ref{ch16a})-(\ref{ch17}) by the equations
\begin{equation}
\xi=\xi_0+r \cos\phi , \:
\upsilon=\upsilon_0+r \sin\phi , \:
\tau=\tau_0,\:\zeta=\zeta_0,
\label{ch115a}
\end{equation}
where $u_0=x_0+\alpha y_0+\beta 
z_0+\gamma t_0$, and $\xi_0, \upsilon_0, \tau_0, \zeta_0$ are calculated
from $x_0, y_0, z_0, t_0$ according to Eqs. (\ref{ch115}).

Then
\begin{equation}
\oint_{C_\parallel}\frac{du}{u-u_0}=2\pi\tilde e_1. 
\label{ch116}
\end{equation}
The expression of $\oint_\Gamma du/(u-u_0)$ can be written 
with the aid of the functional int($M,C$) defined in Eq. (\ref{118}) as
\begin{equation}
\oint_\Gamma\frac{du}{u-u_0}=
2\pi\tilde e_1 \;{\rm int}(u_{0\xi\upsilon},\Gamma_{\xi\upsilon}) ,
\label{ch119}
\end{equation}
where $u_{0\xi\upsilon}$ and $\Gamma_{\xi\upsilon}$ are respectively the
projections of the point $u_0$ and of 
the loop $\Gamma$ on the plane $\xi \upsilon$.

If $f(u)$ is an analytic polar fourcomplex function which can be expanded in a
series as written in Eq. (\ref{g14}), and the expansion holds on the curve
$\Gamma$ and on a surface spanning $\Gamma$, then from Eqs. (\ref{ch112b}) and
(\ref{ch119}) it follows that
\begin{equation}
\oint_\Gamma \frac{f(u)du}{u-u_0}=
2\pi\tilde e_1 \;{\rm int}(u_{0\xi\upsilon},\Gamma_{\xi\upsilon})f(u_0) ,
\label{ch120}
\end{equation}
where $\Gamma_{\xi\upsilon}$ is the projection of 
the curve $\Gamma$ on the plane $\xi \upsilon$,
as shown in Fig. 5.
Substituting in the right-hand side of 
Eq. (\ref{ch120}) the expression of $f(u)$ in terms of the real 
components $P, Q, R, S$, Eq. (\ref{g16}), yields
\begin{equation}
\oint_\Gamma \frac{f(u)du}{u-u_0}=\pi 
\left[(\beta-1)(Q-S)+(\alpha-\gamma)(P-R)\right] 
{\rm int}\left(u_{0\xi\upsilon},\Gamma_{\xi\upsilon}\right),
\label{ch121}
\end{equation}
where $P, Q, R, S$ are the values of the components of $f$ at $u=u_0$.
If the integral is written as
\begin{equation}
\oint_\Gamma \frac{f(u)du}{u-u_0}=
I+\alpha I_\alpha+\beta I_\beta+\gamma I_\gamma,
\label{ch121b}
\end{equation}
it results from Eq. (\ref{ch121}) that
\begin{equation}
I+ I_\alpha+ I_\beta+ I_\gamma=0.
\label{ch121c}
\end{equation}

If $f(u)$ can be expanded as written in Eq. (\ref{g14}) on 
$\Gamma$ and on a surface spanning $\Gamma$, then from Eqs. (\ref{ch112b}) and
(\ref{ch119}) it also results that
\begin{equation}
\oint_\Gamma \frac{f(u)du}{(u-u_0)^{m+1}}=
2\frac{\pi}{m!}\tilde e_1 \;{\rm int}(u_{0\xi\upsilon},\Gamma_{\xi\upsilon}) 
\;f^{(m)}(u_0) ,
\label{ch122}
\end{equation}
where the fact has been used  that the derivative $f^{(m)}(u_0)$ of order $m$
of $f(u)$ at $u=u_0$ is related to the expansion coefficient in Eq. (\ref{g14})
according to Eq. (\ref{g15}).

If a function $f(u)$ is expanded in positive and negative powers of $u-u_j$,
where $u_j$ are polar fourcomplex constants, $j$ being an index, the integral
of $f$ on a closed loop $\Gamma$ is determined by the terms in the expansion of
$f$ which are of the form $a_j/(u-u_j)$,
\begin{equation}
f(u)=\cdots+\sum_j\frac{a_j}{u-u_j}+\cdots
\label{ch123}
\end{equation}
Then the integral of $f$ on a closed loop $\Gamma$ is
\begin{equation}
\oint_\Gamma f(u) du = 
2\pi\tilde e_1 \sum_j{\rm int}(u_{j\xi\upsilon},\Gamma_{\xi\upsilon})a_j .
\label{ch124}
\end{equation}

\subsection{Factorization of polar fourcomplex polynomials}

A polynomial of degree $m$ of the polar fourcomplex variable 
$u=x+\alpha y+\beta z+\gamma t$ has the form
\begin{equation}
P_m(u)=u^m+a_1 u^{m-1}+\cdots+a_{m-1} u +a_m ,
\label{ch125}
\end{equation}
where the constants are in general polar fourcomplex numbers.

If $a_m=a_{mx}+\alpha a_{my}+\beta a_{mz}+\gamma a_{mt}$, and with the
notations of Eqs. (\ref{ch8a}) and (\ref{ch76}) applied for $0, 1, \cdots, m$ ,
the polynomial $P_m(u)$ can be written as
\begin{eqnarray}
\lefteqn{P_m= \left[v_+^m 
+A_1 v_+^{m-1}+\cdots+A_{m-1} v_++ A_m \right] e_+
+\left[v_-^m 
+A_1^{\prime\prime} v_-^{m-1} +\cdots+A_{m-1}^{\prime\prime} 
v_-+ A_m^{\prime\prime} \right]e_-\nonumber}\\
 & &+\left[
\left(e_1v_1+\tilde e_1\tilde v_1\right)^{m}
+\sum_{l=1}^m\left(e_1A_{l1}+\tilde e_1\tilde A_{l1}\right)
\left(e_1v_1+\tilde e_1\tilde v_1\right)^{m-l}\right],\nonumber\\
&&
\label{ch126}
\end{eqnarray}
where the constants $A_{l+}, A_{l-}, A_{l1}, \tilde A_{l1}$ 
are real numbers.
The polynomial of degree $m$ in $(e_1v_1+\tilde e_1\tilde v_1)$
can always be written as a product of linear factors of the form
$[e_1(v_1-v_{1l})+\tilde e_1 (\tilde v_1-\tilde v_{1l})]$, where the
constants $v_{1l}, \tilde v_{1l}$ are real.
The two polynomials of degree $m$ with real coefficients in Eq. (\ref{ch126})
which are multiplied by $e_+$ and $e_-$ can be written as a product
of linear or quadratic factors with real coefficients, or as a product of
linear factors which, if imaginary, appear always in complex conjugate pairs.
Using the latter form for the simplicity of notations, the polynomial $P_m$
can be written as
\begin{equation}
P_m=\prod_{l=1}^m (v_+-s_{l+})e_+
+\prod_{l=1}^m (v_--s_{l-})e_-
+\prod_{l=1}^m
\left[e_1(v_1-v_{1l})+
\tilde e_1(\tilde v_1-\tilde v_{1l})\right],
\label{ch127}
\end{equation}
where the quantities $s_{l+}$ appear always in complex conjugate pairs, and the
same is true for the quantities $s_{l-}$.
Due to the properties in Eqs. (\ref{ch22})-(\ref{ch23b}),
the polynomial $P_m(u)$ can be written as a product of factors of
the form  
\begin{equation}
P_m(u)=\prod_{l=1}^m \left[(v_+-s_{l+})e_+
+\left(v_--s_{l-}\right)e_-
+\left(e_1(v_1-v_{1l})+
\tilde e_1(\tilde v_1-\tilde v_{1l})\right)
\right].
\label{ch128}
\end{equation}
These relations can be written with the aid of Eq. (\ref{ch24}) as
\begin{eqnarray}
P_m(u)=\prod_{p=1}^m (u-u_p) ,
\label{ch128c}
\end{eqnarray}
where
\begin{eqnarray}
u_p=s_{p+} e_+
+s_{p-}e_-
+e_1v_{1p}+
\tilde e_1\tilde v_{1p}, p=1,...,m.
\label{ch128d}
\end{eqnarray}
The roots $s_{p+},
s_{p-}, v_{1p}e_1+\tilde v_{1p}
\tilde e_1$ of the corresponding polynomials in Eq.
(\ref{ch127}) may be ordered arbitrarily.
This means that Eq. (\ref{ch128d}) gives sets of $m$ roots
$u_1,...,u_m$ of the polynomial $P_m(u)$, 
corresponding to the various ways in which the roots
$s_{p+}, s_{p-},v_{1p} e_1+\tilde v_{1p}
\tilde e_1$ are ordered according to $p$ in each
group. Thus, while the hypercomplex components in Eq. (\ref{ch126}) taken
separately have unique factorizations, the polynomial $P_m(u)$ can be written
in many different ways as a product of linear factors. 
The result of the polar fourcomplex integration, Eq. (\ref{ch124}), is however
unique.

If $P(u)=u^2-1$, 
the factorization in Eq. (\ref{ch128c}) is $u^2-1=(u-u_1)(u-u_2)$, where 
$u_1=\pm e_+\pm e_-\pm e_1, u_2=-u_1$, so that
there are 4 distinct factorizations of $u^2-1$,
\begin{eqnarray}
\begin{array}{l}
u^2-1=(u-1)(u+1),\\
u^2-1=(u-\beta)(u+\beta),\\
u^2-1=\left(u-\frac{1+\alpha-\beta+\gamma}{2}\right)
\left(u+\frac{1+\alpha-\beta+\gamma}{2}\right),\\
u^2-1=\left(u-\frac{-1+\alpha+\beta+\gamma}{2}\right)
\left(u+\frac{-1+\alpha+\beta+\gamma}{2}\right).
\end{array}
\label{ch129}
\end{eqnarray}

It can be checked that 
$\left\{\pm e_+\pm e_-\pm e_1\right\}^2= 
e_++e_-+e_1=1$.

\subsection{Representation of polar fourcomplex numbers 
by irreducible matrices}

If $T$ is the unitary matrix,
\begin{equation}
T =\left(
\begin{array}{cccc}
\frac{1}{{2}}&\frac{1}{{2}}    &\frac{1}{{2}}    &\frac{1}{{2}}    \\
\frac{1}{{2}}  &-\frac{1}{{2}} &\frac{1}{{2}}    &-\frac{1}{{2}}   \\
\frac{1}{\sqrt{2}}& 0          & -\frac{1}{\sqrt{2}}  & 0               \\
0                 & \frac{1}{\sqrt{2}} & 0       & -\frac{1}{\sqrt{2}}  \\
\end{array}
\right),
\label{ch129x}
\end{equation}
it can be shown 
that the matrix $T U T^{-1}$ has the form 
\begin{equation}
T U T^{-1}=\left(
\begin{array}{ccc}
x+y+z+t &    0    &  0    \\
0       & x-y+z-t &  0    \\
0       &    0    &  V_1  \\
\end{array}
\right),
\label{ch129y}
\end{equation}
where $U$ is the matrix in Eq. (\ref{ch23}) used to represent the polar fourcomplex
number $u$. In Eq. (\ref{ch129y}), $V_1$ is the matrix
\begin{equation}
V_1=\left(
\begin{array}{cc}
x-z    &   y-t   \\
-y+t   &   x-z   \\
\end{array}\right).
\label{ch130x}
\end{equation}
The relations between the variables $x+y+z+t, x-y+z-t, x-z, y-t$ for the
multiplication of polar fourcomplex numbers have been written in Eqs.
(\ref{ch15}),(\ref{ch16}), (\ref{ch19a}), (\ref{ch19b}). The matrix $T U
T^{-1}$ provides an irreducible representation
\cite{4} of the polar fourcomplex number $u$ in terms of matrices with real
coefficients.

\section{Conclusions}

In the case of the circular fourcomplex numbers, the operations of addition and
multiplication have a simple geometric interpretation based on the amplitudes
$\rho$, the azimuthal angles $\phi, \chi$, and the planar angles $\psi$.
Exponential form exist for the circular fourcomplex numbers, involving the
variables $\upsilon, \phi, \chi$ and $\psi$. The circular fourcomplex functions
defined by series of powers are analytic, and the partial derivatives of the
components of the circular fourcomplex functions are closely related. The
integrals of circular fourcomplex functions are independent of path in regions
where the functions are regular. The fact that the exponential forms of the
circular fourcomplex numbers depend on the cyclic variables $\phi, \chi$ leads
to the concept of pole and residue for integrals on closed paths. The
polynomials of circular fourcomplex variables can always be written as products
of linear factors.

In the case of the hyperbolic fourcomplex numbers, an exponential form exists
for the hyperbolic fourcomplex numbers, involving the amplitude $\mu$ and the
arguments $y_1, z_1, t_1$.  The hyperbolic fourcomplex functions defined by
series of powers are analytic, and the partial derivatives of the components of
the hyperbolic fourcomplex functions are closely related. The integrals of
hyperbolic fourcomplex functions are independent of path in regions where the
functions are regular.  The polynomials of tricomplex variables can be written
as products of linear or quadratic factors.

In the case of the planar fourcomplex numbers,
the operations of addition and multiplication of the planar fourcomplex numbers
introduced in this 
work have a geometric interpretation based on the amplitude $\rho$, the
azimuthal angles $\phi, \chi$, and the planar angle $\psi$. An exponential form
exists for the planar fourcomplex numbers, involving the variables $\rho, \phi,
\chi$ and $\psi$. The planar fourcomplex functions defined by series of powers
are analytic, and the partial derivatives of the components of the planar
fourcomplex functions are closely related. The integrals of planar fourcomplex
functions are independent of path in regions where the functions are regular.
The fact that the exponential form of the planar fourcomplex numbers depends on
the cyclic variables $\phi, \chi$ leads to the concept of pole and residue for
integrals on closed paths. The polynomials of planar fourcomplex variables can
always be written as products of linear or quadratic factors.

In the case of the polar fourcomplex numbers, the operations of addition and
multiplication have a geometric interpretation based on the amplitude $\sigma$,
the polar angles $\theta_+,\theta_-$ and the azimuthal angle $\phi$. An
exponential form exists for the polar fourcomplex numbers, involving the
variables $\sigma, \phi$ and the angles $\theta_+, \theta_-$. The polar
fourcomplex functions defined by series of powers are analytic, and the partial
derivatives of the components of the polar fourcomplex functions are closely
related. The integrals of polar fourcomplex functions are independent of path
in regions where the functions are regular. The fact that the exponential form
of the polar fourcomplex numbers depends on the cyclic variable $\phi$ leads to
the concept of pole and residue for integrals on closed paths. The polynomials
of polar fourcomplex variables can be written as products of linear or
quadratic factors.

\newpage

FIGURE CAPTIONS\\

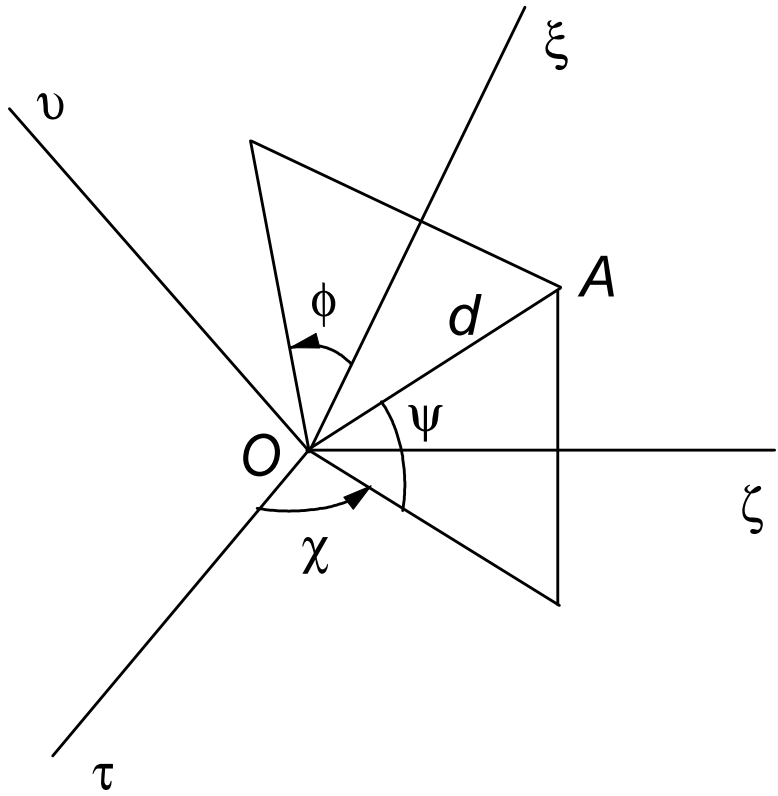
Fig. 1. Azimuthal angles $\phi, \chi$ and planar angle $\psi$ of the
fourcomplex number $x+\alpha y +\beta z+\gamma t$, represented by the point
$A$, situated at a distance $d$ from the origin $O$.\\

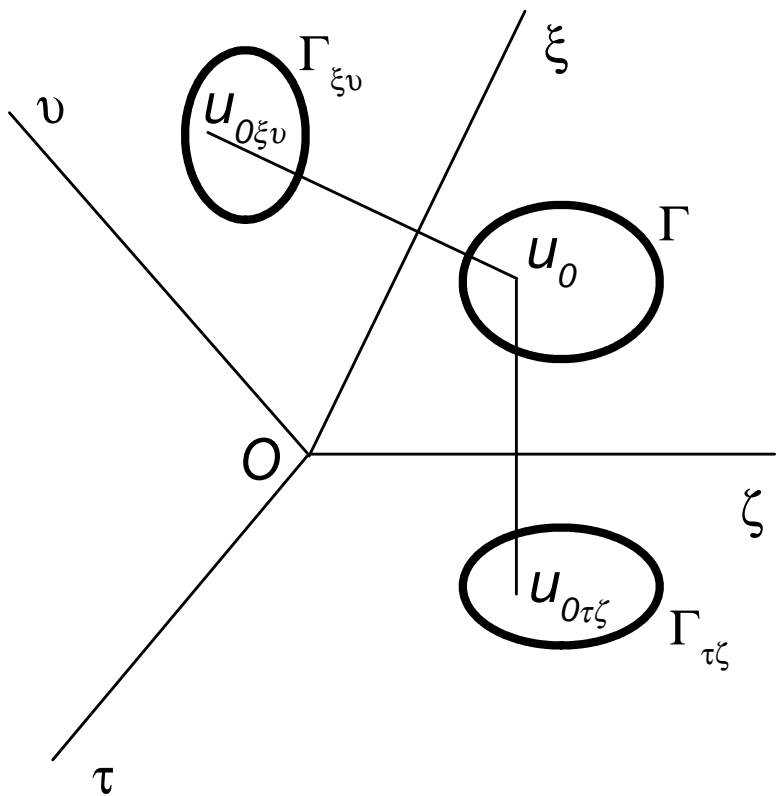
Fig. 2. Integration path $\Gamma$ and the pole $u_0$, and their projections
$\Gamma_{\xi\upsilon}, \Gamma_{\tau\zeta}$ and $u_{0\xi\upsilon},
u_{0\tau\zeta}$ on the planes $\xi \upsilon$ and respectively $\tau \zeta$.\\

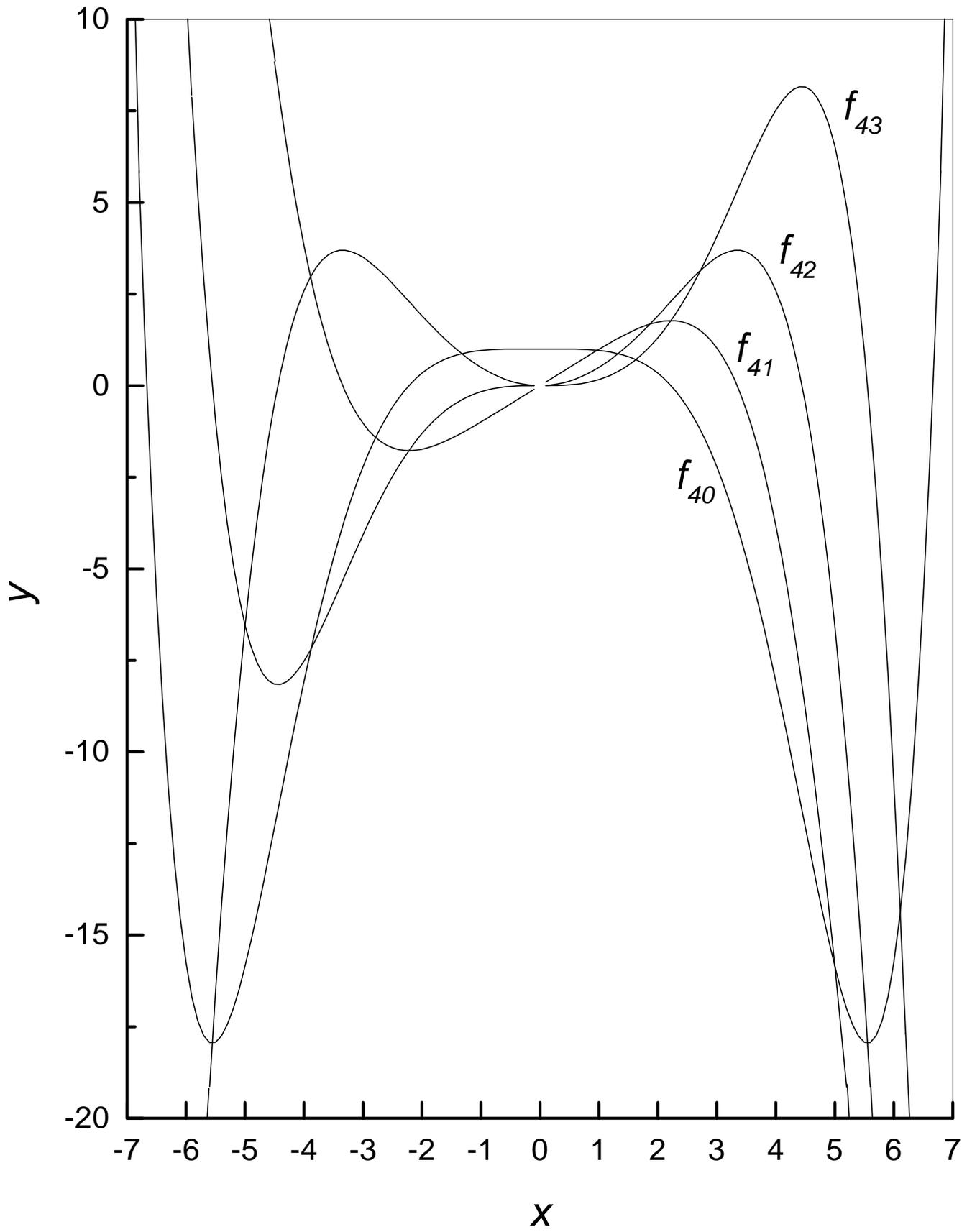
Fig. 3. The planar fourdimensional cosexponential functions $f_{40}, f_{41},
f_{42}, f_{43}$.\\

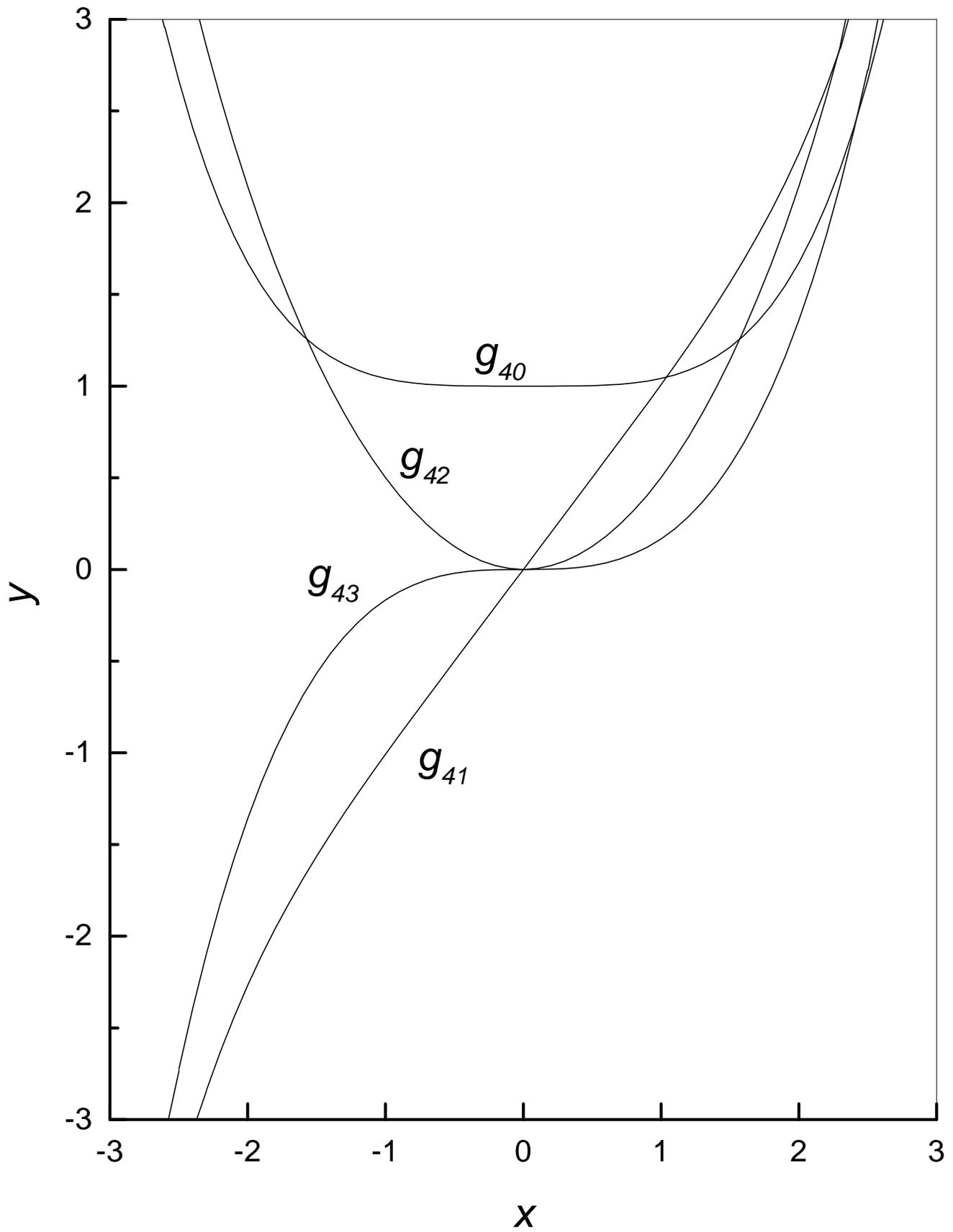
Fig. 4. The polar fourdimensional cosexponential functions $g_{40}, g_{41},
g_{42}, g_{43}$.\\

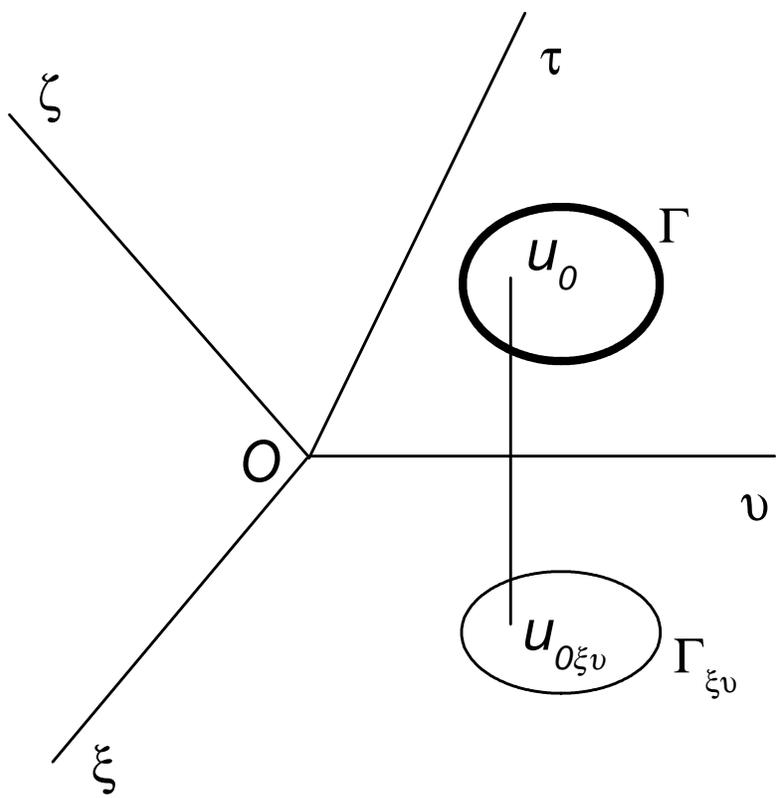
Fig. 5. Integration path $\Gamma$ and the pole $u_0$, and their projections
$\Gamma_{\xi\upsilon}$ and $u_{0\xi\upsilon}$ on the plane $\xi \upsilon$.\\

\end{document}